\documentclass[review]{elsarticle}
\usepackage[utf8]{inputenc}

\usepackage{grffile}
\usepackage{lineno}
\usepackage{times}
\usepackage{calc}
\usepackage{graphics}
\usepackage{epsfig}
\usepackage{amssymb}
\usepackage{amsmath}
\usepackage{flafter}
\usepackage{float}
\usepackage{latexsym}
\usepackage{psfrag}
\usepackage{dsfont}
\usepackage{rotating}
\usepackage{longtable}
\usepackage{cleveref}
\usepackage{color}
\usepackage{longtable,lscape}
\usepackage{multirow}
\usepackage{tabularx}
\usepackage{booktabs}
\usepackage{subcaption}
\usepackage{array}
\newcolumntype{M}[1]{>{\centering\arraybackslash}m{#1}}

\usepackage{pgfplots}
\usepackage{tikzscale}
\usepackage{pgfkeys}
\usepackage{pgfplotstable}
\usetikzlibrary{matrix}
\usepgfplotslibrary{groupplots}
\pgfplotsset{compat=newest}
\pgfplotsset{}
\pgfplotsset{
	tick label style={font=\footnotesize},
	label style={font=\footnotesize},
	legend style={font=\footnotesize}
}
\pgfplotsset{major grid style={loosely dotted, thin, gray}}
\pgfkeys{
	/includetikzgraphic/.is family, /includetikzgraphic,
	default/.style = {
		width = 0.85\textwidth,
		ratio = 0.85},
	width/.estore in = \itgWidth,
	ratio/.estore in = \itgRatio
}
\newlength\figurewidth
\newlength\figureheight
\newcommand{\includetikzgraphic}[2][]{
	\pgfkeys{/includetikzgraphic, default, #1}
	\setlength\figurewidth{\itgWidth}
	\setlength\figureheight{\itgRatio\figurewidth}
	\includegraphics{#2}
}
\modulolinenumbers[5]

\journal{Journal of \LaTeX\ Templates}

\newtheorem{remark}{Remark}
\renewcommand{\vec}[1]{\mathbf{#1}}

\newcommand{\eps}{\varepsilon}

\begin{document}

\begin{frontmatter}

\title{A Neural Network based Shock Detection and Localization Approach for Discontinuous Galerkin Methods}

\author{Andrea D. Beck$^a$\fnref{myfootnote}}
\author{Jonas Zeifang$^a$}
\author{Anna Schwarz$^a$}
\author{David G. Flad$^b$}
\address{$^a$ Institute of Aerodynamics and Gas Dynamics, University of Stuttgart, Stuttgart, Germany}
\address{$^b$ {NASA Ames Research Center, Moffett Field, CA, USA}}
\fntext[myfootnote]{Corresponding Author; $\;$A. Beck and J. Zeifang agree to share first authorship}





\begin{abstract}
The stable and accurate approximation of discontinuities such as shocks on a finite computational mesh is a challenging task. Detection of shocks or strong discontinuities in the flow solution is typically achieved through a priori troubled cell indicators, which guide the subsequent action of an appropriate shock capturing mechanism. Arriving at a stable and accurate solution often requires empirically based parameter tuning and adjustments of the indicator settings to the discretization and solution at hand.
	In this work, we propose to separate the task of shock detection and shock capturing more strongly and aim to develop a shock indicator that is robust, accurate, requires minimal user input and is suitable for high order element-based methods like discontinuous Galerkin and flux reconstruction methods. The novel indicator is learned from analytical data through a supervised learning strategy; its input is given by the high order solution field, its output is an element-local map of the shock position. We use state of the art methods from edge detection in image analysis based on deep convolutional multiscale networks and deep supervision to train the indicators. The resulting networks are then used as black box indicators, showing their robustness and accuracy on well established canonical testcases. All simulations are run ab initio using the developed indicators, showing that they provide also stability during the strongly transient phases. In particular for high order schemes with large cells and considerable inner-cell resolution capabilities,  we demonstrate how the additional accurate prediction of the position of the shock front can be exploited to guide inner-element shock capturing strategies.

\end{abstract}


\end{frontmatter}

\section{\label{sec:intro}Introduction}
\subsection{\label{subsec:sota}Problem Definition and State of the Art}
The occurrence of strong compression waves in compressible fluid flow is a well-known and characteristic phenomenon~\cite{courant1999supersonic}. Mathematically, these shocks are an expression of the non-linearity of the underlying governing equations (e.g. the Euler equations), which can produce discontinuous solution even for initially smooth data~\cite{riemann1860fortpflanzung}. This possibility for discontinuities often clashes with the discretization strategy for the numerical approximation of these flows, where a (at least locally) smooth solution representation is sought. Discretization strategies with ansatz spaces that include discontinuities like low order finite volume (FV) schemes are naturally better equipped to handle shocks, leading to the dominance of FV and related schemes for compressible aerodynamics, but encounter the same mismatch for orders greater than one~\cite{leveque1998nonlinear}.\\
 Based on the rationale that the solution is predominantly smooth except for limited regions, high order discretizations are often favored and have gained considerable interest for practical applications. These schemes provide low approximation errors and excellent scale resolving capabilities, leading to efficient methods in particular for smooth problems~\cite{wang2013high}. In regions with discontinuities, these schemes then encounter spurious oscillations, often called Gibb's instability~\cite{gottlieb1997gibbs}, which can lead to a loss of stability and accuracy. Avoiding instabilities and treating these solutions in a numerically stable and accurate way is the task of the so-called shock capturing methods. A conceptually different approach, labeled shock-fitting, is significantly less widespread due to its complexities on unstructured grids~\cite{paciorri2009shock}.\\ Depending on the baseline discretization scheme chosen, a number of different shock capturing methods exist, for example, for finite volume discretizations, a TVD limiting strategy or an essentially non-oscillatory (W)ENO reconstruction are typical~\cite{pirozzoli2011numerical,shu1988efficient,shu2003high,leveque2002finite}. For the class of flux reconstruction (FR), spectral difference (SD) and discontinuous Galerkin (DG) schemes which are the focus of this work, the methods  typically fall into two categories: In the first category, once a shock is detected, a local diffusion operator of prescribed strength is added to regularize the solution. This strategy is often called the "artificial viscosity shock capturing" method~\cite{premasuthan2014computation,lv2016entropy,FEISTAUER20101612,BARTER20101810,persson2006,HARTMANN2002508}. In the second category, the discretization operator is switched \emph{locally} to a more or guaranteed robust one, for example, through h-refinement and/or p-coarsening or a reconstruction procedure. By appropriate limiting, this operator may then additionally  ensure the boundedness of the solution. We will refer to these methods rather broadly as "h/p-shock capturing" in the following. Examples of the combination of high order discretization schemes with these shock capturing methods can ,for example, be found in~\cite{10.2307/2008501,BURBEAU2001111,yang2009parameter,QIU2005642,dumbser2014posteriori,doi:10.1137/04061372X,BALSARA2007586,ZHONG2013397}. Before applying any of these strategies however, the occurrence and location of a shock must be detected. This, in itself, is a non-trivial task and is typically achieved through \emph {parameter-dependent indicator functions} which are based on physical considerations~\cite{jameson1981}, modal smoothness estimates~\cite{huerta2012,klockner2011viscous,persson2006} or image detection ideas~\cite{sheshadri2014shock}. So at least in theory, in this scenario, the shock detection task can be decoupled from the numerical solution stabilization mechanisms. However, in practical applications, this is most often not the case. Instead, the parameters of the indicator, which informs the local shock capturing action of the scheme, are tuned to provide a stable simulation result. Thus, the indicator not primarily signals a discontinuity but is tuned to \emph{predict} an underresolved state that will likely lead to an invalid solution of the discretized PDE in an a priori fashion. In this way, the "correct" parameter settings for each indicator become dependent not only on the underlying flow field but also heavily on the discretization and its robustness. This dilemma is reflected in the fact that indicators used in this setting are often labeled \emph{troubled cell} indicators, where the potential trouble is not the existence of a sharp gradient or discontinuity in the solution but the \emph{failure of the current discretization to deal with it robustly}. In our presented approach here, we return to the original separation of tasks. \\
 A posteriori methods are less frequently applied, a recent example can be found in~\cite{dumbser2014posteriori}. Here, the occurrence of an unstable solution in troubled cells is detected after the corresponding time step. The solution is then rolled back to a previous, stable state, and a more robust discretization, typically a lower order scheme, is chosen to recompute the solution in a recursive manner until a stable solution is achieved. Within this method, the number of roll-back steps and the shock detection criteria are tuneable parameters. Their optimal values are not known a priori and need to be guided by the user.\\
 As can be surmised from the discussion so far, the detection (or prediction) of troubled cells from the solution of the discretized PDE usually contains at least some level of empiricism and parameter tuning. This somewhat unsatisfactory state results from the complexity of the problem due to the inherent non-linearities in the equation itself and its discretization. Given the rise of the interest in data-driven methods in fluid mechanics with applications to a range of non-linear problems (e.g. subspace model design, active control~\cite{beck2019deep,rabault2019artificial}), it is not surprising that these methods have also been investigated for the problem at hand. Before presenting our approach, we thus briefly survey the existing literature on shock localization and solution stabilization with data-based approaches and machine learning techniques in the next section. \\

\subsection{\label{subsec:sotaml} Shock Detection and Capturing with Data-based Methods}

In this section, we present a brief summary of shock detection and localization methods with a focus on data-based approaches. The initial concept of shock detection based on image analysis as well as application examples of a classical Sobel filtering can be found in~\cite{rusanov1973processing,vorozhtsov1987shock}. In~\cite{liou1995image}, the authors propose a staged approach, where step edges are located through an approximation of second order derivatives of the pressure field. In the regions of interest, a shock function is then evaluated, which incorporates physical considerations like Mach number and pressure gradients. Contours of density gradient or normal Mach number was also used successfully in~\cite{pagendarm1993algorithm,lovely1999shock}; a collection of results of image-based methods can be found in~\cite{WU2013501}. A recent example of shock detection based on modern neural network architectures is presented in~\cite{10.1007/978-3-319-61358-1_16}, where the network was trained to predict the Schlieren value from local velocity strain data; an improved method with higher detection accuracy was shown in~\cite{LIU20191}. While the previous examples only consider shock detection, there have also been some efforts to combine data-based methods for combined shock detection / shock capturing, i.e. as troubled cell indicators. In~\cite{hanveiga:hal-01856358}, the authors classify a troubled cell as one in which the unlimited solution deviates by a user-specified  metric from one limited by the methods proposed in~\cite{krivodonova2007limiters}. They then train a multilayer perceptron (MLP) as a classifier to predict the occurrence of troubled cells for 1D advection and Euler equations. While their network is capable of essentially mimicking the underlying limiter from~\cite{krivodonova2007limiters}, the number of neurons and the associated cost of the network  are substantial. In~\cite{ray2018,ray2019}, the authors develop an MLP based classifier to distinguish smooth and non-smooth regions of the flow and show its robustness by applying it to a range of test problems. An extension of the classifier towards the prediction of an artifical viscosity is presented in~\cite{yu2018data}. Recently, Morgan et al. also proposed an MLP based shock detector for a high order finite element scheme~\cite{morgan2020machine}.

\subsection{\label{subsec:motivation}Motivation and Outline}
In the approach considered here, we return to the original separation of tasks for shock capturing: shock detection and solution stabilization are considered two different, \emph{independent} entities, and we make the choice of treating them as such. For the second task, we rely on a guaranteed stable approximation and employ a robust finite volume formulation with TVD limiters on inner-element sub-cells~\cite{sonntag2017efficient,krais2019flexi}, but note that other solution stabilization methods, e.g. artificial viscosity based ones, are equally possible. For the purpose of this work, we will treat this shock capturing method as a black box with fixed settings and not attempt any parameter tuning. We will instead focus on the first task: Reliably and accurately detecting the occurrence of shocks or steep gradients in an a priori manner based on the current solution alone without the need for parameter adjustment. In order to do this, we propose a novel strategy for developing shock location sensors based on ideas from machine learning, in particular edge detection. We follow the \emph{supervised learning} paradigm, where a (computationally expensive) offline training process based on labeled data produces a (cheap) approximation to the underlying mapping between input and output data. In the context of edge detection, we interpret the flow solution data within one grid element as the image and a possible shock or strong gradient as an edge to be detected and localized within this image. Thus, the input data into our method consists of the conservative solution variables in a given mesh cell, the output is whether a shock is present or not and where it is located. The mapping from input to output data itself is expressed as a convolutional neural network (CNN). The approach we propose is particularly attractive for high order schemes like DG, FR or SD schemes, as it takes advantage of the full high order information contained in a single grid element. The presented novel indicators  offer the advantage of being parameter-free and do not require tuning by the user. We will demonstrate their applicability to a range of canonical flow problems. In addition, our results will show that they are highly insensitive towards changes in grid resolution. Finally, we demonstrate that with a slight modification of the network training, the resulting indicators can now be used not only to flag cells that contain shocks but also to locate the shock precisely within a cell - this local information can then inform the subsequent shock capturing strategy and help target it precisely to where it is needed. Again, this information is particularly beneficial for high order schemes with typically rather large grid cells, where a localized shock capturing on the scale of the smallest resolved length scale (roughly estimated in 1D as $\Delta x/p$, where $\Delta x$ is the length of the considered cell and $p$ is the degree of the polynomial ansatz) is desirable. While we demonstrate the potential benefits of this approach in the context of the h/p-shock capturing with the chosen local finite volume based solution stabilization method, it can readily be transferred to other approaches and could for example inform a local artificial viscosity approach or guide element-local r-adaptation.\\

For the remainder of this paper, we refer to the task of detection whether a shock exists in a grid cell or not as "shock detection" and the corresponding developed method as ANNSI (artificial neural network based shock indicator); the task of locating the exact position of the shock front within an element is called "shock localization" and our proposed method ANNSL (artificial neural network based shock localization).\\

We will proceed as follows: In Sec.~\ref{sec:NumericalMethod}, we introduce the numerical method employed in this work and the present shock capturing method used later in a black box manner. We also review two classical troubled cell indicators which serve as reference for the evaluation of the ANNSI results. In Sec.~\ref{sec:data}, a short summary on neural networks is given, we define the training data and the network architectures used for the ANNSI and ANNSL indicators, and we close the section by reporting the results of the training process. In the following section (Sec.~\ref{sec:NNshock}), we apply the novel ANNSI indicator to well known test problems and compare the results with two classical troubled cell indicators. Once we have established the accuracy and robustness of the ANNSI indicator, we focus on its in-build extension to localize shocks within the grid elements. In Sec.~\ref{sec:NNshock_local}, we therefore apply the ANNSL indicator to several test problems and show how the information about the inner cell shock position can be exploited. Finally, in Sec.~\ref{sec:conout}, we conclude and give an outlook.\\

\section{Numerical Method}\label{sec:NumericalMethod}
\subsection{Governing Equations}
In this work, we consider the usual compressible Euler equations of fluid dynamics which are given by
\begin{align}\label{eq:euler}
	\partial_t\underbrace{\begin{pmatrix}\rho\\\rho\vec u\\ E\end{pmatrix}}_{=:\vec w} + \nabla\cdot\underbrace{\begin{pmatrix}\rho \vec u\\ \rho \vec u \otimes \vec u + p~\vec{I}\\ \vec u \left(E + p\right)
	\end{pmatrix}}_{=:\vec{F}(\vec w)} = 0,
\end{align}
where $\rho$ denotes density, $\vec u$ velocity, $E$ energy density and the pressure $p$ is computed via the equation of state for a perfect gas
\begin{align*}
	p(\rho,\rho\vec u,E) := (\gamma - 1) \left(E - \frac{1}{2} \rho \|\vec u\|^2_2 \right),
\end{align*}
with $\gamma=1.4$ being the isentropic expansion coefficient.
\subsection{High Order Discontinuous Galerkin Spectral Element Method}\label{subsec:dg}
We discretize the equations (Eq.~\eqref{eq:euler}) with a Discontinuous Galerkin Spectral Element Method (DGSEM)~\cite{kopriva2009,hindenlang2012explicit,krais2019flexi}, offering the possibility of choosing arbitrary high approximation order. DGSEM is based on the weak formulation for a piecewise smooth function $\vec{w}_h$
\begin{equation}\label{eq:dg:weak}
	\frac{\partial}{\partial t}\int\limits_{E} \vec{w}_h\,  \phi(\vec{x})\, d\vec{x}  +\oint\limits_{\partial E} \vec{{F}}^{*}_n  \phi(\vec{x})\,ds - \int\limits_{E}
	\vec{ F }(\vec{w}_h)\cdot\vec{\nabla}_x\,\phi(\vec{x})\, d\vec{x} =0,
\end{equation}
on every element $E$ for every polynomial test function $\phi(\vec{x})$ of degree $\mathcal{N}$. Note that all elements and thus Eq.~\eqref{eq:dg:weak} are mapped onto a reference space $[-1,1]^d$. DGSEM uses an inner-element solution representation based on the tensor product of nodal Lagrange basis functions and couples neighboring elements weakly through a numerical flux $\vec{F}^*_n$. Choosing the same $\mathcal{N}+1$ Legendre-Gauss points for quadrature and interpolation reduces the number of required operations per degree of freedom. With the restriction to tensor-product elements, the total number of degrees of freedom per grid cell results in $(\mathcal{N}+1)^d$. The final discrete $d$-dimensional spatial operator can be constructed through multiple one-dimensional operations. For the remainder of this text, we limit our discussion to 2D problems, i.e. $d=2$.
According to the method of lines, integration in time is then achieved by any suitable explicit or implicit method. For further details the reader is referred to~\cite{krais2019flexi},~\cite{beck2014high} and~\cite{hindenlang2012explicit}; the full code framework is available as open source software\footnote{ www.flexi-project.org, GNU GPL v3.0}.

\subsection{Shock Capturing for the DGSEM}\label{subsec:fv}
As stated earlier, high order methods tend to produce oscillations near discontinuities, which conflict with the locally smooth approximation space. Thus, a form of stabilization or regularization has to be introduced. In this work, we apply a finite volume sub-cell scheme~\cite{sonntagdiss} at elements which contain a discontinuity. In these cells, the polynomial solution representation with $(\mathcal{N}+1)^d$ degrees of freedom is projected onto $(\mathcal{N}+1)^d$ equidistantly distributed piecewise constant sub-cells, i.e. an equivalent FV solution representation is constructed. The projection is realized through the Vandermonde matrix $\mathcal{\vec{V}}_{\text{FV}}$. For each sub-cell, a second order TVD finite volume scheme is applied, and the coupling between discontinuous Galerkin elements and finite volume elements is achieved in a natural way via numerical flux functions. This type of shock capturing introduces a stable representation of discontinuities through the reduction of the order of the solution representation ($p$-adaptation), while at the same time preventing a significant loss in resolution by increasing the local element numbers through the introduction of sub-cells ($h$-adaptation). If a cell no longer contains a discontinuity, it is switched back to the DG representation by multiplication with $\mathcal{\vec{V}}_{\text{FV}}^{-1}$. Note that with regards to the discussion in Sec.~\ref{subsec:sota} and~\ref{subsec:motivation}, we treat this shock capturing scheme as a black box and use it as is for the remainder of this work.

\subsection{Shock Indicators}\label{subsec:indicators}
\begin{figure}
	\centering
	\includegraphics[width=0.4\textwidth]{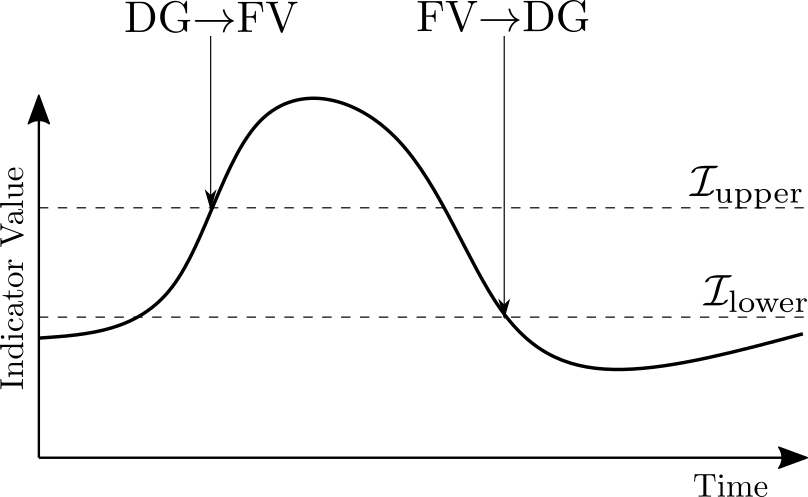}
	\caption{ Indicator based switching for shock capturing in the hybrid DG/FV sub-cell method.\label{fig:indicators}}
\end{figure}

To indicate when to switch from/to the polynomial representation of the discontinuous Galerkin method to/from the piecewise constant finite volume representation, a shock detection function is required. Besides the novel shock indicator introduced later in Sec.~\ref{sec:NeuralNetworks}, a wide range of indicators exits which can be used to identify elements containing a discontinuity, see e.g.~\cite{dolejvsi2003,ducros1999,huerta2012,jameson1981,persson2006} and the discussion in Sec.~\ref{subsec:sota}. Typically the indicators rely on the definition of a user-defined value against which the calculated indicator value can be compared. Due to the interactions of solution and discretization, these values are specific to a certain setup. For our calculations, we define an upper and a lower threshold introducing a delay in the switching procedure to avoid an oscillatory behavior if the indicator value almost matches the threshold. This means that a cell is only switched to the finite volume representation if the indicator value is greater than the upper threshold $\mathcal{I}_\text{upper}$ and is switched back to the polynomial representation if the indicator value is below the lower threshold $\mathcal{I}_\text{low}$. Fig.~\ref{fig:indicators} depicts a typical evolution of the indicator and its thresholds.

 We compare our ANNSI indicator against two well-known formulations from literature:
\begin{itemize}
	\item Indicator based on a modal smoothness analysis in analogy to~\cite{huerta2012,persson2006} taken from~\cite{sonntagdiss}:\newline
	Given the solution in modal polynomial space, this indicator evaluates the relative contribution of the highest modes for the definition of the indicator value $\mathcal{I}_\text{modal}$.
	\begin{align*}
		\mathcal{I}_\text{modal}=\text{log}_{10}\text{max}\Bigg\{ \frac{([\vec{w}_\text{modal}]_i^i,[\vec{w}_\text{modal}]_i^i)_{L_2}}{([\vec{w}_\text{modal}]_0^i,[\vec{w}_\text{modal}]_0^i)_{L_2}},i=\mathcal{N}-2,\mathcal{N}-1,\mathcal{N}\Bigg\},
	\end{align*}
	with $(\cdot,\cdot)_{L_2}$ being the $L_2$ inner product and the truncation to modes between $a$ and $b$ being
	\begin{align*}
		[\vec{w}_\text{modal}]_a^b=\sum_{i=n(a)-1}^{n(b)-1}\vec{w}_{\text{modal},i}\vec\psi_i,
	\end{align*}
	where $n(b)$ denotes the total number of the Legendre basis functions $\psi_i$ for the polynomial degree $b$.
	\item Indicator based on the jump inside an element inspired by~\cite{jameson1981} given in~\cite{sonntagdiss}:\newline
	For two dimensional simulations the indicator value can be calculated as
	\begin{align*}
		\mathcal{I}_\text{jump}=\frac{1}{V_E}\sum_{i,j=0}^{\mathcal{N}}\frac{w_{\text{min},ij}-2w_{ij}+w_{\text{max},ij}}{w_{\text{min},ij}+2w_{ij}+w_{\text{max},ij}}V_{ij},
	\end{align*}
	where $V_E$ and $V_{ij}$ denote the volume of the element and the finite volume sub-cell with the index $i,j$, respectively. With $w$ being a user defined variable of $\vec{w}$ e.g. pressure or density, $w_{\text{min},ij}$ is defined as $w_{\text{min},ij}=\text{min}(w_{i\pm id,j\pm jd},d=1,2)$, denoting the minimal value of the chosen variable neighboring the current degree of freedom. $w_{\text{max},ij}$ is defined in an analogous way.
\end{itemize}

\section{\label{sec:data}Neural Networks Design, Data Generation and Training}\label{sec:NeuralNetworks}
As outlined in Sec.~\ref{subsec:motivation}, our goal is to develop shock detection and localization indicators through methods from supervised machine learning, in particular image-based edge detection methods. Details on these methods and their training process will be given in this section. Before discussing the details, it is helpful to outline the steps: For the sake of completeness, we briefly introduce general neural networks in Sec.~\ref{subsec:nnintro}. For the specific task of edge detection, we then discuss the holistically nested edge detection (HED) method in Sec.~\ref{sec:NNArchitecture} which improved the state of the art when it was introduced in 2015 and is still considered among the best edge detection methods, and present our designed variant of HED. We discuss how we employ these networks as indicators for both shock detection and shock localization. The data set for the supervised learning of variants suitable for DG methods is presented in Sec.~\ref{subsec:trainingdata}, followed by information on the training process and the achieved results in Sec.~\ref{sec:NNTraining}.
\subsection{Function Approximation via Neural Networks\label{subsec:nnintro}}
\begin{figure}
	\centering
	\includegraphics[]{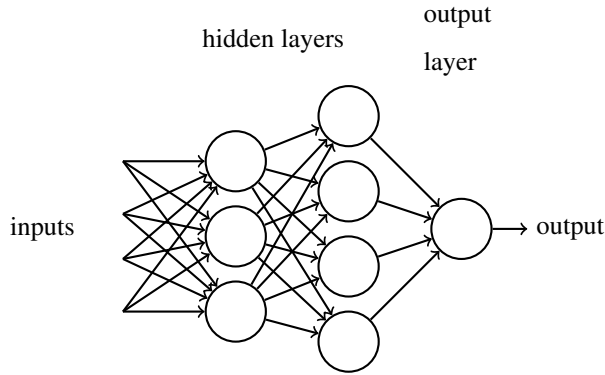}
	\caption{DAC of a classical MLP neural network.\label{fig:nn}}
\end{figure}

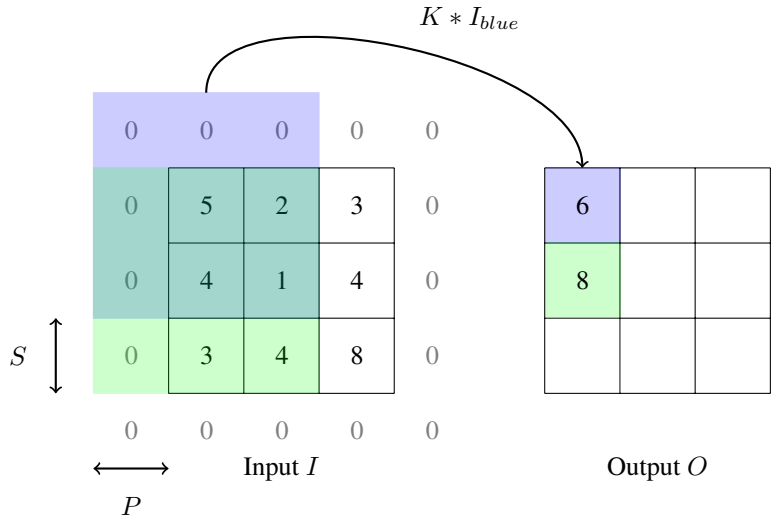
\begin{figure}
	\centering
	\begin{tikzpicture}
		\matrix[matrix of nodes,inner sep=0pt,anchor=south west,nodes={inner sep=0pt,text width=1cm,align=center,minimum height=1cm},text=gray] at (-1,-1) {
			0 & 0 & 0 & 0 & 0 \\
			0 &  &  &  & 0 \\
		    0 &  &  &  & 0 \\
			0 &  &  &  & 0 \\
			0 & 0  & 0 & 0 & 0 \\};
		\draw[step=1cm,color=black] (0,0) grid (3,3);
		\matrix[matrix of nodes,inner sep=0pt,anchor=south west,nodes={inner sep=0pt,text width=1cm,align=center,minimum height=1cm}] at (0,0) {
			5 & 2 & 3 & \\
			4 & 1 & 4 & \\
			3 & 4 & 8 & \\};
		\filldraw[color=blue,fill=blue,opacity=0.2] (-1,1) rectangle (2,4);
		\filldraw[color=green,fill=green,opacity=0.2] (-1,0) rectangle (2,3);
		\draw[step=1cm,color=black] (5,0) grid (8,3);
		\matrix[matrix of nodes,inner sep=0pt,anchor=south west,nodes={inner sep=0pt,text width=1cm,align=center,minimum height=1cm}] at (5,0) {
			|[fill=blue,opacity=0.2]| &  &  & \\
			|[fill=green,opacity=0.2]| &  &  & \\
			|[fill=white,opacity=0.05]|  &  &  & \\};
		\matrix[matrix of nodes,inner sep=0pt,anchor=south west,nodes={inner sep=0pt,text width=1cm,align=center,minimum height=1cm}] at (5,0) {
			6 &  &  & \\
			8 &  &  & \\
			|[fill=white,opacity=0.05]| &  &  & \\};
		\draw[thick, ->] (0.5,4) .. controls +(up:15mm) and +(up:15mm) .. (5.5,3);
		\node[] at (4,5) {$K \ast I_{blue}$};
		\draw[thick, <->] (-1,-1) -- (0,-1);
		\node[] at (-0.5,-1.5) {$P$};
		\draw[thick, <->] (-1.5,0) -- (-1.5,1);
		\node[] at (-2,0.5) {$S$};
		\node[] at (1.5,-1) {Input $I$};
		\node[] at (6.5,-1) {Output $O$};
\end{tikzpicture}
	\caption{Illustration of a convolutional operation. The filter kernel $K$ is of size $3\times 3$, the colors indicate the corresponding filter input and output neurons. For demonstration purposes, the kernel is chosen the identity matrix, for completeness stride (S) and padding (P) are chosen as  $S$ = $P$ = 1.\label{fig:conv}}
\end{figure}

The number of possible artificial neural network (ANN) architectures is large and currently growing~\cite{nnzoo}. We will restrict our discussion herein to the common forms of multilayer perceptrons and convolutional networks.
In their most general form, ANNs are multivariate compound functions that map an input vector $\vec X$ to an output vector $\vec{ Y}$. The mapping itself consists of a sequence of linear and non-linear functions with adjustable parameters. The common visual representation of a neural network as a directed acyclic graph (DAC) encodes the way in which these functions are stacked together, see Fig.~\ref{fig:nn} for an example of an MLP network. Here, arrows between two adjacent layers of neurons represent the entries of the linear weight matrices, with subsequent non-linear activation functions applied to each neuron before feeding into the next layer. In supervised learning, the weights of the network are adjusted through incremental learning from labeled training data, where the approximation quality of the neural network is evaluated by a cost function. The gradient of the cost function w.r.t. the parameters then guides the optimization process to determine the best functional fit to the data. Under mild assumptions on the underlying mapping and the architecture of the network, neural networks with a single or more hidden layers containing a finite number of neurons, can be shown to be universal function approximators~\cite{barron1993universal,cybenko1989approximation,HORNIK1991251,NIPS2017_7203}. For more details on neural networks, their properties and the training process in particular for deep learning, we refer the reader to~\cite{haykin1994neural,schmidhuber2015deep}.
An alternative architecture to the multilayer perceptrons with dense layer connections are the so-called convolutional neural networks. They were originally developed for computer vision tasks~\cite{lecun1990handwritten,lecun1998gradient,krizhevsky2012imagenet}. This architecture retains the layer structure but replaces the dense matrix connections between the individual layers by a local convolution of the input with a filter (see Fig.~\ref{fig:conv}). The entries of these filter kernels are then the parameters to be optimized. These trainable, local filter functions equip the CNNs with some favorable properties over the classical MLPs. First, the sparseness of the inter-layer connection greatly reduces the amount of parameters, making training of large networks feasible. Secondly, due to the local convolution, they are invariant to spatially shifted inputs. As one can interpret a CNN as a method of finding a suitable reduced order basis based on data snapshots, it is not surprising that CNNs have shown to be particularly powerful when the underlying data has a hierarchical structure, as for example in image recognition and segmentation~\cite{lecun1995convolutional} or turbulence modelling~\cite{beck2019deep}.

\subsection{Neural Network Architecture}\label{sec:NNArchitecture}
Due to these favorable properties of CNNs, they are widely used in edge detection (and other image-based tasks). As we interpret the task of shock detection and localization as one from this field, our network design is inspired by the holistically-nested edge detection (HED) network from~\cite{xie2015}. We use the same architecture for both tasks, i.e. the network design for ANNSI and ANNSL is \emph{identical}. Differences between ANNSI and ANNSL exists only in the training data and will be discussed later on.

\subsubsection{Edge Detection through HED\label{subsec:hed}}
The detection of the position and shape of a shock wave within an element is akin to edge detection in images in that a steep gradient or discontinuity separates the original image into distinct regions. We thus use ideas from this area to generate a network capable of localizing shocks in the 2D "image" made up by the solution points in each element. One particularly successful approach that has advanced the state of the art on reference data set is the HED algorithm proposed in~\cite{xie2015}. It has a range of properties that make it suitable for the task at hand. The first one is indicated by the "holistic" attribute, which refers to the fact that HED outputs an image of the same size as the input, with the edges marked on it. This is directly useful for the solution representation and the inner-element interpolation points that exist in element-based high order discretizations. The second feature, covered by the term "nested", refers to the architecture of the HED network. First of all, it is based on a convolutional neural network design which makes it shift-invariant, which is clearly desirable for shock localization. Due to the limited kernel size in each layer of these networks, the total receptive field of each layer increases along the DAC. Thus, different layers act on different scales of the input in a local-to-global manner. This multiscale property of CNNs is used in HED by also considering the side outputs of each layer in the cost function (and thus their parameters in the network training), recognizing that an edge should be scale-invariant. Each layer is thus thought of as a singleton network with the task of locating the edge at the given scale. This notion in HED was inspired by deeply-supervised nets~\cite{lee2015deeply}. In addition to supervision, the predictions of the sideoutputs are combined in a weighted manner and fused with the overall network output, again taking into account the scale invariance of edges. Although HED incorporates these powerful features, its actual network design can easily be derived from standard deep CNNs and is easy to describe and implement.
\subsubsection{ANNSI/ANNSL Architecture}
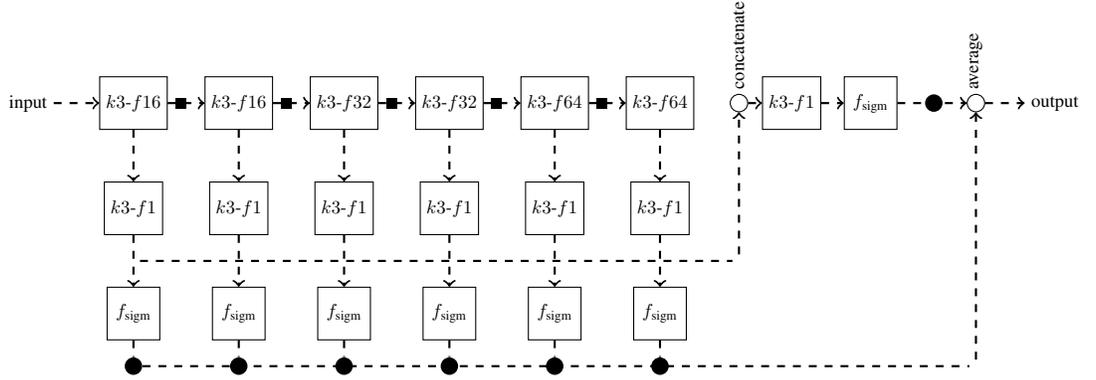
\begin{figure}[!htb]
	\centering
	\begin{tikzpicture}
	[node distance = 2.0cm,every node/.style={scale=0.7},squarednode/.style={rectangle, draw=black, fill=white, minimum size=10mm}, roundnode/.style={circle, draw=black, fill=white, minimum size=1mm},dot/.style={circle, fill=black, minimum size=0.5mm}]
	\node[] at (0,10.3) (00) {input};
	\node[squarednode,right of=00] (01) {$k3$-$f16$};
	\node[squarednode,right of=01] (02) {$k3$-$f16$};
	\node[squarednode,right of=02] (03) {$k3$-$f32$};
	\node[squarednode,right of=03] (04) {$k3$-$f32$};
	\node[squarednode,right of=04] (05) {$k3$-$f64$};
	\node[squarednode,right of=05] (06) {$k3$-$f64$};
	\node[roundnode, right of=06,node distance=1.5cm] (07) {};
	\node[squarednode,right of=07,node distance=1.0cm] (08) {$k3$-$f1$};
	\node[squarednode,right of=08,node distance=1.5cm] (09) {$f_\text{sigm}$};
	\node[dot,right of=09,node distance=1.2cm] (dot) {};
	\node[roundnode, right of=09,node distance=2.0cm] (099) {};
	\node[right of=099,node distance=1.5cm] (output) {output};

	\node[squarednode,below of=01] (11) {$k3$-$f1$};
	\node[squarednode,below of=02] (12) {$k3$-$f1$};
	\node[squarednode,below of=03] (13) {$k3$-$f1$};
	\node[squarednode,below of=04] (14) {$k3$-$f1$};
	\node[squarednode,below of=05] (15) {$k3$-$f1$};
	\node[squarednode,below of=06] (16) {$k3$-$f1$};

	\node[squarednode,below of=11] (21) {$f_\text{sigm}$};
	\node[squarednode,below of=12] (22) {$f_\text{sigm}$};
	\node[squarednode,below of=13] (23) {$f_\text{sigm}$};
	\node[squarednode,below of=14] (24) {$f_\text{sigm}$};
	\node[squarednode,below of=15] (25) {$f_\text{sigm}$};
	\node[squarednode,below of=16] (26) {$f_\text{sigm}$};

	\draw[thick,dashed, ->] (00) -- (01);
	\draw[thick,dashed, ->] (01) -- (02);
	\draw[thick,dashed, ->] (02) -- (03);
	\draw[thick,dashed, ->] (03) -- (04);
	\draw[thick,dashed, ->] (04) -- (05);
	\draw[thick,dashed, ->] (05) -- (06);
	\draw[thick,dashed, ->] (07) -- (08);
	\draw[thick,dashed, ->] (08) -- (09);
	\draw[thick,dashed, ->] (09) -- (099);
	\draw[thick,dashed, ->] (099) -- (output);

	\draw[thick,dashed, ->] (01) -- (11);
	\draw[thick,dashed, ->] (02) -- (12);
	\draw[thick,dashed, ->] (03) -- (13);
	\draw[thick,dashed, ->] (04) -- (14);
	\draw[thick,dashed, ->] (05) -- (15);
	\draw[thick,dashed, ->] (06) -- (16);

	\draw[thick,dashed, ->] (11) -- (21);
	\draw[thick,dashed, ->] (12) -- (22);
	\draw[thick,dashed, ->] (13) -- (23);
	\draw[thick,dashed, ->] (14) -- (24);
	\draw[thick,dashed, ->] (15) -- (25);
	\draw[thick,dashed, ->] (16) -- (26);

	\node[below of=11,node distance=1.0cm] (m) {};
	\node[below of=07,node distance=3.0cm] (mm) {};
	\draw[thick,dashed, - ] (m) -- (mm);
	\draw[thick,dashed, ->] (mm)-- (07);

	\node[dot,below of=21,node distance=1.0cm] (11m) {};
	\node[dot,below of=22,node distance=1.0cm] (12m) {};
	\node[dot,below of=23,node distance=1.0cm] (13m) {};
	\node[dot,below of=24,node distance=1.0cm] (14m) {};
	\node[dot,below of=25,node distance=1.0cm] (15m) {};
	\node[dot,below of=26,node distance=1.0cm] (16m) {};
	\node[below of=099,node distance=5.0cm] (mmm) {};

	\draw[thick,dashed, - ] (21) -- (11m);
	\draw[thick,dashed, - ] (22) -- (12m);
	\draw[thick,dashed, - ] (23) -- (13m);
	\draw[thick,dashed, - ] (24) -- (14m);
	\draw[thick,dashed, - ] (25) -- (15m);
	\draw[thick,dashed, - ] (26) -- (16m);
	\draw[thick,dashed, - ] (11m) -- (mmm);
	\draw[thick,dashed, ->] (mmm) -- (099);

	\node[above of=07,rotate=90, node distance=1.1cm] (concatenate) {concatenate};
	\node[above of=099,rotate=90, node distance=0.8cm] (average) {average};

	\node[fill=black,rectangle,minimum width=0.2cm,minimum height=0.2cm,inner sep=0pt,right of=01,node distance=0.9cm] (s1) {};
	\node[fill=black,rectangle,minimum width=0.2cm,minimum height=0.2cm,inner sep=0pt,right of=02,node distance=0.9cm] (s1) {};
	\node[fill=black,rectangle,minimum width=0.2cm,minimum height=0.2cm,inner sep=0pt,right of=03,node distance=0.9cm] (s1) {};
	\node[fill=black,rectangle,minimum width=0.2cm,minimum height=0.2cm,inner sep=0pt,right of=04,node distance=0.9cm] (s1) {};
	\node[fill=black,rectangle,minimum width=0.2cm,minimum height=0.2cm,inner sep=0pt,right of=05,node distance=0.9cm] (s1) {};
	\end{tikzpicture}
	\caption{Network architecture used for the shock detection (ANNSI) and localization (ANNSL). Each convolution layer has a kernel size $k$ and $f$ channels. Layers with $f=16,32,64$ use LReLU as activation function; in layers with $f=1$ no activation function is used. Black dots symbolize positions where the output is taken for the evaluation of the loss function $C$ and squares symbolize batch normalization.}\label{fig:ANNSL_Architecture}
\end{figure}

While the network design is inspired by HED, our task at hand differs from the intended applications in terms of scale: While in edge detection in images, the input image is typically of size $\mathcal{O}(100)^2$ pixels (disregarding the color channel) or even larger, in our intended application, it is of size $(\mathcal{N}+1)^2$, which, for typical $\mathcal{N}<10$. Thus, the network design itself has been adjusted to reflect this much smaller input space. \\
An overview of the network architecture used in this work is provided in Fig.~\ref{fig:ANNSL_Architecture}. At its core, the upper lane through the network is a standard deep CNN architecture. The side outputs (middle and lower lane) compute the output predictions of each layer of the first lane, which are turned into edge maps with the same size as the input image.  Information from all of the edge maps is fused with the final output, thus taking advantage of the scale invariance of edges. \\
For the specific network used in our context and shown in Fig.~\ref{fig:ANNSL_Architecture}, we use the following technical details:  In layers with $f=16,32,64$ channels, we use a \emph{leaky Rectified Linear Unit}~\cite{maas2013} (LReLU)
\begin{align*}
f_\text{LReLU}(z)=\text{max}(0,z)+\alpha\text{min}(0,z)\quad\text{with}\quad\alpha=0.2,
\end{align*}
as activation function, whereas in layers with $f=1$ no activation function is applied as they are connected to a \emph{sigmoid} function
\begin{align*}
f_\text{sigm}(z)=\frac{1}{1+e^{-z}}.
\end{align*}
 As common, we use batch normalization~\cite{ioffe2015batch} to enhance the training process and use padding $P=1$ and stride $S=1$ to keep the size of the data constant for the different layers, which is required for the "holistically-nested" approach. The filter kernel size is chosen to be $3$, which we found to work best for the limited size of the input image. In total, the network used has 72477 trainable parameters.
\subsection{Training Data\label{subsec:trainingdata}}

For the training of our neural network, we stick to the \emph{supervised learning} strategy, see e.g.~\cite{lecun2015}. Under the assumption that there exists a functional dependency $\vec{f}:\vec{X}\mapsto \vec{Y}$ of the input vector $\vec X$ and the true desired output $\vec Y$, this paradigm states that one tries to find an approximation $\vec{\hat{f}}:\vec X\mapsto \vec{\hat{Y}}$, such that the difference between the truth $\vec Y$ and the predicted output $\vec{\hat{Y}}$ is minimal. The metric for this difference is called cost function $C:=C(\vec{Y},\vec{\hat{Y}})$, and a single evaluation of $\vec{\hat{f}}:\vec X\mapsto \vec{\hat{Y}}$ given the current mapping encoded in the network is called a forward pass. During the backward pass, the gradient of the cost function w.r.t. all the parameters of the network is computed, and the weights are updated following a gradient descent optimization~\citep{rumelhart1986learning,werbos}. Often, this weight update is computed not for single training samples (or all of them), but for so-called \emph{mini batches}~\citep{ioffe2015batch}. The size of these mini batches is chosen empirically and should provide a compromise between convergence rate, stability, memory demands and accuracy; one training iteration over all training samples is called an epoch. The magnitude of the weight increment depends on a scaling parameter termed~\emph{learning rate}. Since the sought mapping is typically non-convex with respect to the parameters of the network, neural networks are considerably data-hungry in their training process.
\subsubsection{Training Data for ANNSI\label{subsec:trainingdata_annsi}}
For the training of a shock detector, this requires the definition of a sufficiently broad spectrum of input images $\vec{X}$ and a corresponding label indicating the true desired output image $\vec Y$ for each date. Note that since the network architecture is of convolutional type, the input into the network are not restricted to vectors. On the contrary, 2D (or higher dimensional volumetric data) are naturally suited as inputs and outputs of these nets. Thus, with our chosen limitation to 2D problems, the input consists of square images with $(\mathcal{N}+1)^2$ pixels, where each pixel corresponds to a nodal degree of freedom of the DG approximation. In other words, in true image-based analysis fashion, we interpret the solution within an element as an image and the solution points as its pixels. \\
 Inspired by~\cite{ray2018,ray2019}, we use a set of smooth and non-smooth analytical functions to represent possible solutions of a simulation with the DGSEM and define on the analytic level if and where a shock is present in the current grid cell. We map all physical elements to reference space in $[-1,1]^2$, thus, for training, all inputs and outputs are given on the \emph{same reference square}. The analytical functions given in Tbl.~\ref{Tbl:exactfunc_2D} in \ref{app:data} are evaluated on Cartesian meshes covering the interval $[-1,1]^2$ with $n_e^2$ elements. A class label of 0 refers to a no-shock solution, a label of 1 to a shock solution. In each of the $n_e^2$ elements, the function is approximated by the tensorproduct basis with $(\mathcal{N}+1)^2$ solution points located at the Legendre Gauss nodes as described above. For the case of  $n_e=1$, the domain coincides with the reference element $[-1,1]^2$, for other choices of $n_e$, the resulting solution is mapped back to the reference element. This provides an efficient and easy way to introduce the effect of different spatial resolutions into the training process.

 With this information, we define the true desired output $\vec{Y}$, which is a so-called \emph{binary shock edge map}, having the same size as the input image $\vec{X}$. The entries of the edge map are $1$ for a degree of freedom / solution point being located at the discontinuity, otherwise it is labeled as $0$.  However, the analytical position of the shock does not generally coincide with the location of the solution points, instead, it will typically pass between two adjacent points. In this case, the two solution points adjacent to the shock front per each coordinate direction are given the $1$ label. This is akin to the definition of a shock bounding box on a sub-element level. Thus, labeled true data $\vec{Y}$ for each input $\vec{X}$ consists of the binary edge shock map of size  $(\mathcal{N}+1)^2$, where each entry is in $[0,1]$ to indicate that pixel being associated with a shock front or not. A class label (element contains a shock or not) can be derived from this for the whole element by $1\text{ if any}(\vec{Y})=1\text{, else }0$.

The total amount of data chosen from the different classes are summarized in Tbl.~\ref{Tbl:SamplesDetection}. The number of training cells with and without shocks are approximately balanced to enhance training. Note since the input tensors are dependent on the choice of the polynomial degree and to take the different spatial resolution capabilities into account, we choose to train a specific network for each $\mathcal{N}$.\\
\begin{table}[!thb]
	\begin{minipage}[c]{0.4\linewidth}
		\centering
		\subcaptionbox{$N=5$}{
		\begin{tabular}{lrrrr}
			$\#$     &   \multicolumn{2}{c}{Training Set} &   \multicolumn{2}{c}{Validation Set} \\ \toprule &   class $0$ &   class $1$ & class $0$ & class $1$ \\
			\midrule
			1        &                              20311 &                                                0 &        2046 &                                   0 \\
			2        &                              66288 &                                                0 &        6638 &                                   0 \\
			3        &                               6402 &                                                0 &         619 &                                   0 \\
			4        &                                 51 &                                            37266 &           6 &                                3758 \\
			5        &                               2377 &                                            16297 &         247 &                                1618 \\
			6        &                                102 &                                            32308 &           9 &                                3206 \\
			7        &                                 59 &                                            18551 &           8 &                                1853 \\
			\midrule &                              95590 &                                           104422 &        9573 &                               11530 \\
			\bottomrule
		\end{tabular}}
	\end{minipage}
	\hspace{20mm}
	\begin{minipage}[c]{0.4\linewidth}
		\centering
		\subcaptionbox{$N=9$}{
		\begin{tabular}{lrrrr}
			$\#$     &   \multicolumn{2}{c}{Training Set} &   \multicolumn{2}{c}{Validation Set} \\ \toprule &   class $0$ &   class $1$ & class $0$ & class $1$ \\
			\midrule
			1        &                              20400 &                                                0 &        2388 &                                   0 \\
			2        &                              66400 &                                                0 &        7289 &                                   0 \\
			3        &                               6400 &                                                0 &        1195 &                                   0 \\
			4        &                                 40 &                                            40023 &          95 &                                6494 \\
			5        &                               2067 &                                            15660 &         251 &                                2184 \\
			6        &                                116 &                                            31094 &          21 &                                4206 \\
			7        &                                 80 &                                            18520 &          2 &                                2938 \\
			\midrule &                              95503 &                                           105297 &       11712 &                               15822 \\
			\bottomrule
		\end{tabular}}
	\end{minipage}
	\caption{Number of samples chosen from the different analytical functions in the training and validation sets for the polynomial degrees a) $\mathcal{N}=5$ and b) $\mathcal{N}=9$ for the shock detector ANNSI. \label{Tbl:SamplesDetection}}
\end{table}

\begin{remark}
In the presented approach, the shock detection per element is achieved by training the network to predict the binary edge map within an element and decide that a shock is present if any of the pixels in the edge map show $1$. This might seems as a waste of information and training, as the predicted matrix is condensed to a single bit per cell. We choose this route for a number of reasons:
\begin{itemize}
\item  As discussed in Sec.~\ref{subsec:sotaml}, a number of publications report useful and accurate shock detectors based on MLP architectures. In our investigations however, the results obtained for shock detection with the ANNSI approach showed greater accuracy and robustness than training an element-wise classifier directly. In addition, we are also interested in the inner-cell shock location, see next item.
	\item While we choose to condense the edge map to a single value, its full information is available during the computation. For shock capturing schemes which do not change the localization of the solution nodes, this information can directly be used for shock localization. Such schemes are e.g. artificial viscosity approaches (see e.g.~\cite{BARTER20101810}) or sub-cell shock capturing on Legendre-Gauss nodes~\cite{sonntag2017efficient}. Thus, our choice of input and output data provides flexibility in using the full edge map information in ways that fit the chosen shock capturing approach.
	\item It avoids having to develop two different neural network architectures, thereby reducing the complexity of the approach.
\end{itemize}
\end{remark}

\subsubsection{Training Data for Shock Localization (ANNSL)}\label{sec:Training_Localization}
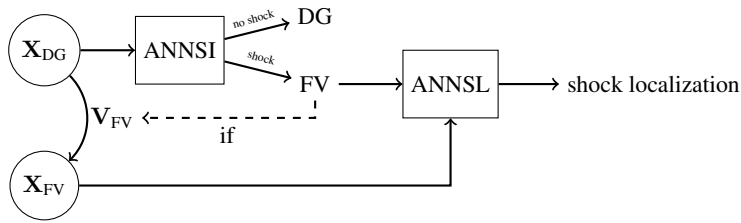
\begin{figure}[!htb]
	\centering
	\begin{tikzpicture}
	[node distance = 2.0cm,every node/.style={scale=0.9},squarednode/.style={rectangle, draw=black, fill=white, minimum size=10mm}, roundnode/.style={circle, draw=black, fill=white, minimum size=1mm},dot/.style={circle, fill=black, minimum size=0.5mm}]
	\node[roundnode] at (0,0) (DGin) {$\vec{X}_\text{DG}$};
	\node[roundnode, below of=DGin] (FVin) {$\vec{X}_\text{FV}$};
	\draw[thick,->, bend angle=45, bend left]  (DGin) to (FVin);
	\node[below of=DGin, right of=DGin, node distance=1.0cm] (Vdm) {$\mathcal{\vec{V}}_\text{FV}$};
	\node[squarednode,right of=DGin] (ANNSI) {ANNSI};
	\node[right of=ANNSI,above of=ANNSI,yshift=-1.5cm] (DGout) {DG};
	\node[below of=DGout,node distance=1.0cm] (FVout) {FV};
	\node[squarednode,right of=FVout] (ANNSL) {ANNSL};
	\node[right of=ANNSL, xshift=1.0cm] (local) {shock localization};

	\draw[thick,->] (DGin) to (ANNSI);
	\draw[thick,->] (ANNSI) -- (DGout) node[midway,above,rotate=16] {\tiny{no shock}};
	\draw[thick,->] (ANNSI) -- (FVout) node[midway,above,rotate=-16] {\tiny{shock}};
	\draw[thick,->] (FVout) to (ANNSL);
	\draw[thick,->] (ANNSL) to (local);
	\draw[thick,->] (FVin) -| (ANNSL);

	\draw[dashed,thick,->] (FVout) node [below, transform shape,xshift=-1.3cm,yshift=-0.5cm] {if}  |- (Vdm);

	\end{tikzpicture}
	\caption{Schematic overview of the workflow for shock localization.} \label{fig:Scheme_Localization}
\end{figure}

As discussed in the previous section, the  ANNSI network outputs a prediction of a shock map within each element, where the pixels correspond to the position of the Legendre Gauss solution points of degree $\mathcal{N}$ that host the DG solution polynomial. In theory, based on this information, any shock capturing scheme could be guided to improve the inner cell treatment of shocks. In the following, we discuss our strategy for exploiting this information in the context of our hybrid DGSEM / FV scheme: Keeping the network architecture and input functions fixed, we re-train the network but shift the position of the input and output pixels from $(\mathcal{N}+1)^2$ Legendre Gauss nodes to $(\mathcal{N}+1)^2$ equispaced nodes. This just requires evaluating the functions from Tbl.~\ref{Tbl:exactfunc_2D} on the new node distribution. The reason for this shift lies in a peculiarity of the chosen shock capturing strategy described in Sec.~\ref{subsec:fv}: The sub-cell finite volume scheme is formulated on an equispaced grid, not the Legendre Gauss grid. With increasing degree $\mathcal{N}$, the discrepancy between the equidistant and Legendre Gauss points increases, and local correspondence is lost. This makes transferring positional information from the FV to the DG grid and vice versa less accurate and robust. In other words, locating the shock on the DG grid is not very useful for shock capturing on the equispaced grid. For this reason, we train the ANNSL network directly on the FV grid. \\
In summary, our strategy for shock detection and localization works as follows: First, all cells are treated as DG cells with the corresponding solution points and the ANNSI network (trained on Legendre Gauss nodes) acts as a shock detector. Cells that are marked by ANNSI as shock-free remain DG elements for the next timestep. Those that are marked as containing shocks are treated by the FV shock-capturing scheme: The solution is projected onto the FV space using the Vandermonde matrix $\mathcal{\vec{V}}_{\text{FV}}$. This solution on equispaced points now informs the ANNSL network, which predicts the shock location within the respective element. Fig.~\ref{fig:Scheme_Localization} provides a schematic overview of this procedure.

We also note that it is possible to directly apply the ANNSL to the FV representation in all elements of the domain without pre-selection by ANNSI: This gives slightly less robust, however comparable results but requires the computational effort of computing the projection onto the FV space in all cells at each timestep. Since this approach however makes the initial ANNSI network superfluous, it can result in an overall speed-up. Since our main focus and motivation in this work is on accuracy and robustness of the approach and not so much on efficiency, we did not investigate this in greater detail.

To summarize our training rationale, we train two identical networks for the prediction of a binary edge map on two different sets of input data: the DG solution points (ANNSI) and equispaced FV (ANNSL) solution points. Both networks thus provide information about the  inner-cell shock position by marking the associated nodes / pixels. The choice to retrain the network on equispaced points is solely driven by the specific shock capturing scheme employed here. We choose to use the ANNSI network to pre-select elements containing shocks to show its capabilities as a pure shock detection scheme as well as to demonstrate the flexibility of our method with regards to the underlying properties of the discretization. The exact shock location on the FV grid is then conducted by the ANNSL net, which has been trained with suitable input data.
The training data is summarized in Tbl.~\ref{Tbl:SamplesLocalization}, using approximately the same number of samples per class as for the ANNSI net in Tbl.~\ref{Tbl:SamplesDetection}.

\begin{table}[!thb]
	\centering
	\begin{minipage}[c]{0.4\linewidth}
		\begin{tabular}{lrrrr}
			$\#$     &   \multicolumn{2}{c}{Training Set} &   \multicolumn{2}{c}{Validation Set} \\ \toprule &   class $0$ &   class $1$ & class $0$ & class $1$ \\
			\midrule
			1        &                              21327 &                                                0 &        4253 &                                   0 \\
			2        &                              69190 &                                                0 &       13769 &                                   0 \\
			3        &                               6706 &                                                0 &        1345 &                                   0 \\
			4        &                               4248 &                                            35364 &         861 &                                7115 \\
			5        &                               2281 &                                            17239 &         469 &                                3386 \\
			6        &                                 80 &                                            33858 &          17 &                                6922 \\
			7        &                                253 &                                            19454 &          37 &                                3826 \\
			\midrule &                             104085 &                                           105915 &       20751 &                               21249 \\
			\bottomrule
		\end{tabular}
	\end{minipage}
	\caption{Number of samples chosen from the different analytical functions in the training and validation sets for the polynomial degree $\mathcal{N}=9$ for the shock localization with ANNSL. \label{Tbl:SamplesLocalization}}
\end{table}

\subsection{Neural Network Training}\label{sec:NNTraining}
\begin{figure}
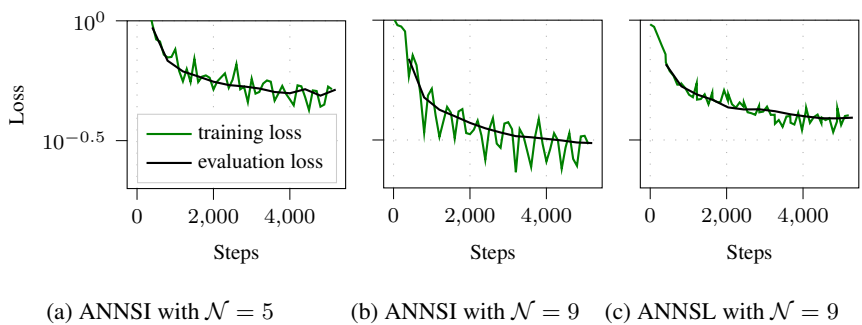

	\begin{minipage}[c]{0.37\linewidth}
	        \vspace{-0.1cm}
		\subcaptionbox{ANNSI with $\mathcal{N}=5$}{\includetikzgraphic[width=\linewidth]{./figures/tikz/loss_HED_N5.tikz}}
	\end{minipage}
	\hspace{-0.2cm}
	\begin{minipage}[c]{0.37\linewidth}
	        \hspace{+0.15cm}
		\subcaptionbox{ANNSI with $\mathcal{N}=9$}{\includetikzgraphic[width=\linewidth]{./figures/tikz/loss_HED_N9.tikz}}
	\end{minipage}
	\hspace{-1.50cm}
	\begin{minipage}[c]{0.37\linewidth}
	        \hspace{+0.40cm}
		\subcaptionbox{ANNSL with $\mathcal{N}=9$}{\includetikzgraphic[width=\linewidth]{./figures/tikz/loss_HED_N9_FV.tikz}}
	\end{minipage}
	\caption{Training and validation loss during the training of ANNSI (a/b) and ANNSL (c) network.}\label{fig:TrainingLoss}
\end{figure}

As discussed in Sec.~\ref{subsec:hed}, the HED architecture evaluates the edge prediction on different scales through the incorporation of side outputs into the cost function. Quantities that contribute to the cost function are marked as black dots in Fig.~\ref{fig:ANNSL_Architecture}. Following~\cite{liu2017}, an adapted weighted cross-entropy function is chosen
\begin{align}\label{eq:costfunction}
	C=-\frac{1}{m}\sum_{m}\sum_{l=1}^{L}\left[\frac{1}{S_Y}\sum_{s=1}^{S_Y}\Lambda\vec{Y}_s\text{ln}(\vec{\hat{Y}}_s^l)+\lambda(1-\Lambda)(1-\vec{Y}_s)\text{ln}(1-\hat{\vec{Y}}_s^l)\right]_m,
\end{align}
with $L$ being the number of input layers to the loss function (black dots in Fig.~\ref{fig:ANNSL_Architecture}), $S_Y$ being the size of the output vector, the weighting factor $\lambda=1.1$ and $\Lambda$ being defined as
\begin{align*}
	\Lambda=\frac{\#\text{total number of labels with class 1}}{\#\text{total number of labels}},
\end{align*}
for the mini-batch size $m$.
We start the training with an initial guess for the degrees of freedom of the neural network: All biases are chosen to be zero and the weights are initialized with randomly drawn numbers from a normal distribution with zero mean and standard deviation of $\sigma=0.01$, except for the last layer where the standard deviation of the normal distribution is chosen to be $\sigma=0.2$, following~\cite{liu2017}. The remaining parameters for the training process are the same for ANNSI and ANNSL and are summarized in Tbl.~\ref{Tbl:TrainingANNSI}. We exemplary report the results for ANNSI trained for $\mathcal{N}=5$ and $\mathcal{N}=9$ to cover both a medium and a high order scheme, as well as for ANNSL with $\mathcal{N}=9$.
\begin{table}[!thb]
	\centering
	\begin{tabularx}{0.8\textwidth}{|l|X|}
		\hline
		training optimizer & Adam~\cite{kingma2014} \\ \hline
		mini-batch size & $m=500$ \\ \hline
		epochs & $n_\text{epoch}=120$ \\ \hline
		learning rate & $\text{LR}=0.01$\newline reduced every $15^{th}$ epoch by factor $0.5$ \\
		\hline
	\end{tabularx}
	\caption{Parameters for the training of the ANNSI and ANNSL network}\label{Tbl:TrainingANNSI}
\end{table}
The influence of the training procedure on the loss evaluated with Eq.~\eqref{eq:costfunction} is visualized in Fig.~\ref{fig:TrainingLoss}. It shows a monotone decreasing validation loss, indicating that no overfitting takes place.
In order to quantify the accuracy of the trained networks, we choose the $F_1$ score, a custom metric for binary classifiers which incorporates true/false positives and negatives into the metric. The ANNSI networks reach a $F_1$ score of $96.461\%$ for $\mathcal{N}=5$ and $96.679\%$ for $\mathcal{N}=9$. With the ANNSL network for $\mathcal{N}=9$ a $F_1$ score of $95.028\%$ is reached. During training, we achieved $F_1$ scores of over $90\%$ consistently with much smaller networks. We choose to focus on the best network in this work only to evaluate the full potential of our approach.

\section{\label{sec:NNshock} Results: Shock Detection}
With the ANNSI networks for $\mathcal{N}=5$ and $\mathcal{N}=9$ from Sec.~\ref{sec:NNTraining} in place, we now apply the trained shock indicators to a range of 2D shock problems. We first describe the selected flow cases briefly in Sec.~\ref{subsec:testcases}. With these definitions out of the way, we then discuss the results for the ANNSI approach and contrast it against the established indicators from Sec.~\ref{subsec:indicators}. We first focus here on the shock detection aspect part of the proposed indicator, i.e. the marking of grid elements where shocks likely reside. The aspect of shock localization within an element will be discussed in Sec.~\ref{sec:NNshock_local}. For all computations presented herein, we rely on the discretization scheme, including the TVD shock capturing method, presented in Sec.~\ref{subsec:dg} and~\ref{subsec:fv}. The polynomial ansatz degree of the scheme is either $\mathcal{N}=5$ and $\mathcal{N}=9$ as indicated; the scheme is a standard DGSEM collocation formulation on Legendre-Gauss points, with a Roe or HLLE scheme for the numerical flux functions of both DG and FV elements. A minmod limiter ensures the TVD property of the linear FV reconstruction polynomials. The solution is advanced in time by an explicit high order Runge-Kutta scheme~\cite{kennedy2000low}; the typical CFL condition limits the allowable global timestep~\cite{krais2019flexi}. The density is chosen as an input quantity into ANNSI.\\
All simulations with ANNSI or ANNSL shown in this manuscript are run ab initio, i.e. from the given initial conditions at time $t_{start}=0.0$ to the indicated $t_{end}$, using the proposed indicators only. While we focus on discussing the solution quality at a given time typical for the selected test case, this shows that ANNSI and ANNSL approaches are capable of dealing with highly transient situations without the need for any special treatment.
\subsection{\label{subsec:testcases} Test Cases}
\subsubsection{\label{subsec:RP} 2D Riemann Problems}
In~\cite{schulz1993numerical}, the definition of the 1D Riemann (RP) problem for the Euler system is extended to two dimensions. The authors identify 16 distinct configurations in which constant states in the four domain quadrants are separated initially by an elementary wave and propose them as canonical test cases. We choose three configurations (4,6,12), which among them include the three typical wave phenomena, shocks, rarefactions and contact discontinuities. The initial conditions and end times of the RPs are given in Tbl.~\ref{tab:RmP_initial_conditions}; they define the solution in the four quadrants of the domain (numbered from the upper right quadrant counterclockwise). The boundary conditions on all domain faces are set to the exact solutions.

\begin{figure}[!htbp]
	\begin{minipage}[c]{0.32\linewidth}
		\caption*{Configuration 4}
		\includegraphics[width=\linewidth]{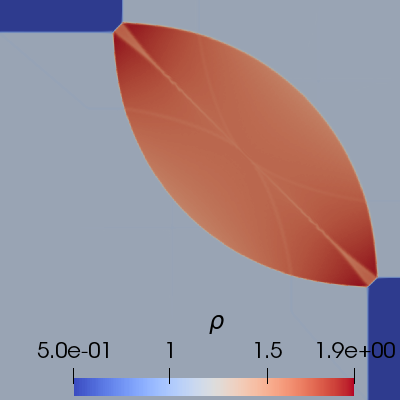}
	\end{minipage}
	\begin{minipage}[c]{0.32\linewidth}
		\caption*{Configuration 6}
		\includegraphics[width=\linewidth]{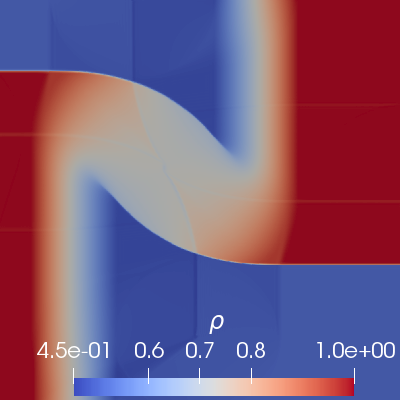}
	\end{minipage}
	\begin{minipage}[c]{0.32\linewidth}
		\caption*{Configuration 12}
		\includegraphics[width=\linewidth]{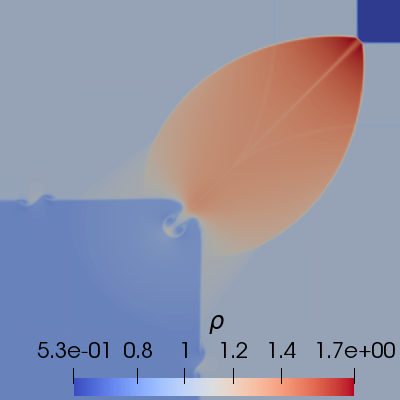}
	\end{minipage}
	\caption{Reference solution for the density contours of configuration 4, 6 and 12 of the two-dimensional Riemann problem. Computed on a $100 \times 100$ mesh with a polynomial degree of $\mathcal{N}=5$. \label{fig:2D_RMP_density}}
\end{figure}

\begin{table}[!htbp]
	\centering
	\begin{tabular}{llll}
		\toprule
		\# & Configuration 4 & Configuration 6 & Configuration 12 \\
		\midrule
		1 & (1.1, 0, 0, 1.1) & (1, 0, 0, 1) & (0.5313, 0, 0, 0.4) \\
		2 & (0.5065, 0.8939, 0, 0.35) & (0.5, -0.8708, 0, 0.3636) & (1, 0.7276, 0 ,1) \\
		3 & (1.1, 0.8939, 0.8939, 1.1) & (1, -0.8708, 0.7977, 1) & (0.8, 0, 0, 1) \\
		4 & (0.5065, 0, 0.8939, 0.35) & (0.5, 0, 0.7977, 0.3636) & (1, 0, 0.7276, 1) \\
		\bottomrule
		$t_{end}$& 0.25 & 0.2 & 0.25 \\
		\bottomrule
	\end{tabular}
	\caption{Initial conditions of $(\rho, v_1, v_2, p)$ in the four quadrants of the domain of the two-dimensional Riemann problems~\cite{schulz1993numerical}. \label{tab:RmP_initial_conditions}}
\end{table}
Fig.~\ref{fig:2D_RMP_density} shows the resulting density contours for the three configurations as a reference, computed on a fine grid with the numerical scheme described above. Note that along the contact discontinuities in configuration 12, a Kelvin-Helmholtz instability triggers the roll-up of vortical structures. For the comparisons of the shock indicators later on, we define a baseline regular Cartesian mesh $M_{RP}$ with $50\times 50$ elements for all configurations.
\subsubsection{\label{subsec:dmr} Double Mach Reflection Problem}
\begin{figure}[!htbp]
	\includegraphics[width=\linewidth]{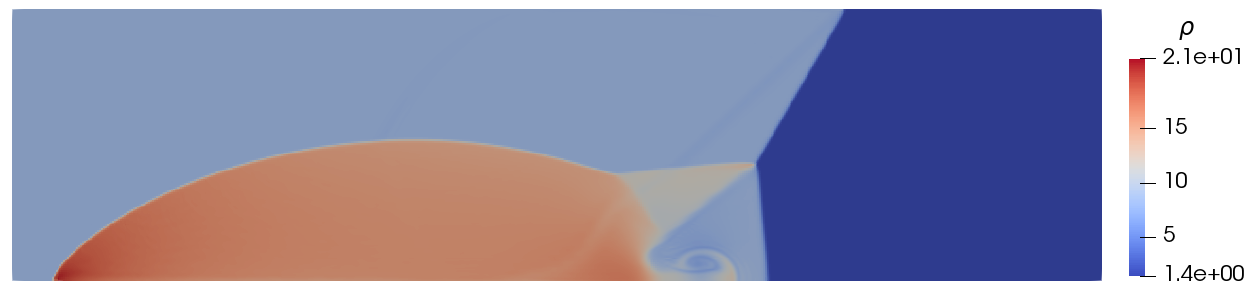}
	\caption{Reference solution of the density of DMR computed on a $96 \times 24$ mesh at $t_{end}=0.2$. \label{fig:dmr_density}}
\end{figure}

This testcase consists of a strong shock at Mach 10 impinging on a wedge at an oblique angle of $30^{\circ}$. Our setup is consistent with the description in~\cite{woodward1984numerical}, with slip-wall boundary conditions at the top and bottom, a Dirichlet inflow on the left and a supersonic outflow at the right hand side of the domain. Fig.~\ref{fig:dmr_density} shows the reference solution at $t_\text{end}=0.2$ computed on a fine grid. We define a coarse baseline grid $M_{DMR}$ with $49\times 12$ elements for later use in this study. The numerical flux is approximated by a Roe scheme.

\subsubsection{\label{subsec:ffs} Forward Facing Step Problem}
\begin{figure}
	\begin{minipage}[c]{0.8\linewidth}
		\includegraphics[trim={8cm 1cm 2cm 0},clip,width=\linewidth]{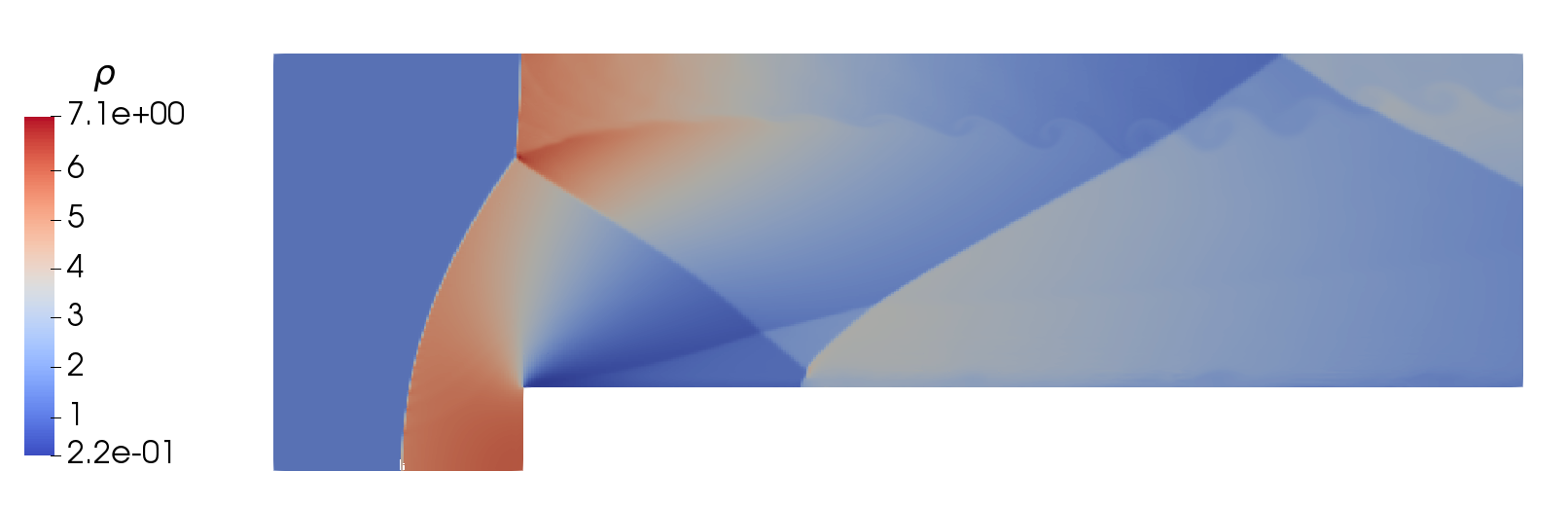}
	\end{minipage}
	\begin{minipage}[c]{0.1\linewidth}
		\begin{tikzpicture}[line cap=round,line join=round,x=1.15cm,y=1.15cm,axis/.style={->}]
		\node at (0 ,-0.4) {\includegraphics[scale=0.1]{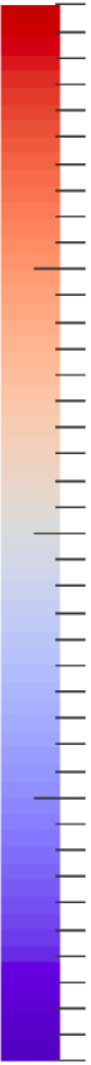}};
		\node at (0.0,1.6 ) {$~~~\rho~~~$};
		\node at (0.5,1.2 ) {$7.1$};
		\node at (0.5,0.4 ) {$5.5$};
		\node at (0.5,-0.4 ) {$3.7$};
		\node at (0.5,-1.2 ) {$1.8$};
		\node at (0.5,-2.0 ) {$0.22$};
		\end{tikzpicture}
	\end{minipage}
	\caption{Reference solution of the density of the FFS at $t_{end}=4$.}\label{fig:FFS}
\end{figure}

This test setup, introduced in~\cite{woodward1984numerical} describes the flow over a forward facing step (FFS) with a Mach number of $\mathrm{Ma}=3.0$. The expected flow phenomena are mainly a developing bow shock in front of the step and the reflections of this shock from the upper and lower wall. Inflow and outflow are described with Dirichlet boundary conditions, other boundaries are modeled as walls. In Fig.~\ref{fig:FFS} the density at $t_\text{end}=4$ is visualized for a calculation with $50$ elements in the $x$-direction and $25$ and $20$ elements in the two blocks in $y$-direction using $\mathcal{N}=9$ and an HLLE Riemann solver. For more details on the used mesh see also Fig.~\ref{fig:ANNSL_FFS}.

\subsubsection{\label{subsec:naca} 2D NACA 0012}
\begin{figure}[!htb]
	\begin{minipage}[c]{0.5\linewidth}
		\centering
		\subcaptionbox{Density at $t_{end}=8.0$.}{\includegraphics[width=0.7\linewidth,clip=]{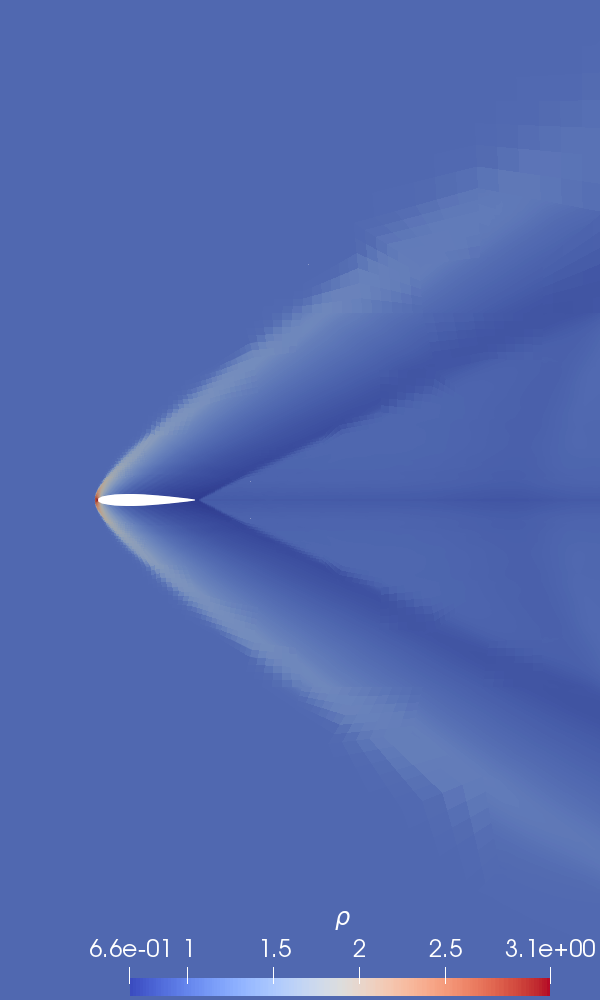}}
	\end{minipage}
	\begin{minipage}[c]{0.5\linewidth}
		\centering
		\subcaptionbox{Unstructured mesh used for the computation.}{\includegraphics[width=0.7\linewidth,clip=]{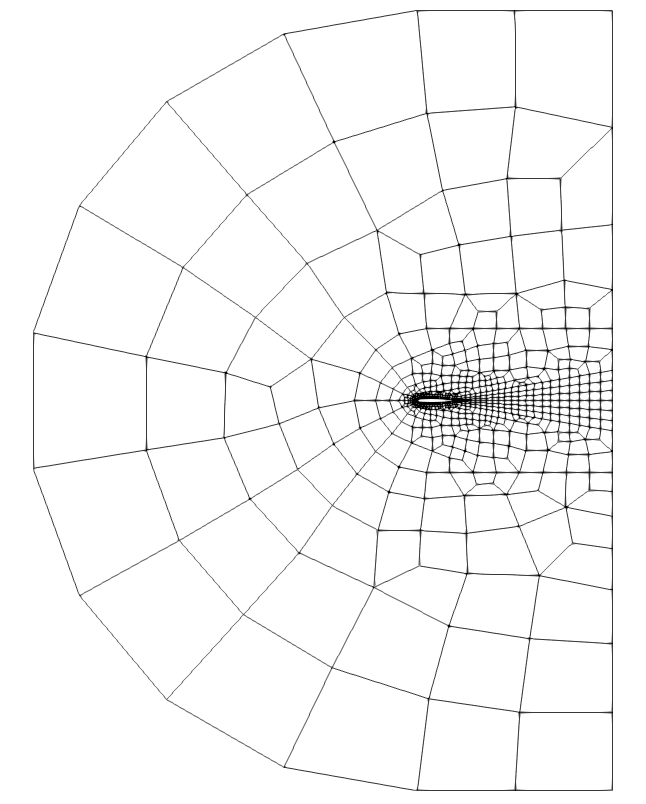}}
	\end{minipage}
	\caption{Density (left) and mesh (right) of the NACA0012 profile. \label{fig:naca0012_density}}
\end{figure}

As a final testcase, we consider the flow over a NACA 0012 airfoil at $Ma=2.0, AoA=0^{\circ}$ on an unstructured grid, depicted in Fig.~\ref{fig:naca0012_density}. This setup results in a detached, curved bow shock and a tail shock. The solution is advanced in time from an initial freestream condition until a steady state is achieved.
The mesh is coarse, without symmetries and does not take a priori knowledge of the shock front position, strength and shape into account, resulting in a deliberately "bad" mesh for this task. This setup was generated to test the shock indicators in a more realistic setting for high order methods than in the previous test cases,  where the shock location might be uncertain or rapidly moving, grid refinement is not available or too costly, and large grid elements are the norm. Fig.~\ref{fig:naca0012_density} shows the computational grid and the flow solution.

\subsection{\label{subsec:n5}ANNSI Results for $\mathcal{N}=5$}
We first consider the results for the $\mathcal{N}=5$ computations and networks, and discuss selected testcases from Sec.~\ref{subsec:testcases} and contrast the results from ANNSI and the troubled cell indicators from Sec.~\ref{subsec:indicators}. It is important to underline that while for the modal and jump indicators 	$\mathcal{I}_\text{modal}$ and $\mathcal{I}_\text{jump}$ the respective parameters defining the thresholds need to be set by the user and possibly adapted for each situation, the ANNSI prediction of the neural network established in Sec.~\ref{sec:NeuralNetworks} is used here 'as is' throughout this chapter, without the need (or possibility) of user intervention or parameter adjustment. \\
In a first step, we apply the ANNSI to a range of test problems to establish its suitability, robustness and accuracy for the intended purpose.

\begin{figure}[!htbp]
	\begin{minipage}[c]{0.03\linewidth}
		\vspace{0.5cm}
		(a) \\
		\vspace{1.7cm}
		(b)
	\end{minipage}
	\begin{minipage}[c]{0.2\linewidth}
		\caption*{Configuration 4}
		\includegraphics[width=\linewidth]{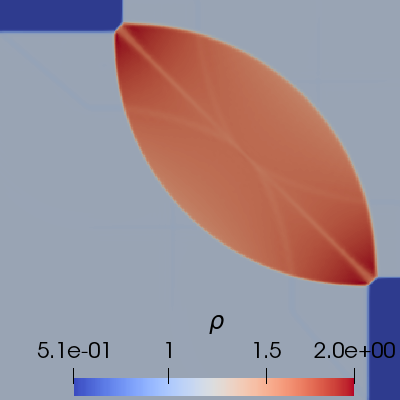}
		\includegraphics[width=\linewidth]{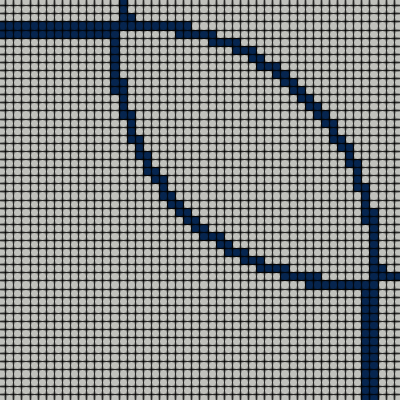}
	\end{minipage}
	\begin{minipage}[c]{0.2\linewidth}
		\caption*{Configuration 6}
		\includegraphics[width=\linewidth]{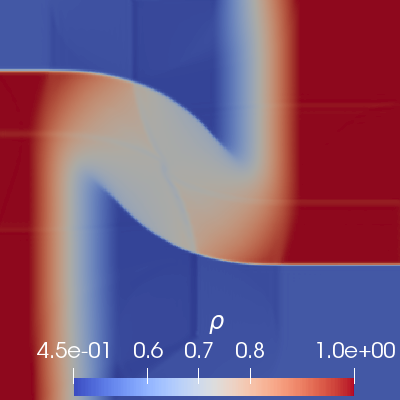}
		\includegraphics[width=\linewidth]{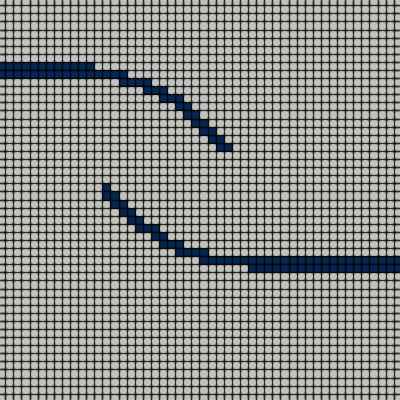}
	\end{minipage}
	\begin{minipage}[c]{0.2\linewidth}
		\caption*{Configuration 12}
		\includegraphics[width=\linewidth]{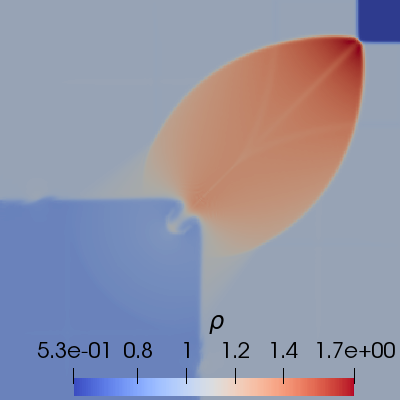}
		\includegraphics[width=\linewidth]{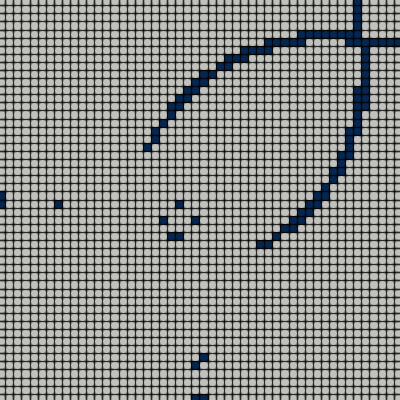}
	\end{minipage}
	\begin{minipage}[c]{0.3\linewidth}
		\caption*{DMR}
		\includegraphics[width=1.5\linewidth]{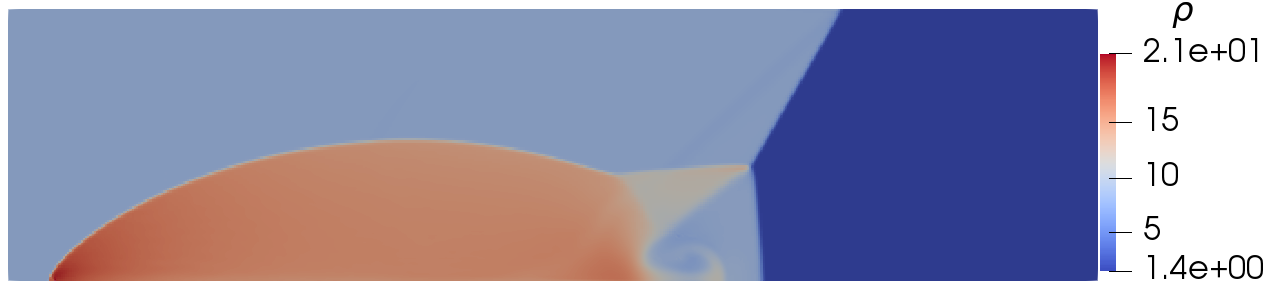}
		\includegraphics[width=1.5\linewidth]{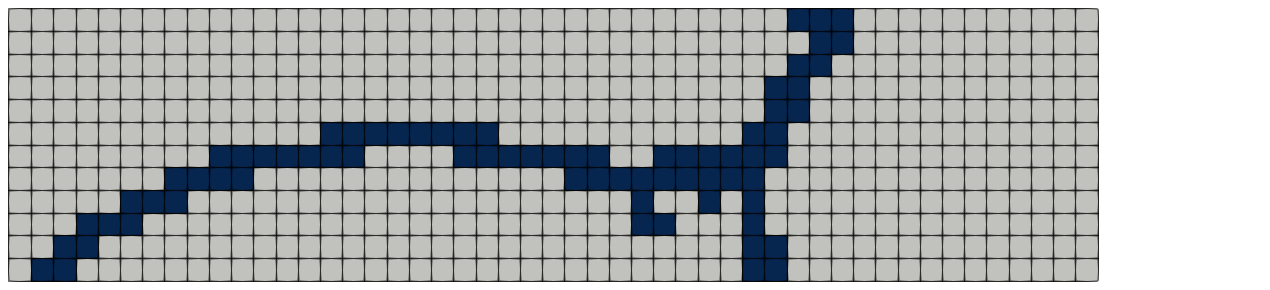}
	\end{minipage}
	\caption{$\mathcal{N}=5$ results: Flagged cells (dark color) of ANNSI for configuration 4 (left), 6 (left of middle) and 12 (right of middle) of the two-dimensional Riemann problems and for the Double Mach Reflection (right). The Riemann problems are calculated on the baseline mesh $M_{RP}$ with $50 \times 50$ mesh and the DMR on the $M_{DMR}$ mesh with $49 \times 12$ elements.\label{fig:comp:n5_annsi}}
\end{figure}
Fig.~\ref{fig:comp:n5_annsi} shows the resulting flow fields in the top row and the flagged grid elements in the bottom row for the Riemann problems and the DMR. All simulations using ANNSI are stable and the resulting flow fields are qualitatively as expected. The flagged cells also coincide with the shock locations, the shock fronts are unbroken and well-captured. For the Riemann problem configuration 12, the contact discontinuity lines are not detected, which is likely due to the roll up of the shear layer and the weak gradients.

Based on these initial findings that the ANNSI method delivers a useful and robust shock detection, we now compare its predictions against the $\mathcal{I}_\text{modal}$ and $\mathcal{I}_\text{jump}$ indicators. As pointed out above, the upper and lower detection thresholds of these indicators are tuneable parameters and no strict quantitative guidelines exist that can help fix them. Instead, the task of the user is to find a setting mainly through trial and error that results in sufficient robustness but do not trigger excessive solution stabilization mechanisms, which can either be computationally expensive or unduly influence smooth regions of the solution - or both. Beyond the academic examples considered here, such parameter variations however are often too costly and cumbersome, and thus often a parameter set based on previous experience is chosen. In order to now compare this approach with the one proposed here, we hand-tuned the $\mathcal{I}_\text{modal}$ and $\mathcal{I}_\text{jump}$ indicators collectively on the four test cases from Fig.~\ref{fig:comp:n5_annsi} (Riemann problems 4,6,12 and DMR) under the following conditions: For each of the two indicators, the optimal parameter settings must provide a stable simulation of all four test cases and must not result in excessive flagging of cells outside of the expected regions.  Our approach resulted in the settings listed in Tbl.~\ref{tab:indi_handtuned}.\\

\begin{table}[!htbp]
	\centering
	\begin{tabular}{lll}
		\toprule
		\# & $\mathcal{I}_\text{modal}$ &$\mathcal{I}_\text{jump}$\\
		\midrule
		$\mathcal{I}_\text{upper}$ & -4.5 & 0.012 \\
		$\mathcal{I}_\text{lower}$ & -4.7 & 0.01 \\

		\bottomrule
	\end{tabular}
	\caption{Settings for $\mathcal{I}_\text{modal}$ and $\mathcal{I}_\text{jump}$, collectively tuned for the RP and DMR cases. \label{tab:indi_handtuned}}
\end{table}
\begin{remark}
The approach of hand-tuning these indicators based on the criteria specified above is somewhat unsatisfactory as it is more or less subjective. However, this is commonly used in practical applications with a priori troubled cell indicators for shock capturing and again reflects the challenges of predicting the solution update of the discretized PDE based on the current solution alone.
\end{remark}

\begin{remark}
Tuning the parameters for each of the cases individually would likely have lead to improved results in some cases. However, our goal here is not to establish the superiority of any of the indicators over the other ones but to compare them qualitatively in a typical setting.
\end{remark}

\begin{figure}[!htbp]
	\begin{minipage}[c]{0.03\linewidth}
		\vspace{2.0cm}
		(a) \\
		\vspace{3.0cm}
		(b) \\
		\vspace{2.5cm}
		(c) \\
		\vspace{1.8cm}
		(d)
	\end{minipage}
	\begin{minipage}[c]{0.3\linewidth}
		\caption*{ANNSI}
		\includegraphics[width=\linewidth]{pictures/RMP4/50x50/RMP4_ANNDI_N5_State_0000000.250000000_FV_Elems.png}
		\includegraphics[width=\linewidth]{pictures/RMP6/50x50/RMP6_ANNDI_N5_State_0000000.200000000_FV_Elems.png}
		\includegraphics[width=\linewidth]{pictures/RMP12/50x50/RMP12_ANNDI_N5_State_0000000.250000000_FV_Elems.png}
		\includegraphics[width=1.08\linewidth]{pictures/DMR/49x12/DMR_ANNDI_N5_State_0000000.200000000_FV_Elems.png}
	\end{minipage}
	\begin{minipage}[c]{0.3\linewidth}
		\caption*{$\mathcal{I}_\text{jump}$}
		\includegraphics[width=\linewidth]{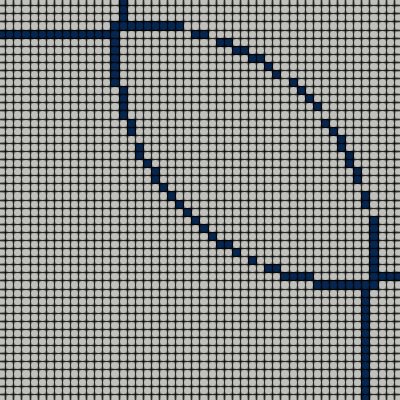}
		\includegraphics[width=\linewidth]{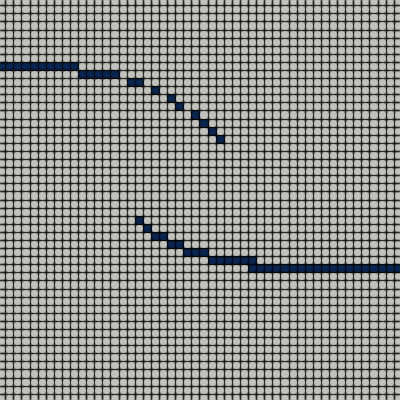}
		\includegraphics[width=\linewidth]{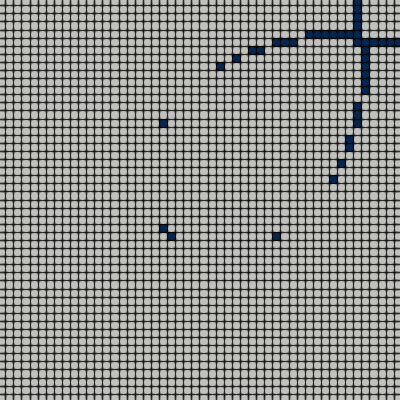}
		\includegraphics[width=0.94\linewidth]{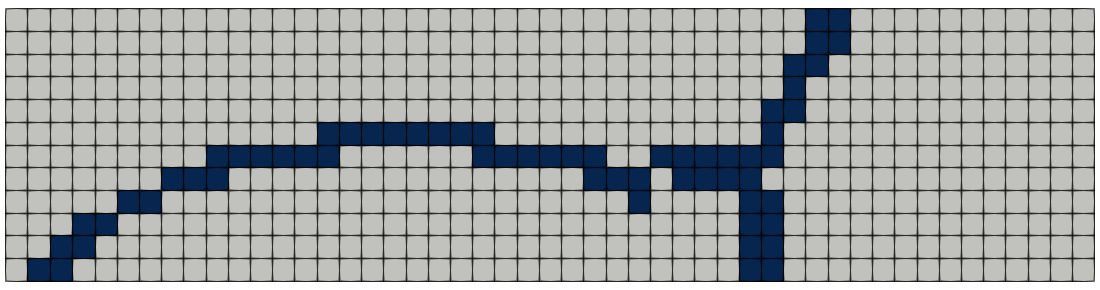}
	\end{minipage}
	\begin{minipage}[c]{0.3\linewidth}
		\caption*{$\mathcal{I}_\text{modal}$}
		\includegraphics[width=\linewidth]{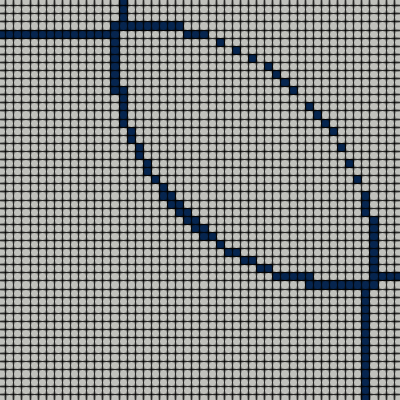}
		\includegraphics[width=\linewidth]{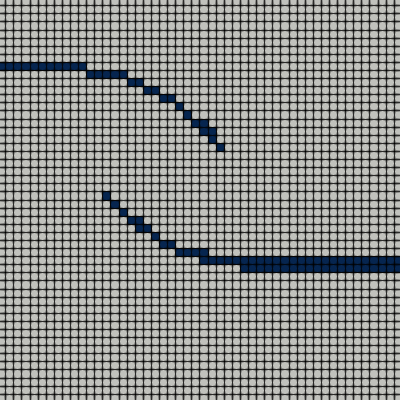}
		\includegraphics[width=\linewidth]{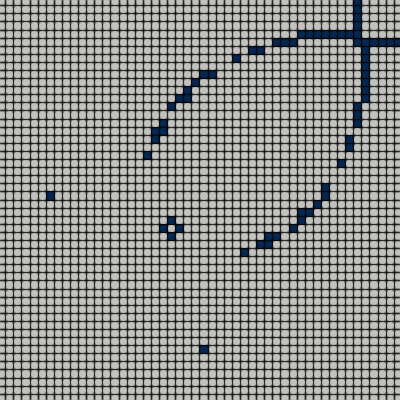}
		\includegraphics[width=0.94\linewidth]{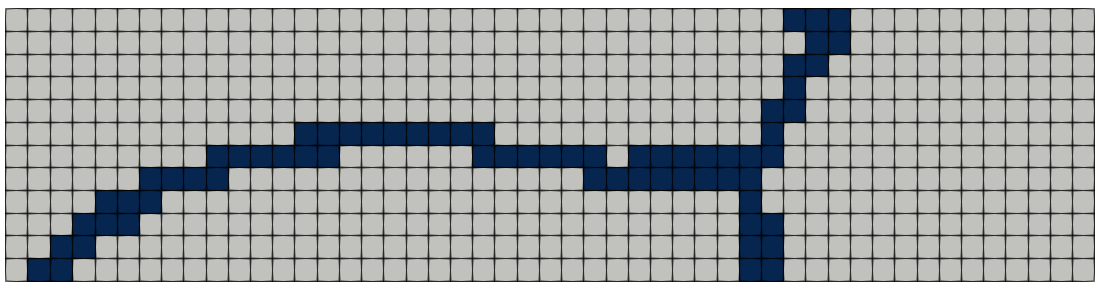}
	\end{minipage}
	\caption{$\mathcal{N}=5$ results: Flagged cells (dark color) of ANNSI (left), the jump (middle) and the modal indicator (right) for configuration 4 (a), 6 (b) and 12 (c) of the two-dimensional Riemann problems and for the Double Mach Reflection (d). The Riemann problems are calculated on the baseline mesh $M_{RP}$ with $50 \times 50$ mesh and the DMR on the $M_{DMR}$ mesh with $49 \times 12$ elements.\label{fig:comp:n5}}
\end{figure}
\begin{figure}[!htbp]
	\begin{minipage}[c]{0.3\linewidth}
		\caption*{ANNSI}
		\includegraphics[width=\linewidth]{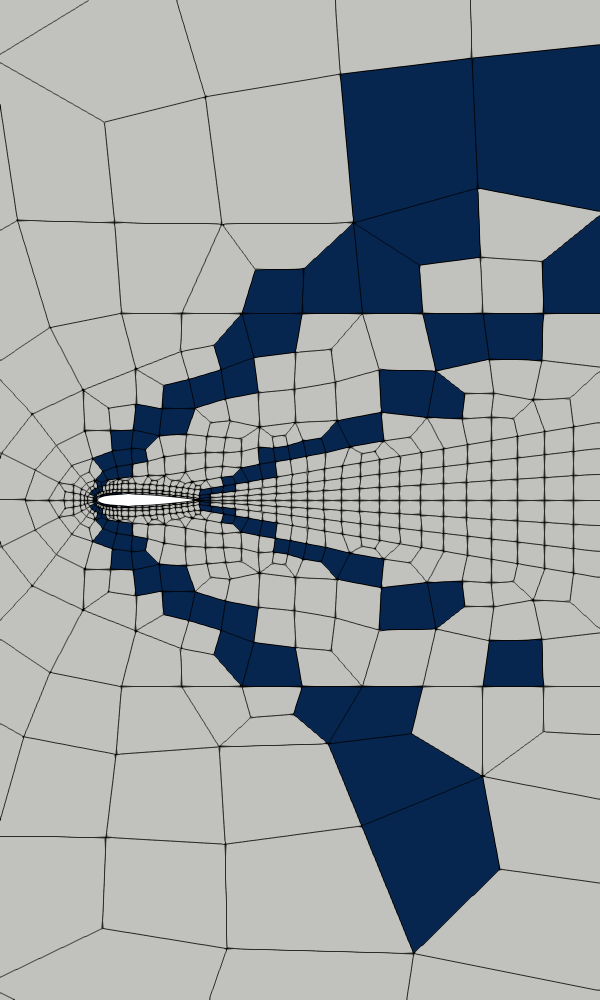}
	\end{minipage}
	\begin{minipage}[c]{0.3\linewidth}
		\caption*{$\mathcal{I}_\text{jump}$}
		\includegraphics[width=\linewidth]{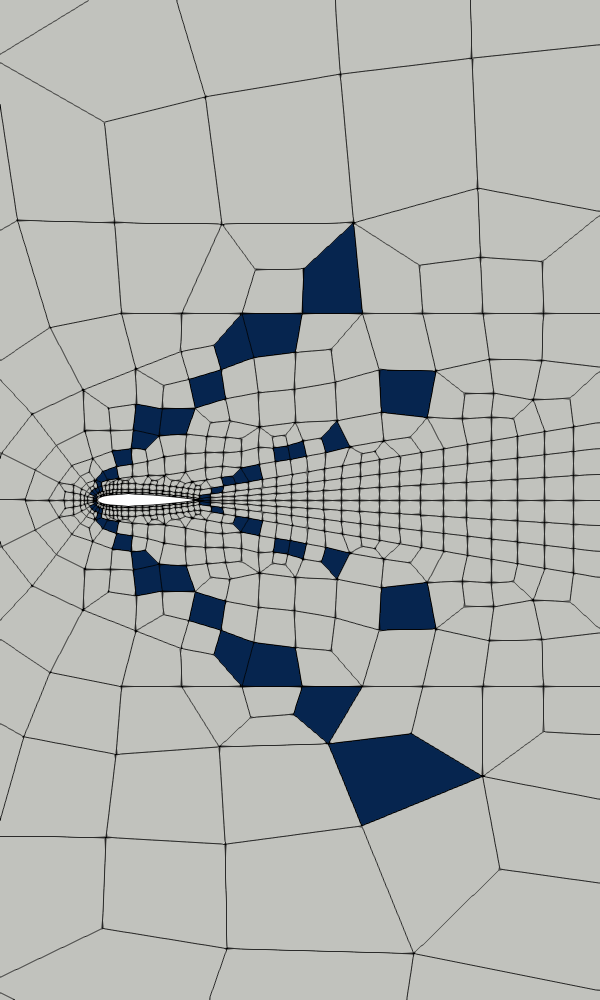}
	\end{minipage}
	\begin{minipage}[c]{0.3\linewidth}
		\caption*{$\mathcal{I}_\text{modal}$}
		\includegraphics[width=\linewidth]{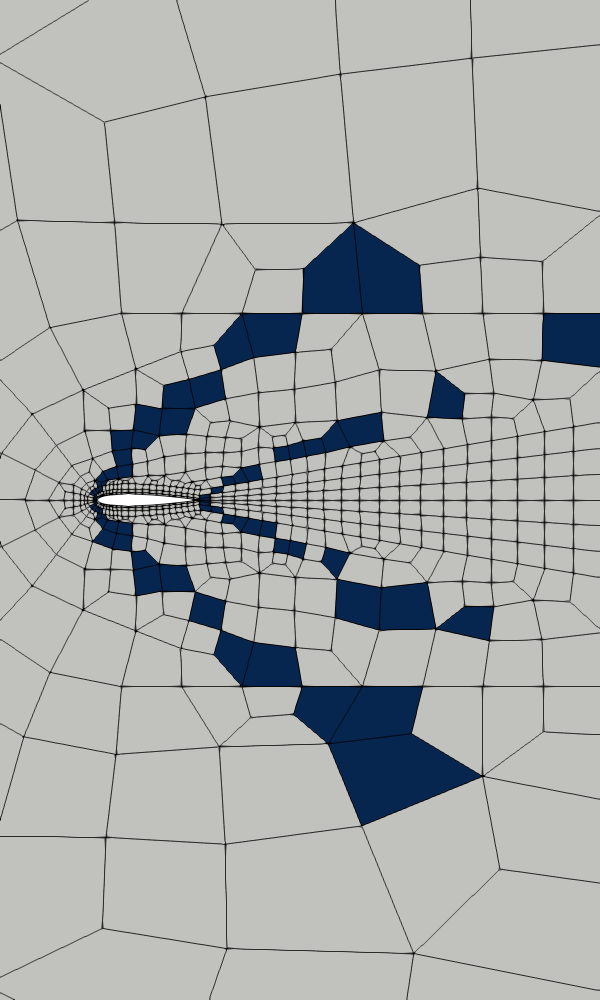}
	\end{minipage}
	\caption{$\mathcal{N}=5$ results: Flagged cells (dark color) of ANNSI (left), the jump (middle) and the modal indicator (right) for $Ma=2.0$ NACA0012 flow.\label{fig:comp:n5_naca}}
\end{figure}

We compare these results with the ANNSI cases in Fig.~\ref{fig:comp:n5}, where we have repeated the relevant results from Fig.~\ref{fig:comp:n5_annsi} as the first column. All three indicators yield comparable results, flagging shock locations but not expansion waves. The ANNSI predictions typically flag two adjacent cells across a shock, while in particular the jump indicator shows a more refined picture. On the other hand, the shock fronts detected by ANNSI are continuous and unbroken in contrast to the  $\mathcal{I}_\text{modal}$ and $\mathcal{I}_\text{jump}$ results, where "holes" in the shock fronts appear.
Applying the same indicators to the NACA0012 case leads to a stable simulation in all cases; the results for the flagged cells are shown in Fig.~\ref{fig:comp:n5_naca}. The ANNSI driven method performs well on this suboptimal grid, yielding continuous shock fronts and a more conservative estimate than the other two. We note that as
the computational grid is not symmetric w.r.t the $x-$axis, so we also do not expect a
necessarily symmetric flagging of the shock fronts.

The results presented in this section can thus be summarized as: The developed ANNSI method is robust, suitable for structured and unstructured grids and of comparable accuracy to other well-tuned indicators.

\subsection{\label{subsec:n9} ANNSI Results for $\mathcal{N}=9$}
In this section we apply the ANNSI network to the same 2D problems as in the previous section.  We perform the calculations with a polynomial degree of $\mathcal{N}=9$ to show the invariance of the performance with respect to the polynomial degree of the ansatz functions.
\begin{figure}[!htbp]
	\begin{minipage}[c]{0.03\linewidth}
		\vspace{0.5cm}
		(a) \\
		\vspace{1.7cm}
		(b)
	\end{minipage}
	\begin{minipage}[c]{0.2\linewidth}
		\caption*{Configuration 4}
		\includegraphics[width=\linewidth]{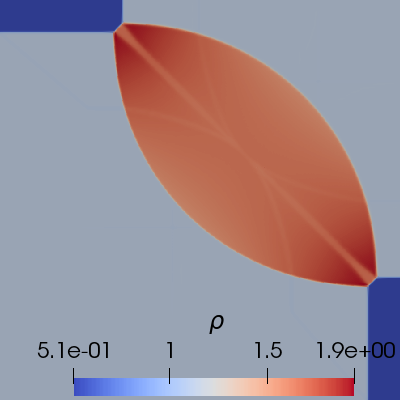}
		\includegraphics[width=\linewidth]{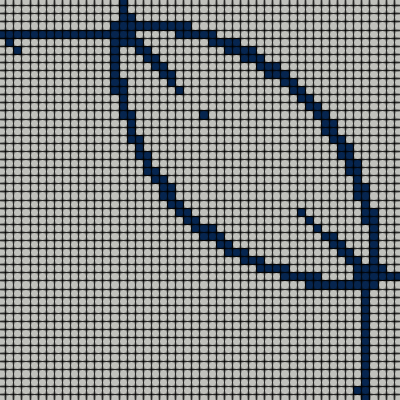}
	\end{minipage}
	\begin{minipage}[c]{0.2\linewidth}
		\caption*{Configuration 6}
		\includegraphics[width=\linewidth]{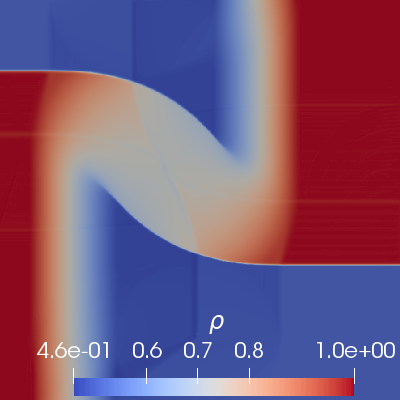}
		\includegraphics[width=\linewidth]{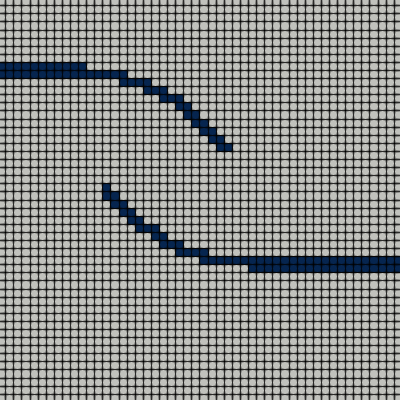}
	\end{minipage}
	\begin{minipage}[c]{0.2\linewidth}
		\caption*{Configuration 12}
		\includegraphics[width=\linewidth]{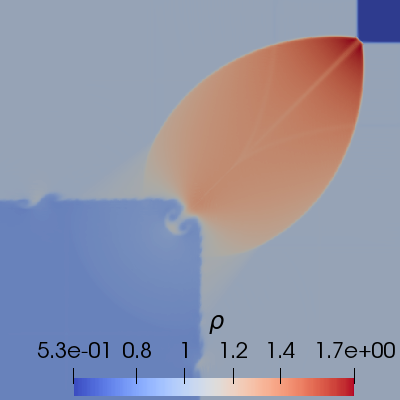}
		\includegraphics[width=\linewidth]{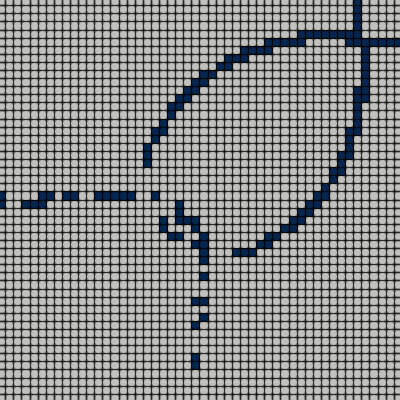}
	\end{minipage}
	\begin{minipage}[c]{0.3\linewidth}
		\caption*{DMR}
		\includegraphics[width=\linewidth]{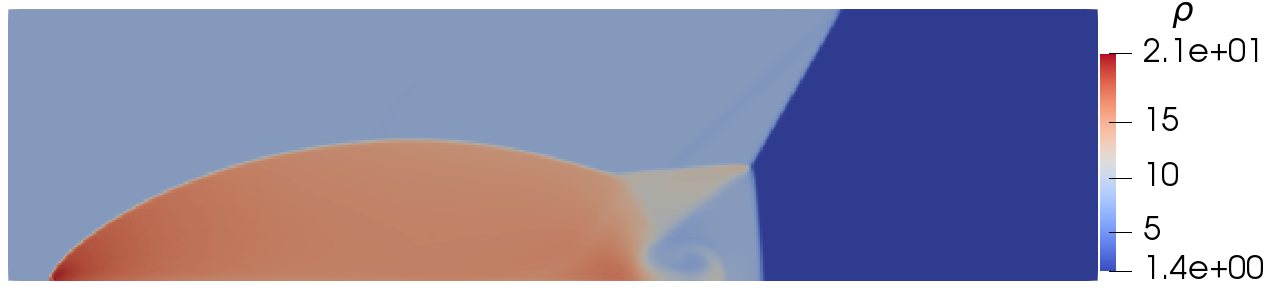}
		\includegraphics[width=\linewidth]{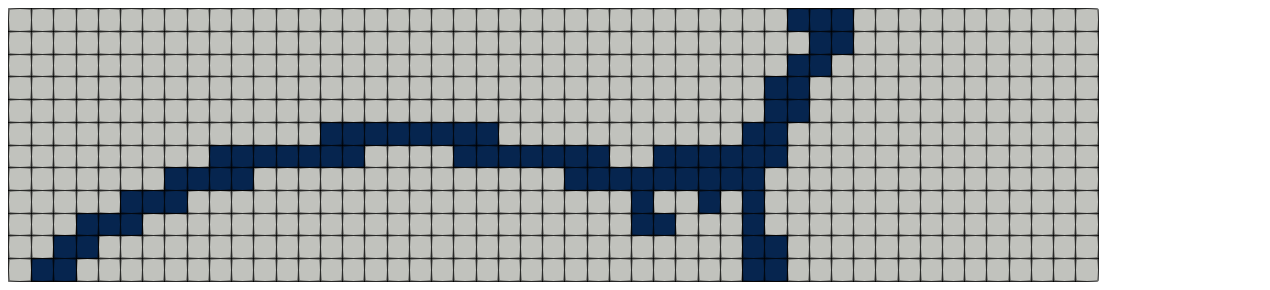}
	\end{minipage}
	\caption{$\mathcal{N}=9$ results: Flagged cells (dark color) of ANNSI for configuration 4 (a), 6 (b) and 12 (c) of the two-dimensional Riemann problems and for the Double Mach Reflection (d). The Riemann problems are calculated on a mesh with half the resolution in each direction ($25 \times 25$ elements) and the DMR on a mesh with $24 \times 6$ elements.\label{fig:comp:n9}}
\end{figure}

In Fig.~\ref{fig:comp:n9} we repeat the calculations from Fig.~\ref{fig:comp:n5_annsi} with $\mathcal{N}=9$ instead of $\mathcal{N}=5$. A very similar behavior can be observed: again all simulations using the ANNSI are stable and the resulting flow fields are qualitatively as expected. The shock localization is correctly recognized and captured with the finite volume sub-cell scheme. Differences can be seen at the contact discontinuities in RP12 which are partially detected by the indicator. We also compared the ANNSI results against those for the selected problems in a repetition of Fig.~\ref{fig:comp:n5}. We found the results from $\mathcal{N}=5$ transferable to the $\mathcal{N}=9$ case.

\subsection{\label{subsec:grid_reso} Influence of Grid Resolution}
Having established the properties of the ANNSI approach on the baseline grids, we now investigate the influence of grid resolution changes by either halving or doubling the grid cells for the $\mathcal{N}=5$ case  per direction. Since the training data for ANNSI was constructed in reference space, we expect its results to be robust against grid resolution changes. For $\mathcal{I}_\text{jump}$ and $\mathcal{I}_\text{modal}$, we again use the found parameter settings from Sec.~\ref{subsec:n5}.

\begin{figure}[!htbp]
	\begin{minipage}[c]{0.03\linewidth}
		\vspace{2.0cm}
		(a) \\
		\vspace{3.0cm}
		(b) \\
		\vspace{2.5cm}
		(c) \\
		\vspace{1.8cm}
		(d)
	\end{minipage}
	\begin{minipage}[c]{0.3\linewidth}
		\caption*{ANNSI}
		\includegraphics[width=\linewidth]{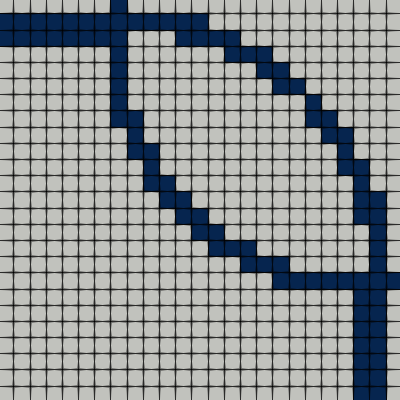}
		\includegraphics[width=\linewidth]{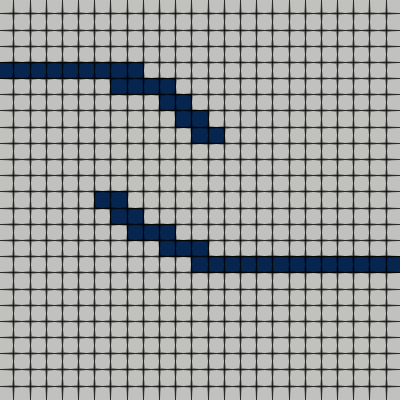}
		\includegraphics[width=\linewidth]{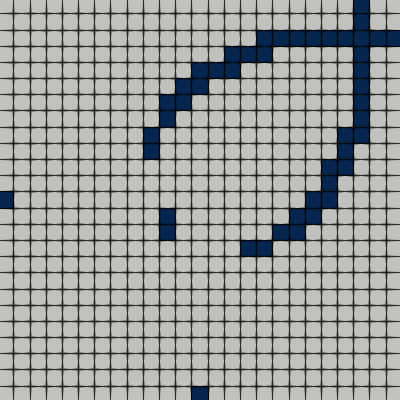}
		\includegraphics[width=\linewidth]{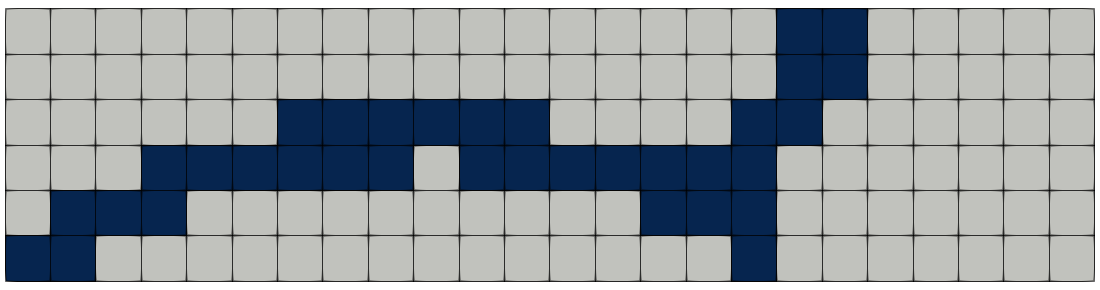}
	\end{minipage}
	\begin{minipage}[c]{0.3\linewidth}
		\caption*{$\mathcal{I}_\text{jump}$}
		\includegraphics[width=\linewidth]{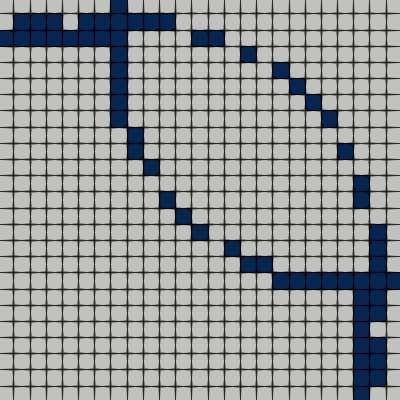}
		\includegraphics[width=\linewidth]{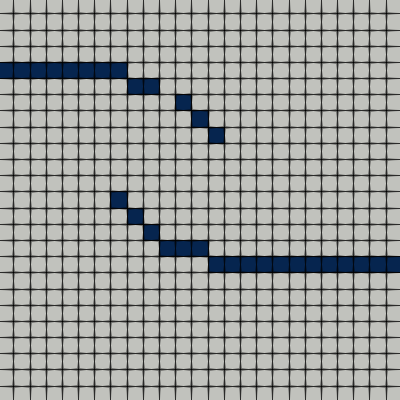}
		\includegraphics[width=\linewidth]{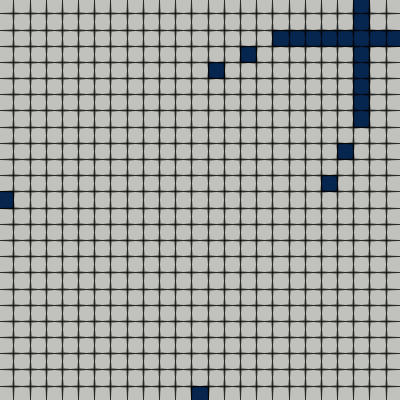}
		\includegraphics[width=\linewidth]{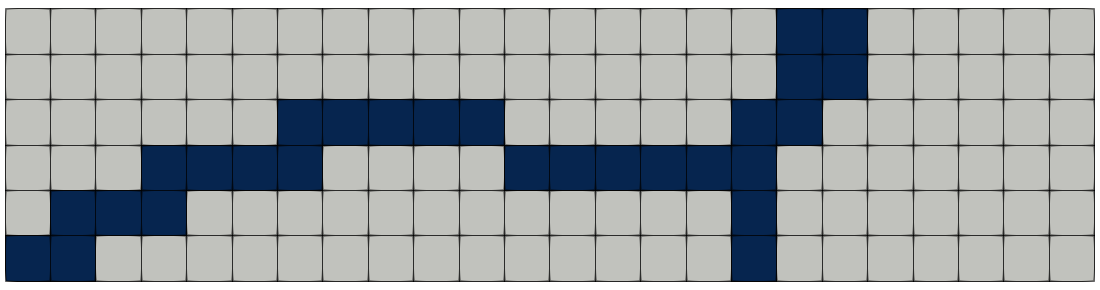}
	\end{minipage}
	\begin{minipage}[c]{0.3\linewidth}
		\caption*{$\mathcal{I}_\text{modal}$}
		\includegraphics[width=\linewidth]{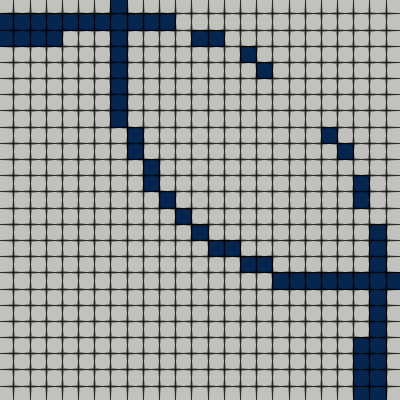}
		\includegraphics[width=\linewidth]{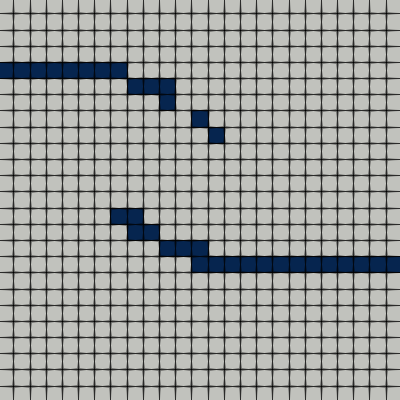}
		\includegraphics[width=\linewidth]{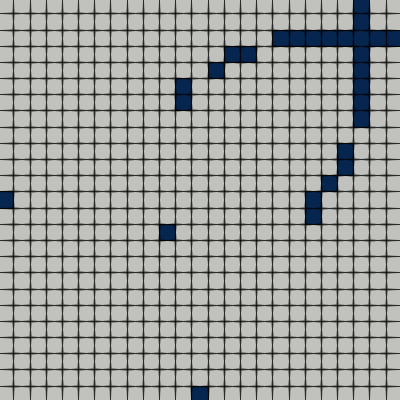}
		\includegraphics[width=\linewidth]{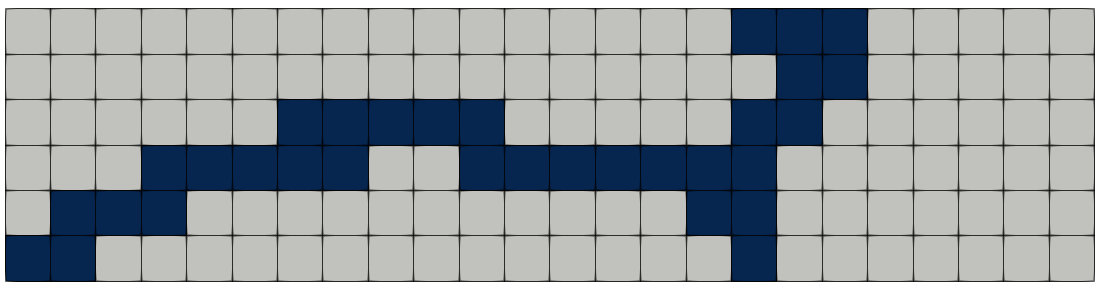}
	\end{minipage}
	\caption{$\mathcal{N}=5$ results, coarse mesh: Flagged cells (dark color) of ANNSI (left), the jump (middle) and the modal indicator (right) for configuration 4 (a), 6 (b) and 12 (c) of the two-dimensional Riemann problems and for the Double Mach Reflection (d). The Riemann problems are calculated on a mesh with half the resolution in each direction ($25 \times 25$ elements) and the DMR on a mesh with $24 \times 6$ elements.\label{fig:comp:n5_coarse}}
\end{figure}

\begin{figure}[!htbp]
	\begin{minipage}[c]{0.03\linewidth}
		\vspace{2.0cm}
		(a) \\
		\vspace{3.0cm}
		(b) \\
		\vspace{2.5cm}
		(c) \\
		\vspace{1.8cm}
		(d)
	\end{minipage}
	\begin{minipage}[c]{0.3\linewidth}
		\caption*{ANNSI}
		\includegraphics[width=\linewidth]{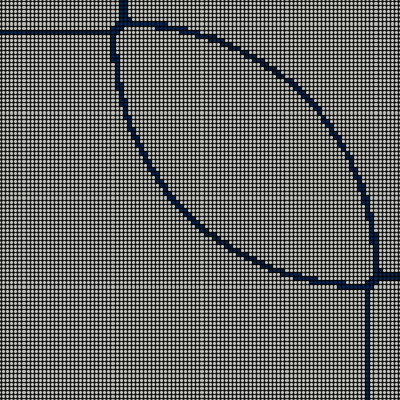}
		\includegraphics[width=\linewidth]{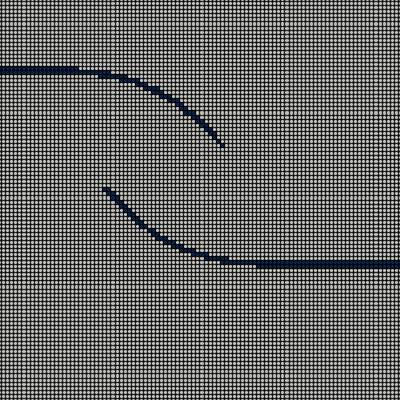}
		\includegraphics[width=\linewidth]{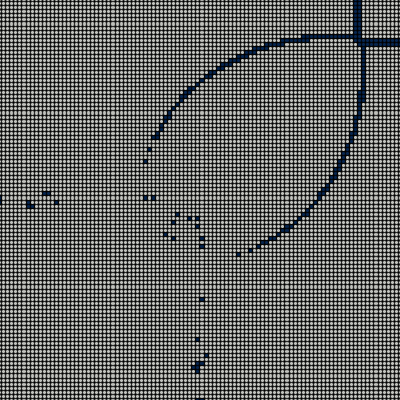}
		\includegraphics[width=1.135\linewidth]{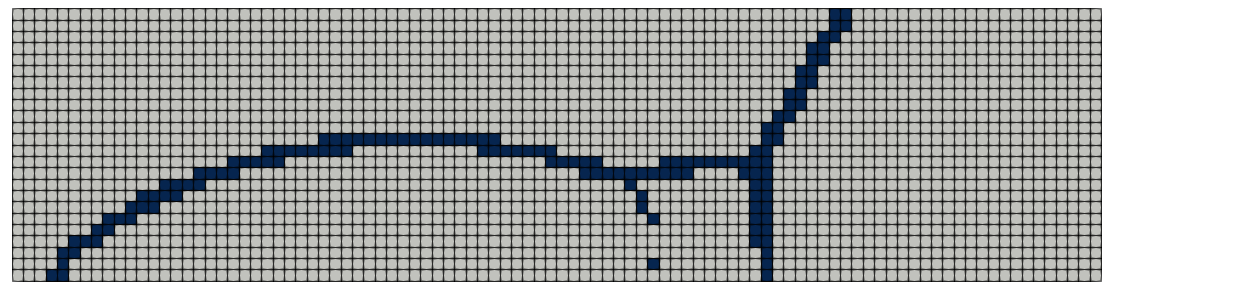}
	\end{minipage}
	\begin{minipage}[c]{0.3\linewidth}
		\caption*{$\mathcal{I}_\text{jump}$}
		\includegraphics[width=\linewidth]{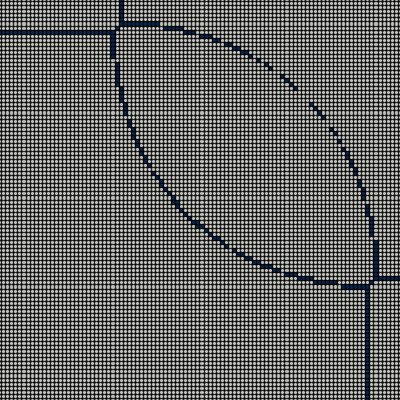}
		\includegraphics[width=\linewidth]{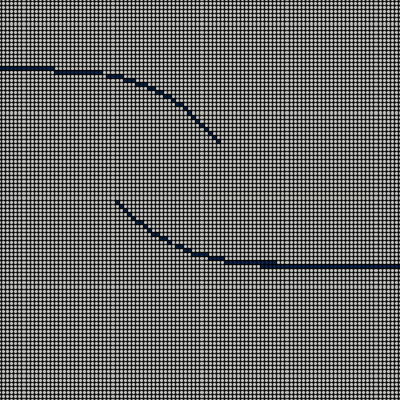}
		\includegraphics[width=\linewidth]{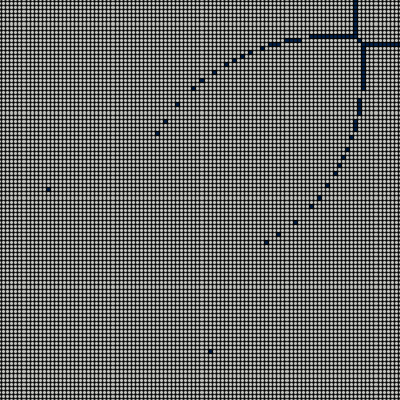}
		\vspace{0.35cm}
	\end{minipage}
	\begin{minipage}[c]{0.3\linewidth}
		\caption*{$\mathcal{I}_\text{modal}$}
		\includegraphics[width=\linewidth]{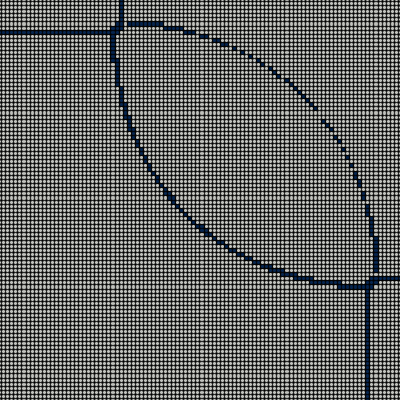}
		\includegraphics[width=\linewidth]{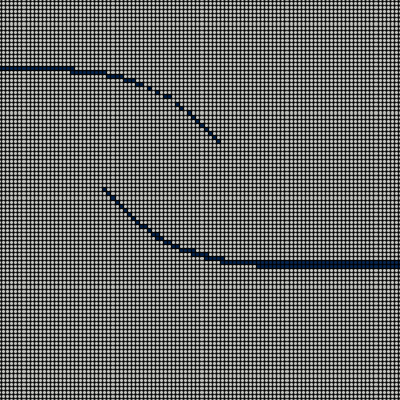}
		\includegraphics[width=\linewidth]{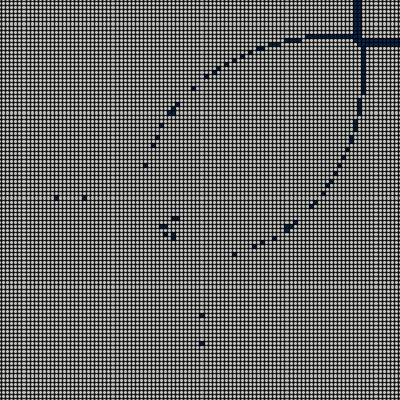}
		\includegraphics[width=\linewidth]{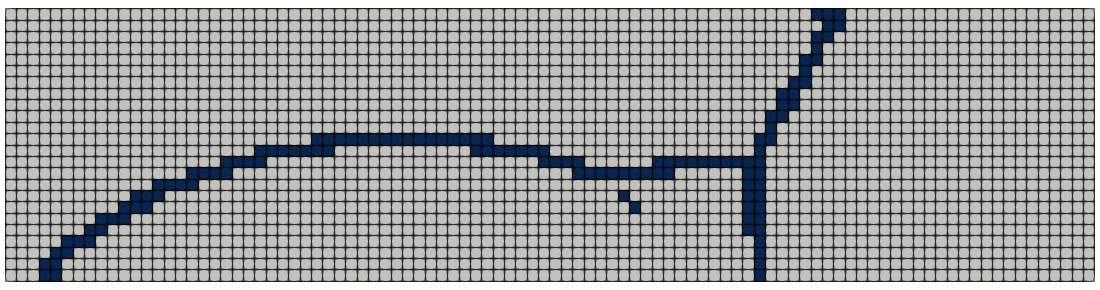}
	\end{minipage}
	\caption{$\mathcal{N}=5$ results, fine mesh: Flagged cells (dark color) of ANNSI (left), the jump (middle) and the modal indicator (right) for configuration 4 (a), 6 (b) and 12 (c) of the two-dimensional Riemann problems and for the Double Mach Reflection (d). The Riemann problems are calculated on a mesh with twice the resolution in each direction ($100 \times 100$ elements) and the DMR on a mesh with $96 \times 24$ elements.\label{fig:comp:n5_fine}}
\end{figure}

Fig.~\ref{fig:comp:n5_coarse} and~\ref{fig:comp:n5_fine} repeat the investigations for the coarse and fine grids. Missing plots indicate that the simulation was not stable.
The coarse grid results confirm our assumption that the ANNSI approach is rather robust against a resolution drop. Comparing Fig.~\ref{fig:comp:n5_coarse} and  Fig.~\ref{fig:comp:n5}, the shock position and extension are predicted in a consistent manner. On the fine grid in Fig.~\ref{fig:comp:n5_fine}, the ANNSI results converge and become sharper as desirable.
\begin{remark}
Note that it is certainly possible to achieve stable results for all testcases with the $\mathcal{I}_\text{jump}$ and $\mathcal{I}_\text{modal}$ indicators under the given conditions (see e.g.~\cite{sonntagdiss,sonntag2017efficient}), but this would require retuning the threshold parameters. Since our focus here is precisely on investigating these sensitivities with regards to the resolution, we did not adjust the parameters of $\mathcal{I}_\text{jump}$ and $\mathcal{I}_\text{modal}$.
\end{remark}

Summarizing the findings in this section, the developed ANNSI approach results in robust solutions for the problems considered. It compares well against other solution-based indicators in terms of accuracy and does not excessively flag elements outside of shock regions and can deal robustly with the initial solution transients and shock movement. It has shown to be insensitive to grid resolution and consistently flags the same regions / features,  thus working as intended as a discretization independent shock detector, and not as a troubled cell indicator. While results of comparable or possibly better quality for specific cases can be achieved by parameter studies and tuning with other a priori indicators, once trained, ANNSI is parameter-free. It can be argued that the time spent training the neural network should be weighted against the time one might invest in parameter tuning, but with large scale application with huge computational costs in mind, it is certainly more economic to do invest the costs "offline" before the actual computation. This aspect will become even more important in the next section, where we discuss the results of the ANNSL extension of ANNSI, where a localization of the shock front within an element is proposed.

\section{\label{sec:NNshock_local}Results: Shock Localization}
In this section, we illustrate how the information about the shock location inside an element generated by the ANNSL network can be exploited. First, the test cases from the previous section are used to illustrate the capabilities of the shock localization procedure. Next, the ANNSL indicator is applied to a new testcase and the possible exploitation of the shock localization in an adaptive mesh refinement (AMR) framework is illustrated. All ANNSL computations are again started from ab initio, with ANNSL being active in every timestep.
\subsection{Shock Localization\label{subsec:annsl1}}
In the previous section, we have shown that the ANNSI shock detector performs well in terms of robustness, accuracy and insensitivity to spatial resolution changes for a range of classical test problems. Here, we use those applications again to illustrate the capabilities of the ANNSL network of identifying the shock position inside an element, i.e. of creating an inner-element bounding box around the shock front and flagging the solution points adjacent to the front. Since this is particularly useful for high order discretizations  such as $\mathcal{N}=9$, we will focus our investigation on this case, using the network trained in Sec.~\ref{sec:NNTraining} on the data from Sec.~\ref{sec:Training_Localization}.
\begin{figure}
	\centering
	\begin{minipage}[c]{0.3\linewidth}
		\caption*{Configuration 4}
		\includegraphics[width=\linewidth]{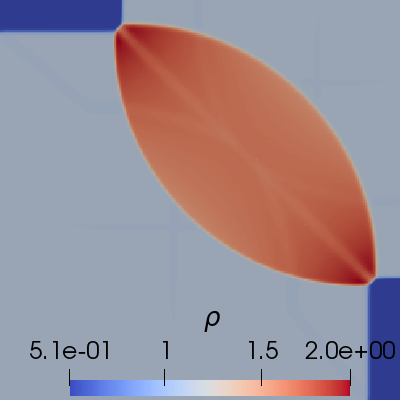}
		\includegraphics[width=\linewidth]{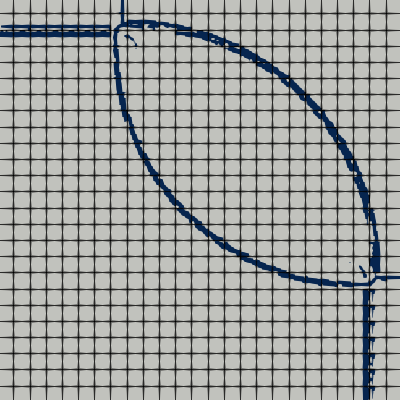}
	\end{minipage}
	\begin{minipage}[c]{0.3\linewidth}
		\caption*{Configuration 6}
		\includegraphics[width=\linewidth]{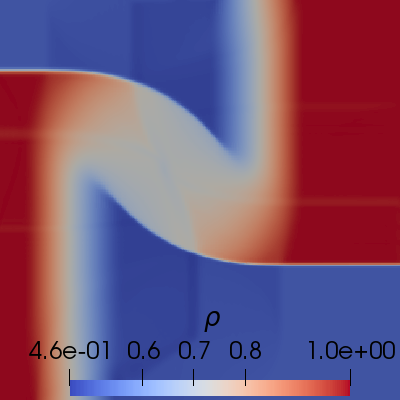}
		\includegraphics[width=\linewidth]{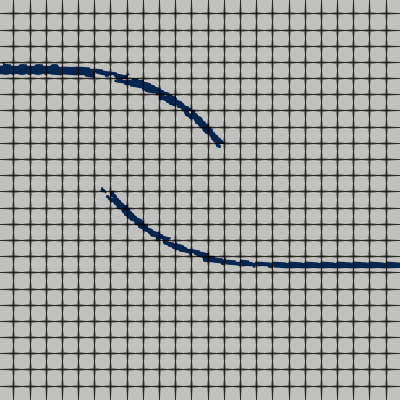}
	\end{minipage}
	\begin{minipage}[c]{0.3\linewidth}
		\caption*{Configuration 12}
		\includegraphics[width=\linewidth]{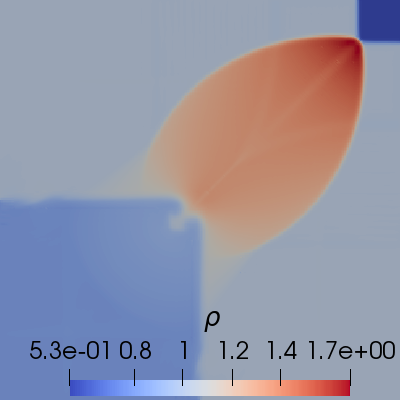}
		\includegraphics[width=\linewidth]{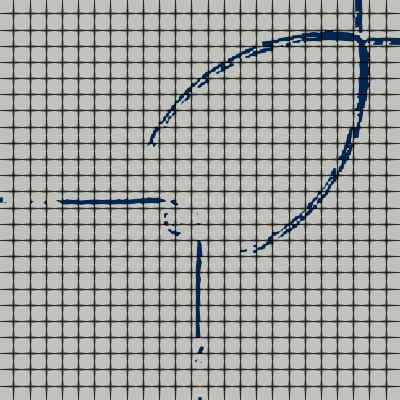}
	\end{minipage}
	\caption{Density (top) and binary shock edge map obtained with the ANNSL network with the underlying grid of size $25 \times 25$ (bottom) of RP configuration 4 (left), 6 (middle) and 12 (right).}\label{fig:ANNSL_RPs}
\end{figure}
\begin{figure}
	\centering
	\begin{minipage}[c]{0.4\linewidth}
		\includegraphics[width=\linewidth]{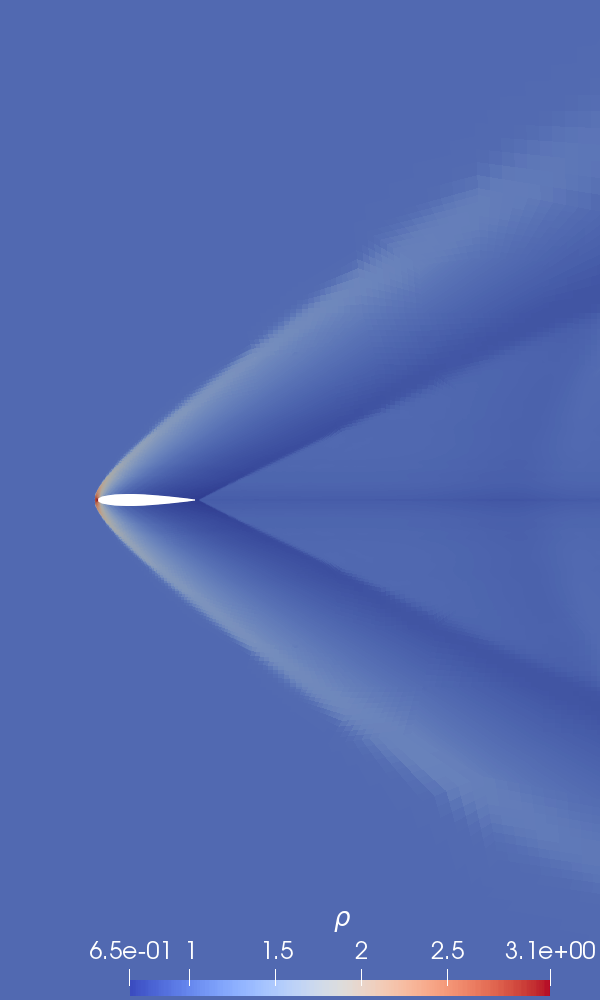}
	\end{minipage}
	\begin{minipage}[c]{0.4\linewidth}
		\includegraphics[width=\linewidth]{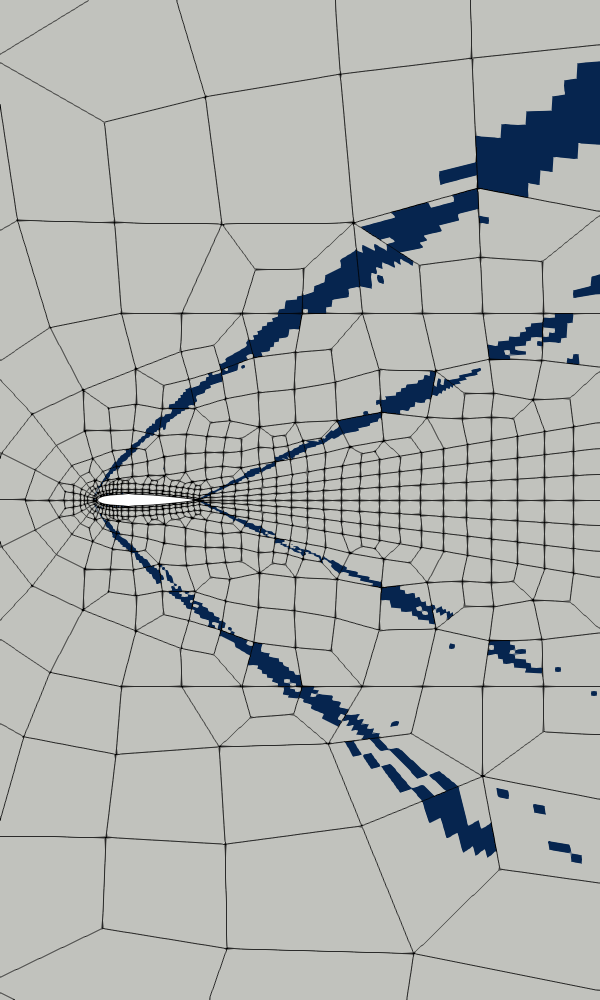}
	\end{minipage}
	\caption{Density (left) and binary shock edge map obtained with the ANNSL network with the underlying grid (right) of the NACA 0012 profile.}\label{fig:ANNSL_NACA}
\end{figure}
\begin{figure}
	\begin{minipage}[c]{0.44\linewidth}
		\includegraphics[trim={0 0 6cm 0},clip,width=\linewidth]{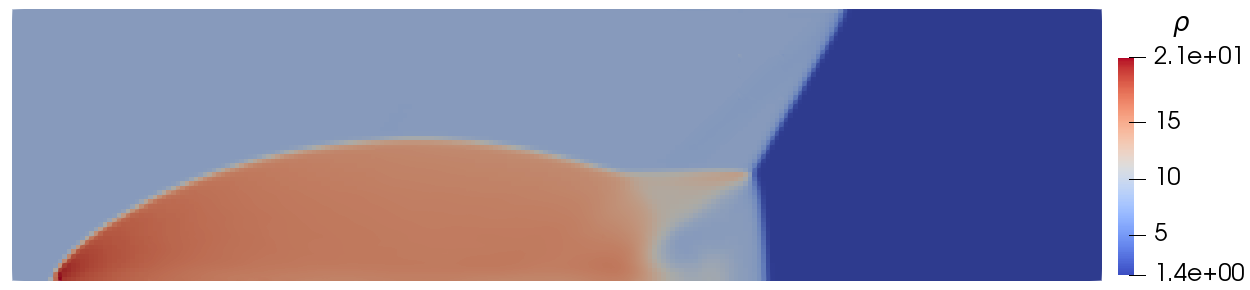}
		\includegraphics[trim={0 0 6cm 0},clip,width=\linewidth]{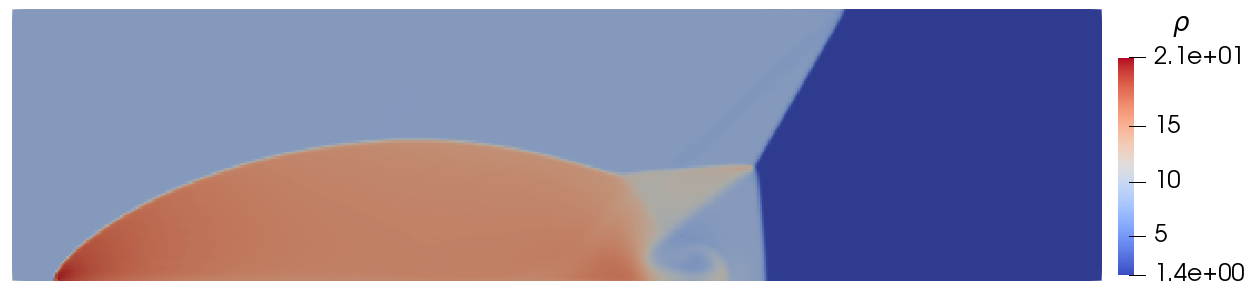}
		\includegraphics[trim={0 0 6cm 0},clip,width=\linewidth]{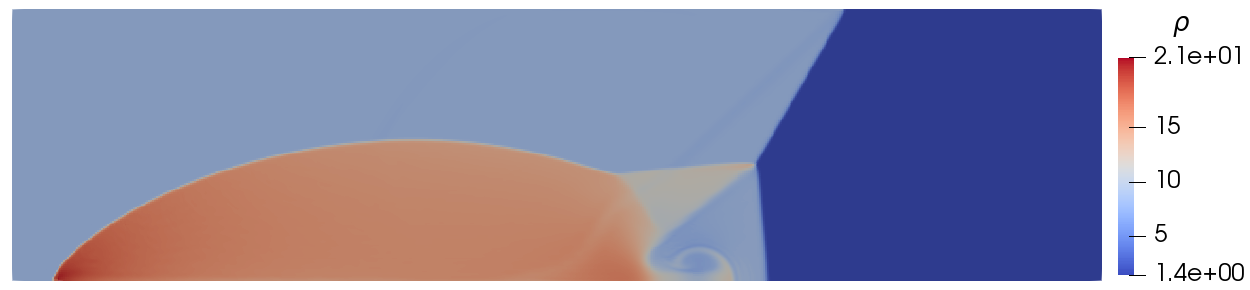}
	\end{minipage}
	\begin{minipage}[c]{0.1\linewidth}
	\begin{tikzpicture}[line cap=round,line join=round,x=1.15cm,y=1.15cm,axis/.style={->}]
	    \node at (0 ,-0.4) {\includegraphics[scale=0.10]{./pictures/bar_color_red.png}};
	    \node at (0.0,1.6 ) {$~~~\rho~~~$};
	    \node at (0.5,1.2 ) {$21$};
	    \node at (0.5,0.4 ) {$15$};
	    \node at (0.5,-0.4 ) {$10$};
	    \node at (0.5,-1.2 ) {$5.0$};
	    \node at (0.5,-2.0 ) {$1.4$};
	\end{tikzpicture}
	\end{minipage}
	\hspace{0.2cm}
	\begin{minipage}[c]{0.51\linewidth}
		\includegraphics[width=\linewidth]{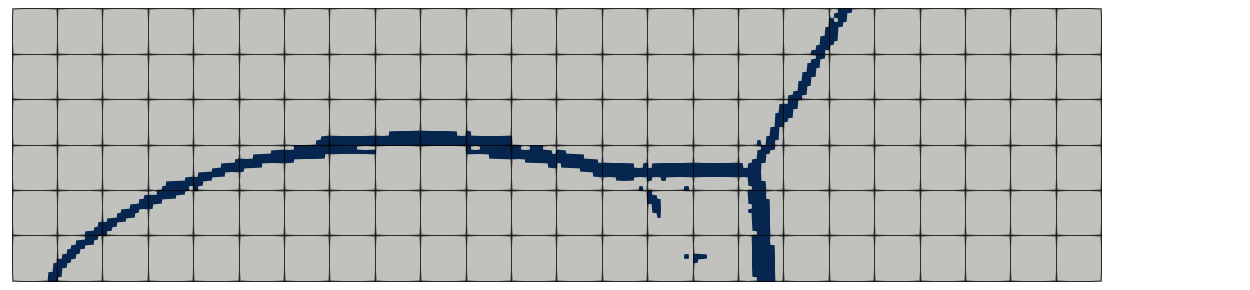}
		\includegraphics[width=\linewidth]{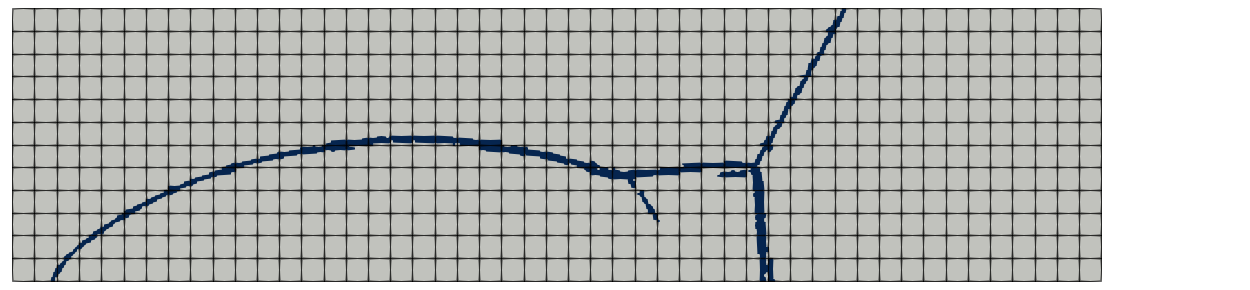}
		\includegraphics[width=\linewidth]{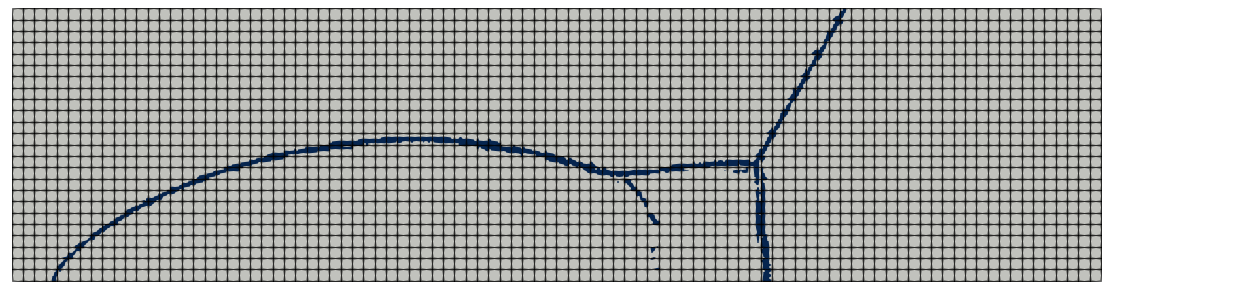}
	\end{minipage}
	\caption{Density (left) and binary shock edge map obtained with the ANNSL network with the underlying grid (right) of the DMR computed on meshes of size $24\times6$ (upper), $49 \times 12$ (middle) and $96 \times 24$ (lower).}\label{fig:ANNSL_DMR}
\end{figure}
Fig.~\ref{fig:ANNSL_RPs} and Fig.~\ref{fig:ANNSL_NACA} apply the ANNSL network to the established test cases. The binary edge map is shown alongside the flow solution, with dark colors indicating that a pixel / solution point is considered to be in the direct vicinity of or on the shock front. One can observe that with the proposed indicator, a very sharp identification of discontinuities is possible. The identification is generally continuous within a cell (no holes or bumps in the shock fronts) and consistent across element boundaries - neighboring element "agree" in the identification of the fronts. All of these properties are desirable in edge detection, however, many more simple edge detection algorithms struggle with one or more of them~\cite{xie2015}. For the intentionally "bad" grid in Fig.~\ref{fig:ANNSL_NACA}, the ANNSL prediction also works very well on coarse cells. The infrequently occurring 'two-stripes' patterns (see e.g. lower region of the bow shock) are due to the discrete resolution of a shock: on the discrete level, it consists of two regions with kinks at the foot and head of the shock and an almost linear part in between. As the linear part is not distinguishable from other smooth solutions, it is not recognized as a shock, and only the shock edges are recognized as parts of a shock.\\
In Fig.~\ref{fig:ANNSL_DMR}, we analyze the convergence behaviour of the shock map with decreasing element size. We observe that the shock localization is stable with respect to the grid size, and the edge map predictions become successively more refined. This feature is an important prerequisite for successful h/p/r-adaption strategies to improve accuracy at the shock. In the next section, we discuss a possibility for taking advantage of this knowledge during the simulation.

\subsection{Neural Network informed H-Adaptation}
\begin{figure}
	\centering
	\begin{minipage}[c]{0.45\linewidth}
		\includegraphics[trim={8cm 1cm 1cm 0},clip,width=\linewidth]{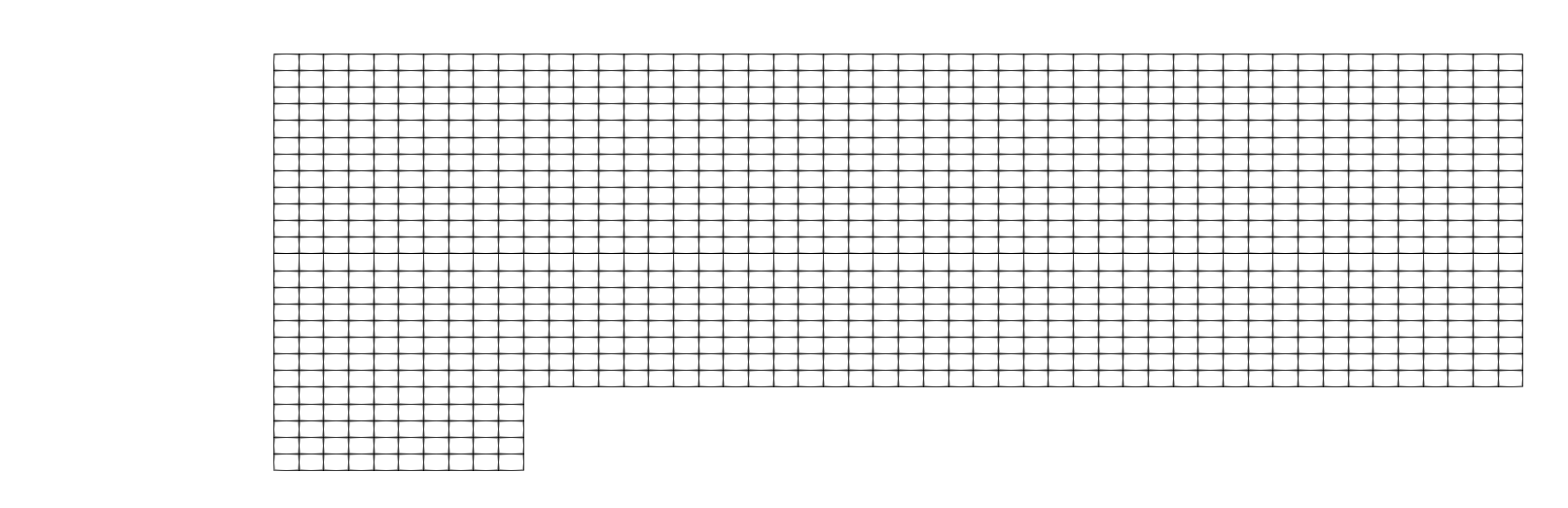}
		\includegraphics[trim={8cm 1cm 1cm 0},clip,width=\linewidth]{pictures/FFS/FFS_Density.png}
		\includegraphics[trim={8cm 1cm 1cm 0},clip,width=\linewidth]{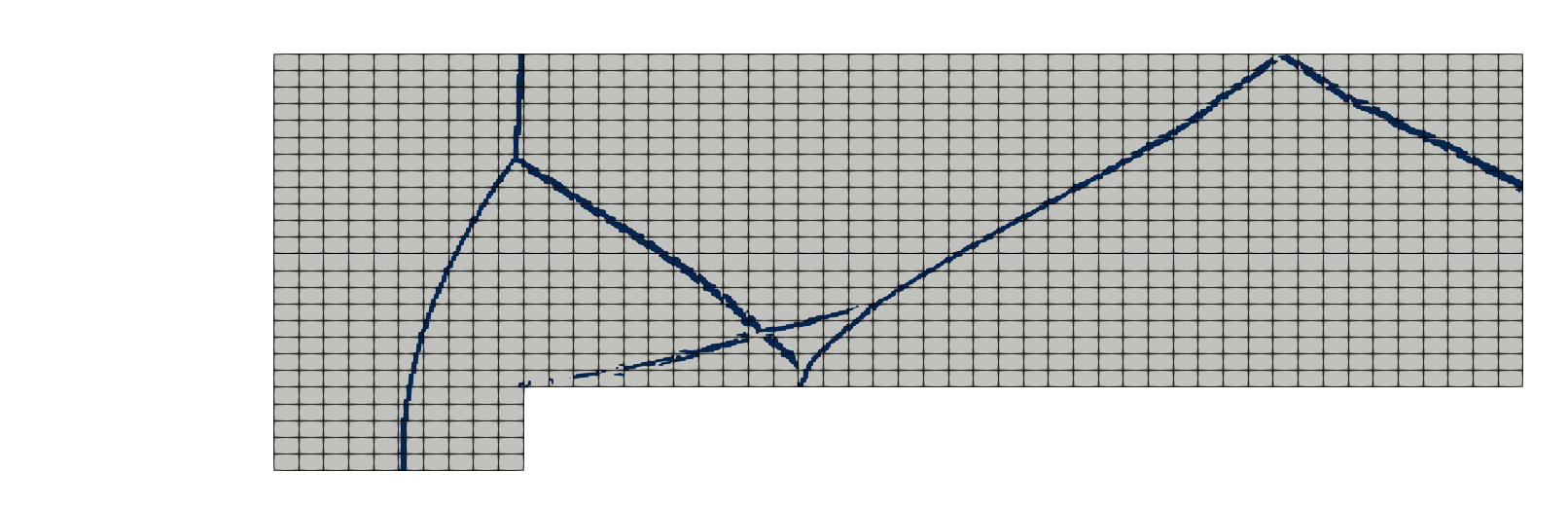}
	\end{minipage}
	\begin{minipage}[c]{0.45\linewidth}
		\includegraphics[trim={8cm 1cm 1cm 0},clip,width=\linewidth]{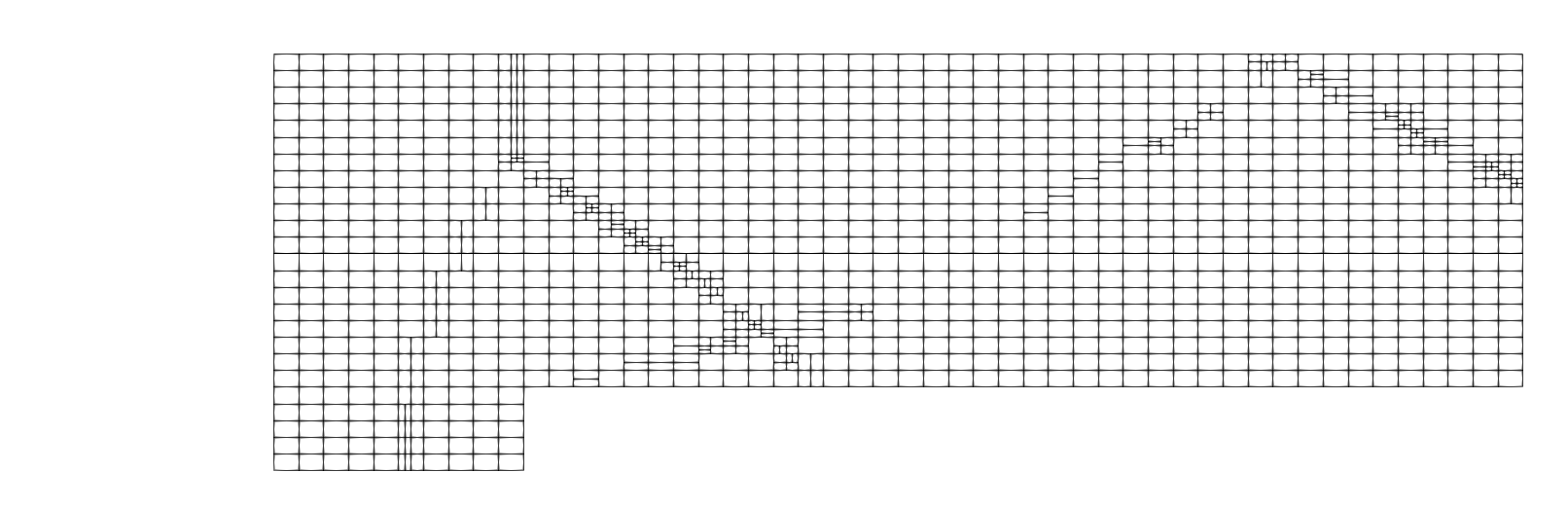}
		\includegraphics[trim={8cm 1cm 1cm 0},clip,width=\linewidth]{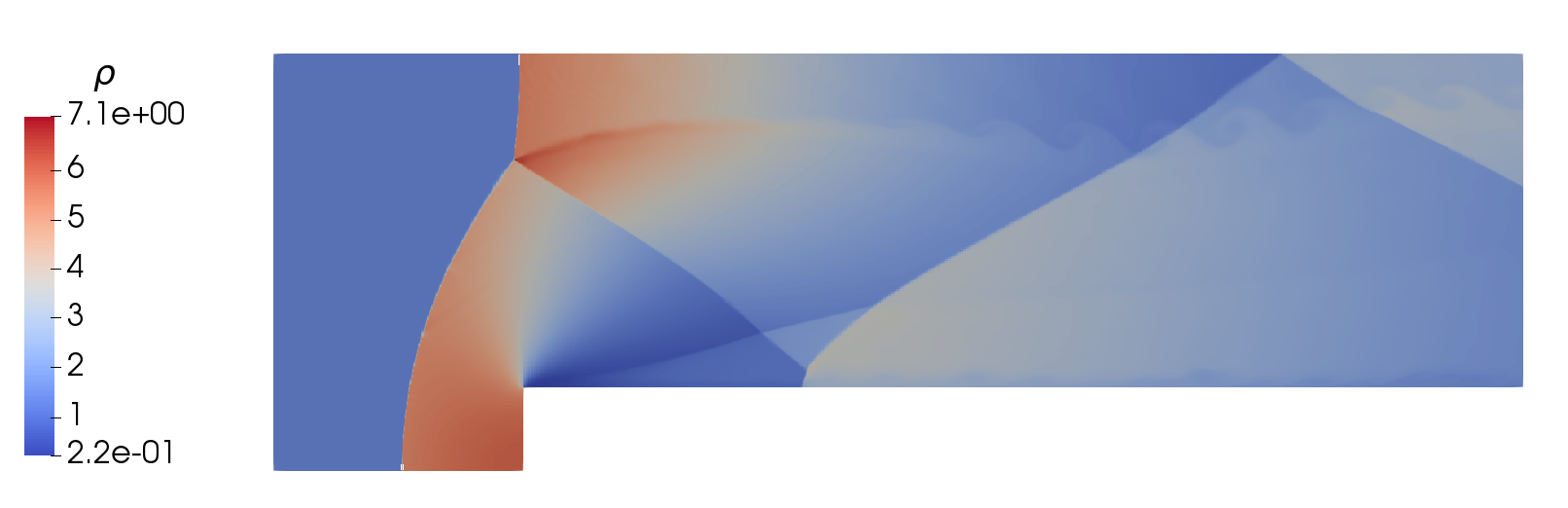}
		\includegraphics[trim={8cm 1cm 1cm 0},clip,width=\linewidth]{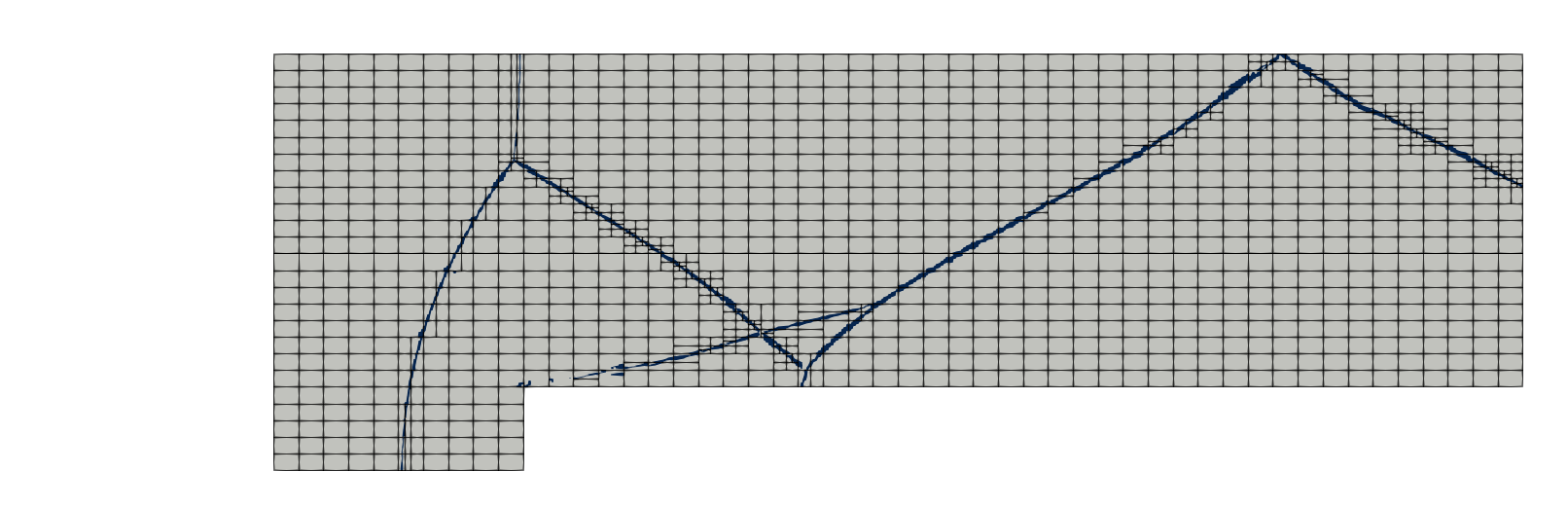}
	\end{minipage}
	\begin{minipage}[c]{0.06\linewidth}
		\begin{tikzpicture}[line cap=round,line join=round,x=1.15cm,y=1.15cm,axis/.style={->}]
		\node at (0 ,-0.4) {\includegraphics[scale=0.1]{./pictures/bar_color_red.png}};
		\node at (0.0,1.6 ) {$~~~\rho~~~$};
		\node at (0.5,1.2 ) {$7.1$};
		\node at (0.5,0.4 ) {$5.5$};
		\node at (0.5,-0.4 ) {$3.7$};
		\node at (0.5,-1.2 ) {$1.8$};
		\node at (0.5,-2.0 ) {$0.22$};
		\end{tikzpicture}
	\end{minipage}
	\caption{Computational grid (upper), density (middle) and binary shock edge map obtained with the ANNSL network with the underlying grid (bottom) of the FFS computed on a direction-wise equidistant mesh (left) and on a h-refined version of it (right).}\label{fig:ANNSL_FFS}
\end{figure}

We use the flow over a forward facing step problem on a regular Cartesian mesh defined in Sec.~\ref{subsec:ffs} to illustrate how the information about the inner-cell shock position can be exploited during the calculation to improve the accuracy of the shock capturing approach. We restrict the discussion to a three-level h-refinement strategy in which a donor cell on the original grid can be halved in any coordinate direction separately or in both at the same time, leading to either 2 or 4 new child cells replacing the original element (anisotropic or isotropic split). This process can be repeated twice, so three levels of grid elements can be created. The original regular grid is shown in the top left column in Fig.~\ref{fig:ANNSL_FFS}. We conduct all simulations with a polynomial degree of $\mathcal{N}=9$ and the ANNSL network from the previous sections and  visualize the results at $t_\text{end}=4$.

The left column of Fig.~\ref{fig:ANNSL_FFS} shows the solution on the baseline mesh, which is  direction-wise equidistant. The corresponding shock localization indicator and the resulting density field are shown below, again, the ANNSL approach leads to a sharp and consistent prediction of the front locations, confirming the findings from Sec.~\ref{subsec:annsl1}. From the binary edge map, two pieces of information can now be extracted: First, the shock width is characterized by the amount of flagged solution points / pixels within one element. Secondly, the orientation of the shock front is available. Several methods to take advantage of this information to inform an h-refinement are possible; we choose a very simple strategy for demonstration purposes here which already produces remarkable results: We define an indicator for anisotropic mesh refinement which exploits both information from the binary edge map
\begin{align}\label{eq:MeshRefIndicator}
	\mathcal{I}_\text{meshref}^{\text{dir}}=\sum_{i=0}^{\mathcal{N}}\sum_{s=1}^{r^{\text{dir}}_i-1}\left(1.7^{s-1}+1.7^s\right),
\end{align}
independently for both the $x$ and $y$-directions, with $r^{\text{dir}}_i$ being the number of solution points with class $1$ of layer $i$ in the corresponding $x$ or $y$-direction, i.e. we count the number of flagged pixels line-wise. After evaluating the indicator we apply the mesh refinement in each direction independently. Last, a constant interpolation on the new grid is done. We apply the indicator again on the refined solution, choosing the threshold values $\mathcal{I}^\text{threshold}_\text{meshref}=33$ and $\mathcal{I}^\text{threshold}_\text{meshref}=172$ for the first and the second refinement. These three levels of grid elements provide a sharp resolution of the shocks. Note that during the refinement it is taken care that only non-conforming interfaces with a 2:1 ratio occur. Details on DGSEM and the finite volume sub-cell method on meshes with non-conforming interfaces can be found in~\cite{krais2019flexi,sonntagdiss}.
Fig.~\ref{fig:ANNSL_FFS} (right column) shows the resulting mesh, the results on the refined mesh and the binary edge map, illustrating the sharper resolution of the shocks due to the refined mesh.
\begin{figure}
	\begin{minipage}[c]{0.44\linewidth}
		\begin{tikzpicture}
		\node[inner sep=0pt] at (0,6.0)
		{\includegraphics[width=\linewidth]{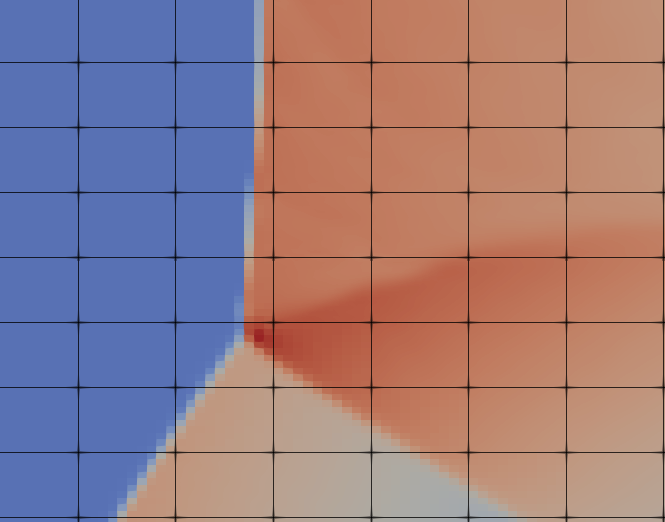}};
		\draw[thick] (-1.7,6.02) rectangle (0.7,8.1);
		\end{tikzpicture}
		\includegraphics[width=\linewidth,trim=0 0 0 0,clip]{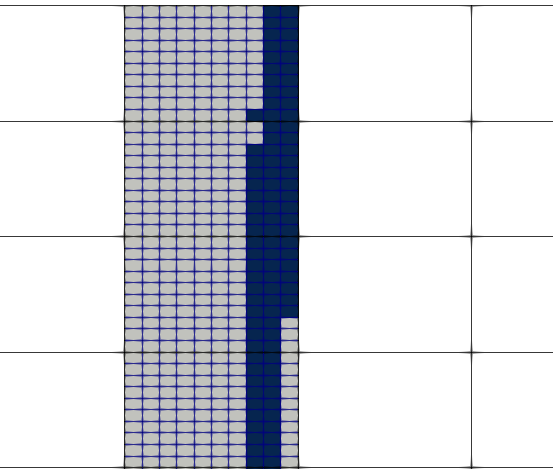}
	\end{minipage}
	\begin{minipage}[c]{0.44\linewidth}
		\begin{tikzpicture}
		\node[inner sep=0pt] at (0,6.0)
		{\includegraphics[width=\linewidth]{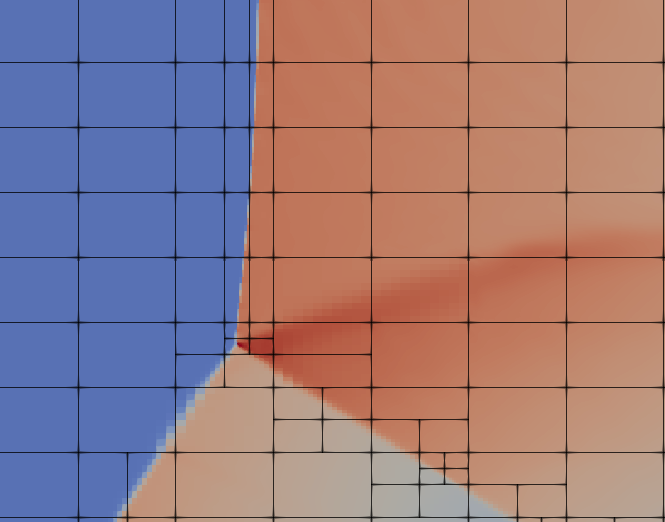}};
		\draw[thick] (-1.7,6.02) rectangle (0.7,8.1);
		\end{tikzpicture}
		\includegraphics[width=\linewidth,trim=0 0 0 0,clip]{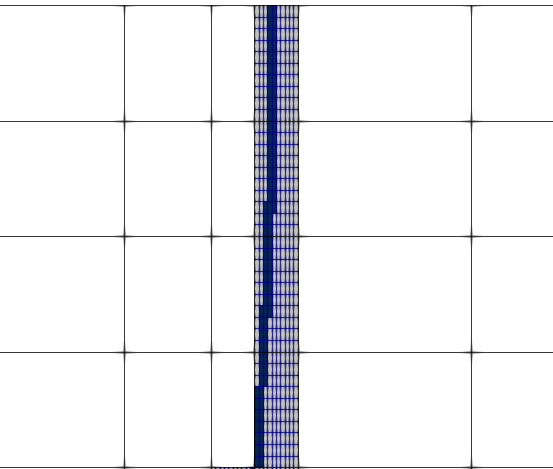}
	\end{minipage}
	\begin{minipage}[c]{0.06\linewidth}
		\begin{tikzpicture}[line cap=round,line join=round,x=1.15cm,y=1.15cm,axis/.style={->}]
		\node at (0 ,-0.4) {\includegraphics[scale=0.1]{./pictures/bar_color_red.png}};
		\node at (0.0,1.6 ) {$~~~\rho~~~$};
		\node at (0.5,1.2 ) {$7.1$};
		\node at (0.5,0.4 ) {$5.5$};
		\node at (0.5,-0.4 ) {$3.7$};
		\node at (0.5,-1.2 ) {$1.8$};
		\node at (0.5,-2.0 ) {$0.22$};
		\end{tikzpicture}
	\end{minipage}
	\caption{Top row: Zoom to density and underlying mesh at upper shock crossing point for FFS with direction-wise equidistant grid (left) and anisotropic refined mesh (right). The region shown in the bottom row is marked by the black rectangle. Bottom row: Detailed view of edge map shock predictions. The grid cells are denoted by black lines. Cells without shocks have a white background. In cells in which a shock is detected, the finite volume subgrid is shown, dark colors indicate the shock location prediction.\label{fig:ANNSL_FFS_Detail}}
\end{figure}
Fig.~\ref{fig:ANNSL_FFS_Detail} (top row) provides a zoomed-in view  on the upper crossing point of the shocks, showing both the density field and the associated mesh. It is easy to recognize how the refinement informed by ANNSL has improved the sharpness of the shock while avoiding unnecessary refinement of cells or along directions. Also, the refinement criterion proposed in Eq.~\eqref{eq:MeshRefIndicator} has guided the selection between anisotropic and isotropic splits, which additionally helps to keep the number of newly introduced cells and thus computational cost low. In the bottom row of Fig.~\ref{fig:ANNSL_FFS_Detail}, the computational grid and the binary edge map in this region are depicted. Cells of the computational grid are denoted by black lines; in elements with white backgrounds, no shocks have been detected and the DG representation of the solution is used. In cells with a light blue background, shocks are present, and the FV representation of the solution on equispaced grid cells is active. In FV cells / pixels of the binary edge map marked with a dark color, the ANNSL method has predicted the occurrence of a shock. As discussed before, the network has been designed to mark the two adjacent pixels to the shock front per direction, i.e. to provide a bounding box of the discontinuity on the scale of smallest resolved length scale within a DG element. The bottom plots in Fig.~\ref{fig:ANNSL_FFS_Detail} shows that ANNSL method successfully provides such a bounding box along the shock front, given a consistent and sharp inner-cell shock localization on both the original (left) and refined grid (right).\\
 Concluding, this example illustrates how the information obtained by the novel indicator for shock localization can be used in an anisotropic adaptive mesh refinement framework to improve the resolution of the shock. Due to the separation of shock detection and shock capturing in our approach, the indicator can be transferred to other discretization schemes and can then help to inform other shock capturing methods.

\section{\label{sec:conout}Conclusion and Outlook}
The stable approximation of shocks or very steep gradients on a finite computational grid is a challenging task due to the non-linear interactions of solution behaviour and discretization. In the most commonly used approach, parameterized a priori indicators are used to not only detect the occurrence of a shock within an element, but also to implicitly predict the stability of the discretized PDE solution at the next time step. This subtle redefinition of tasks is often expressed by the use of the term "troubled cell" indicators. The interactions of solution, discretization and indicator parameters introduces the need for parameter tuning and user experience, and considerable time may be spent to find an optimal indicator setting. \\
In this work, we have proposed a novel indicator based on edge detection methods from machine learning and computer vision. The indicator is developed, tested and applied in a high order discontinuous Galerkin setting, but it can easily be transferred to related discretization schemes. In order to make this indicator (at least more) universal and parameter-independent, we have returned to the original concept of separating the tasks of detection a shock and of its stable and accurate numerical approximation. The latter one is based on our hybrid finite volume / discontinuous Galerkin scheme and is treated as a black box for the purpose of this investigation. The task of detecting and localizing a shock is performed by neural networks designed for edge detection and works on the basis of considering the inner element solution akin to an image and the nodal degrees of freedom as its pixels. The neural networks are trained on analytical smooth and non-smooth functions to provide a binary edge map in each element that bounds the shock fronts in between two solution points. This information can then be used to flag an element containing a shock in the sense of a detector or to guide an inner-cell shock capturing method.  \\
This feature is particularly attractive for high order schemes in non-trivial geometries and for complex flows, where the location of the shocks might not be obvious a priori (so the grid can be generated accordingly) or where shock movement occurs. Due to the typically large grid cells in use with these schemes, locating the shock front accurately and reliably within an element is useful to guide h/p/r-refinement methods, artificial viscosity based approaches, enriched basis functions or other methods to stably approximate the solution. \\
We see our work here as a step towards this goal. Our shock detection indicator has shown to be as robust and accurate as two established shock indicators, without the need for user intervention. It has also proven to be robust against changes in spatial resolution. The novel shock localization algorithm is capable of reliably bounding the shock front within a high order element, even on coarse grids. Its location prediction behaves consistently when refining the grid, opening the way for targeted, multilevel grid refinement. We have demonstrated how to use the information provided by the indicator to improve solution quality considerably by local h-refinement. \\
Based on the current work, a number of possible routes can be investigated in the future: The extension to 3D is a first natural step. Here, a tensorproduct expansion of the training data is an obvious choice. Incorporating more physical information into the training, i.e. through training on e.g. the density and pressure fields can extend the methods towards identifying other flow features. This might become particularly interesting for the notoriously difficult shock-turbulence interaction problems. Also, beyond the binary edge map, predicting a shock strength or other properties can be useful for informing the shock capturing method. Improving network performance and its integration in existing solvers are important technical issues to be addressed.\\

It can be argued that the costs of computing the neural network based indicator are certainly higher than those of established formulations, so the natural question is to ask if it is worth it. While we see our work as a proof of concept and do not directly tackle this issue here, answering this question in the future will depend on more than just the cost per evaluation: For one, the presented indicators can help to reduce the human and computational effort of parameter tuning. Even more significant is the potential for saving effort and improving solution accuracy through exploitation of the precise shock position, in particular in 3D for high order discretizations. Also to be taken into account are moves towards integration of specialized hardware for machine learning applications in HPC systems.\\

\appendix
\section{Generation of 1D Training Data\label{app:data}}
%
 In the following, we detail the generation of training samples for the different analytical functions with $\eps=0.1$, chosen from 7 families:

 \begin{itemize}
 	\item No.~1: Linear functions, initialized on three meshes $n_e=1,10,20$. The probability of choosing each of the meshes is $p_e=0.5,0.3,0.2$, respectively.
 	\item No.~2: Superposition of oscillations with different amplitudes and frequencies, initialized on meshes with $n_e=1,10,20$ and probabilities of $p_e=0.3,0.4,0.3$. Note that $f_\text{Nyquist}=\mathcal{N}/2$ approximates the maximum possible frequency which can be resolved by the Lagrange polynomials.
 	\item No.~3: Exponential functions in both directions, initialized on $n_e=10,20$ with probabilities $p_e=0.6,0.4$. No function evaluation with $n_e=1$ is done in order not to indicate too steep gradients with class $0$.
 	\item No.~4: Four constant values separated with straight and curved lines, initialized on $n_e=1$. The assignment to class $1$ is done if the maximum jump height fulfills the condition $|\text{max}(|u_i|)-\text{min}(u_i)|>\eps~\text{max}(|u_i|)$ and $\text{max}(|u_i|)>0.01$ holds. Degrees of freedom neighboring or being directly at the discontinuity are labeled as $1$. The parameter $d$ is chosen to give a ratio of straight to curved separation lines of 7:3.
 	\item No.~5: Magnitude function with linear gradients connected by kinks, initialized on $n_e=1$. The assignment to class $1$ is done if the condition $a>\eps~\text{max}(|u_i|)$ for the gradient $a$ and $\text{max}(|u_i|)>0.01$ for the maximum value holds and is set for degrees of freedom neighboring or being directly at the discontinuity.
 	\item No.~6: Linear or constant states connected by kinks, initialized on $n_e=1$.  The assignment to class $1$ is done if the maximum gradient fulfills $\text{max}(b_1,b_2)>\eps~\text{max}(|u_i|)$ and $\text{max}(|u_i|)>0.01$ holds and is set for degrees of freedom neighboring or being directly at the discontinuity.
 	\item No.~7: Decaying high frequency oscillations in analogy to Gibb's instability, initialized on $n_e=1$. The assignment to class $1$ is done if the maximum difference in the element fulfills the condition $|\text{max}(|u_i|)-\text{min}(u_i)|>0.2~\text{max}(|u_i|)$ and $\text{max}(|u_i|)>0.01$ holds. Then degrees of freedom having the highest absolute value in the element are labeled as $1$.
 \end{itemize}
Tbl.~\ref{Tbl:exactfunc_2D} summarizes the function families and their parameters. After the generation of the training samples, the inputs are normalized to enhance the convergence of the training process. We first shift the data
\begin{align*}
u_i = \begin{cases}
u_i - \text{min}_x u_i & \text{if} \ \text{min}_x u_i < 0, \\
u_i                    & \text{otherwise}, \end{cases}
\end{align*}
and then scale them to the interval $\mathcal{I}=[0,1]$
\begin{align*}
u_i = \begin{cases}
\frac{u_i}{\text{max}_x \vert{u_i\vert}} & \text{if} \ \text{max}_x \vert{u_i\vert} > 1, \\
u_i                    & \text{otherwise}. \end{cases}
\end{align*}
This enables to train with data in the positive interval $\mathcal{I}=[0,1]$ while retaining the information about relative differences in the solution. To obtain training samples from the analytical functions, they are initialized on Cartesian meshes covering the interval $[-1,1]^2$ with $n_e^2$ elements. To account for under resolved phenomena, we use a polynomial degree in each grid cell of $2\mathcal{N}$ for the evaluation of the analytical functions. Next, the analytical functions given in the Lagrange basis representation of degree $2\mathcal{N}$ are projected onto the Lagrange basis representation of degree $\mathcal{N}$, which is used for the simulation with DGSEM.
Summing up, the input image $\vec{X}$ is build up with the normalized nodal values of the Lagrange polynomials of degree $\mathcal{N}$ resulting in $(\mathcal{N}+1)^2$ values. The true output $\vec{Y}$ is the binary edge map.
\begin{table}[!htb]
	\centering
	\def\arraystretch{1.2}
	\begin{tabularx}{1.0\textwidth}{|l|X|X|r|}
		\hline
		\# & $u(x,y)$          & Parameters & Class  \\ \hline\hline
		1 & $ax + by$           & $a,b \in \mathcal{N}[0,0.2]$ & 0 \\ \hline
		2 & $\sum_{k=1}^{N_f} a_k sin(k\pi x)+ \newline b_k cos(k\pi y) + c$ & $a_k,b_k \in \mathcal{U}[-0.5,0.5],\newline c \in \mathcal{U}[0,1], \newline N_f=1,..,f_{\text{Nyquist}}$ & 0 \\ \hline
		3 & $exp(a_1[(x-a_2)^2+(y-a_3)^2])+ \newline exp(a_4[(x-a_5)^2+(y-a_6)^2])$ & $a_i \in \mathcal{U}[-1,1]$ & 0 \\ \hline
		4 & 4 values $u_i$ in 4 sections defined by\newline
		$y-y_0 = m(x-x_0)^d$ \newline$y-y_0 = -1/m(x-x_0)^d$ & $u_i \in \mathcal{U}[0,1],\ m \in \mathcal{U}[0,10]$,\newline $x_0,y_0 \in \mathcal{U}[-1,1], \newline d = \{1,2\}$ & 0,1 \\ \hline
		5 & $a \vert{(y-y_0)-m(x-x_0)\vert} + c$ & $a \in \mathcal{N}[0,0.4],\ m \in \mathcal{U}[-2,2]$, $x_0,y_0 \in \mathcal{U}[-1,1], \ c \in \mathcal{U}[0,1]$ & 0,1 \\ \hline
		6 & \textbf{if} $a_1>0$ \textbf{then} \newline $a_2 \text{max}(0, b_1(x - x_0)) \newline + a_3 \text{max}(0,b_2(y-y_0)) + c$\newline \textbf{else} \newline$a_2 \text{max}(0, b_1 (x - x_0)+ \newline b_2 (y-y_0)) + c$ & $a_i = \{-1,1\},\newline b_i \in \mathcal{N}[0,0.6] \newline x_0, y_0 \in \mathcal{U}[-0.6,0.6],\newline c \in \mathcal{U}[0,1]$ & 0,1 \\ \hline
		7 & $a_1 \text{sin}(f_\text{Nyquist} \pi (x-x_0)) \newline \text{exp}(a_3(x-x_0)) \newline+ a_2 \text{cos}(f_\text{Nyquist} \pi (y-y_0))\newline \text{exp}(a_3(y-y_0)) + c$ & $c \in \mathcal{U}[0,1], \newline a_1, a_2 \in \mathcal{N}[0,0.4],\newline x_0, y_0 \in \mathcal{U}[-1,1], \newline a_3 \in \mathcal{U}[-2,2]$ & 0,1 \\ \hline
	\end{tabularx}
	\caption{Functions used to generate the training and validation sets on the interval $[-1,1]^2$. Note that $\mathcal{U}[\mathcal{I}]$ denotes the uniform distribution on the interval $\mathcal{I}$ and $\mathcal{N}[\mu,\sqrt{\sigma^2}]$ denotes the normal distribution with mean $\mu$ and variance $\sigma^2$.}
	\label{Tbl:exactfunc_2D}
\end{table}

\newpage
\bibliographystyle{acm}
\bibliography{references}

\begin{thebibliography}{10}

\bibitem{BALSARA2007586}
{\sc Balsara, D.~S., Altmann, C., Munz, C.-D., and Dumbser, M.}
\newblock {A sub-cell based indicator for troubled zones in RKDG schemes and a
  novel class of hybrid RKDG+HWENO schemes}.
\newblock {\em Journal of Computational Physics 226}, 1 (2007), 586--620.

\bibitem{barron1993universal}
{\sc Barron, A.~R.}
\newblock Universal approximation bounds for superpositions of a sigmoidal
  function.
\newblock {\em IEEE Transactions on Information theory 39}, 3 (1993), 930--945.

\bibitem{BARTER20101810}
{\sc Barter, G.~E., and Darmofal, D.~L.}
\newblock {Shock capturing with PDE-based artificial viscosity for DGFEM: Part
  I. Formulation}.
\newblock {\em Journal of Computational Physics 229}, 5 (2010), 1810--1827.

\bibitem{beck2019deep}
{\sc Beck, A., Flad, D., and Munz, C.-D.}
\newblock Deep neural networks for data-driven {LES} closure models.
\newblock {\em Journal of Computational Physics 398\/} (2019), 108910.

\bibitem{beck2014high}
{\sc Beck, A.~D., Bolemann, T., Flad, D., Frank, H., Gassner, G.~J.,
  Hindenlang, F., and Munz, C.-D.}
\newblock High-order discontinuous {G}alerkin spectral element methods for
  transitional and turbulent flow simulations.
\newblock {\em International Journal for Numerical Methods in Fluids 76}, 8
  (2014), 522--548.

\bibitem{BURBEAU2001111}
{\sc Burbeau, A., Sagaut, P., and Bruneau, C.-H.}
\newblock {A Problem-Independent Limiter for High-Order Runge–Kutta
  Discontinuous Galerkin Methods}.
\newblock {\em Journal of Computational Physics 169}, 1 (2001), 111--150.

\bibitem{10.2307/2008501}
{\sc Cockburn, B., Hou, S., and Shu, C.-W.}
\newblock The {R}unge-{K}utta local projection discontinuous {G}alerkin finite
  element method for conservation laws. {IV}: The multidimensional case.
\newblock {\em Mathematics of Computation 54}, 190 (1990), 545--581.

\bibitem{courant1999supersonic}
{\sc Courant, R., and Friedrichs, K.~O.}
\newblock {\em Supersonic flow and shock waves}, vol.~21.
\newblock Springer Science \& Business Media, 1999.

\bibitem{cybenko1989approximation}
{\sc Cybenko, G.}
\newblock Approximation by superpositions of a sigmoidal function.
\newblock {\em Mathematics of control, signals and systems 2}, 4 (1989),
  303--314.

\bibitem{dolejvsi2003}
{\sc Dolej{\v{s}}{\'\i}, V., Feistauer, M., and Schwab, C.}
\newblock On some aspects of the discontinuous {G}alerkin finite element method
  for conservation laws.
\newblock {\em Mathematics and Computers in Simulation 61}, 3-6 (2003),
  333--346.

\bibitem{ducros1999}
{\sc Ducros, F., Ferrand, V., Nicoud, F., Weber, C., Darracq, D., Gacherieu,
  C., and Poinsot, T.}
\newblock Large-eddy simulation of the shock/turbulence interaction.
\newblock {\em Journal of Computational Physics 152}, 2 (1999), 517--549.

\bibitem{dumbser2014posteriori}
{\sc Dumbser, M., Zanotti, O., Loub{\`e}re, R., and Diot, S.}
\newblock A posteriori subcell limiting of the discontinuous {G}alerkin finite
  element method for hyperbolic conservation laws.
\newblock {\em Journal of Computational Physics 278\/} (2014), 47--75.

\bibitem{FEISTAUER20101612}
{\sc Feistauer, M., Ku{\v{c}}era, V., and Prokopov{\'{a}}, J.}
\newblock {Discontinuous {G}alerkin solution of compressible flow in
  time-dependent domains}.
\newblock {\em Mathematics and Computers in Simulation 80}, 8 (2010),
  1612--1623.

\bibitem{gottlieb1997gibbs}
{\sc Gottlieb, D., and Shu, C.-W.}
\newblock On the {G}ibbs phenomenon and its resolution.
\newblock {\em SIAM review 39}, 4 (1997), 644--668.

\bibitem{hanveiga:hal-01856358}
{\sc Han~Veiga, M.~M., and Abgrall, R.}
\newblock {Towards a general stabilisation method for conservation laws using a
  multilayer Perceptron neural network: 1D scalar and system of equations}.
\newblock In {\em {ECCM - ECFD 2018 6th European Conference on Computational
  Mechanics (Solids, Structures and Coupled Problems) 7th European Conference
  on Computational Fluid Dynamics}\/} (Glasgow, United Kingdom, June 2018).

\bibitem{HARTMANN2002508}
{\sc Hartmann, R., and Houston, P.}
\newblock {Adaptive Discontinuous {G}alerkin Finite Element Methods for the
  Compressible Euler Equations}.
\newblock {\em Journal of Computational Physics 183}, 2 (2002), 508--532.

\bibitem{haykin1994neural}
{\sc Haykin, S.}
\newblock {\em Neural networks: a comprehensive foundation}.
\newblock Prentice Hall PTR, 1994.

\bibitem{hindenlang2012explicit}
{\sc Hindenlang, F., Gassner, G.~J., Altmann, C., Beck, A., Staudenmaier, M.,
  and Munz, C.-D.}
\newblock Explicit discontinuous {G}alerkin methods for unsteady problems.
\newblock {\em Computers \& Fluids 61\/} (2012), 86--93.

\bibitem{HORNIK1991251}
{\sc Hornik, K.}
\newblock {Approximation capabilities of multilayer feedforward networks}.
\newblock {\em Neural Networks 4}, 2 (1991), 251--257.

\bibitem{huerta2012}
{\sc Huerta, A., Casoni, E., and Peraire, J.}
\newblock A simple shock-capturing technique for high-order discontinuous
  {Galerkin} methods.
\newblock {\em International journal for numerical methods in fluids 69}, 10
  (2012), 1614--1632.

\bibitem{ioffe2015batch}
{\sc Ioffe, S., and Szegedy, C.}
\newblock Batch normalization: Accelerating deep network training by reducing
  internal covariate shift.
\newblock {\em arXiv preprint arXiv:1502.03167\/} (2015).

\bibitem{jameson1981}
{\sc Jameson, A., Schmidt, W., and Turkel, E.}
\newblock Numerical solution of the {Euler} equations by finite volume methods
  using {Runge} {Kutta} time stepping schemes.
\newblock In {\em 14th fluid and plasma dynamics conference\/} (1981), p.~1259.

\bibitem{kennedy2000low}
{\sc Kennedy, C.~A., Carpenter, M.~H., and Lewis, R.~M.}
\newblock Low-storage, explicit {R}unge--{K}utta schemes for the compressible
  {N}avier--{S}tokes equations.
\newblock {\em Applied numerical mathematics 35}, 3 (2000), 177--219.

\bibitem{kingma2014}
{\sc Kingma, D.~P., and Ba, J.}
\newblock Adam: A method for stochastic optimization.
\newblock {\em arXiv preprint arXiv:1412.6980\/} (2014).

\bibitem{klockner2011viscous}
{\sc Kl{\"o}ckner, A., Warburton, T., and Hesthaven, J.~S.}
\newblock Viscous shock capturing in a time-explicit discontinuous {G}alerkin
  method.
\newblock {\em Mathematical Modelling of Natural Phenomena 6}, 3 (2011),
  57--83.

\bibitem{kopriva2009}
{\sc Kopriva, D.~A.}
\newblock {\em Implementing spectral methods for partial differential
  equations: {Algorithms} for scientists and engineers}.
\newblock Springer Science \& Business Media, 2009.

\bibitem{krais2019flexi}
{\sc Krais, N., Beck, A., Bolemann, T., Frank, H., Flad, D., Gassner, G.,
  Hindenlang, F., Hoffmann, M., Kuhn, T., Sonntag, M., et~al.}
\newblock Flexi: A high order discontinuous galerkin framework for
  hyperbolic-parabolic conservation laws.
\newblock {\em arXiv preprint arXiv:1910.02858\/} (2019).

\bibitem{krivodonova2007limiters}
{\sc Krivodonova, L.}
\newblock Limiters for high-order discontinuous {Galerkin} methods.
\newblock {\em Journal of Computational Physics 226}, 1 (2007), 879--896.

\bibitem{krizhevsky2012imagenet}
{\sc Krizhevsky, A., Sutskever, I., and Hinton, G.~E.}
\newblock Imagenet classification with deep convolutional neural networks.
\newblock In {\em Advances in neural information processing systems\/} (2012),
  pp.~1097--1105.

\bibitem{lecun1995convolutional}
{\sc LeCun, Y., Bengio, Y., et~al.}
\newblock Convolutional networks for images, speech, and time series.
\newblock {\em The handbook of brain theory and neural networks 3361}, 10
  (1995), 1995.

\bibitem{lecun2015}
{\sc LeCun, Y., Bengio, Y., and Hinton, G.}
\newblock Deep learning.
\newblock {\em nature 521}, 7553 (2015), 436--444.

\bibitem{lecun1990handwritten}
{\sc LeCun, Y., Boser, B.~E., Denker, J.~S., Henderson, D., Howard, R.~E.,
  Hubbard, W.~E., and Jackel, L.~D.}
\newblock Handwritten digit recognition with a back-propagation network.
\newblock In {\em Advances in neural information processing systems\/} (1990),
  pp.~396--404.

\bibitem{lecun1998gradient}
{\sc LeCun, Y., Bottou, L., Bengio, Y., Haffner, P., et~al.}
\newblock Gradient-based learning applied to document recognition.
\newblock {\em Proceedings of the IEEE 86}, 11 (1998), 2278--2324.

\bibitem{lee2015deeply}
{\sc Lee, C.-Y., Xie, S., Gallagher, P., Zhang, Z., and Tu, Z.}
\newblock Deeply-supervised nets.
\newblock In {\em Artificial intelligence and statistics\/} (2015),
  pp.~562--570.

\bibitem{leveque1998nonlinear}
{\sc LeVeque, R.~J.}
\newblock Nonlinear conservation laws and finite volume methods.
\newblock In {\em Computational methods for astrophysical fluid flow}.
  Springer, 1998, pp.~1--159.

\bibitem{leveque2002finite}
{\sc LeVeque, R.~J., et~al.}
\newblock {\em Finite volume methods for hyperbolic problems}, vol.~31.
\newblock Cambridge university press, 2002.

\bibitem{liou1995image}
{\sc Liou, S.-P., Singh, A., Mehlig, S., Edwards, D., and Davis, R.}
\newblock An image analysis based approach to shock identification in cfd.
\newblock In {\em 33rd Aerospace Sciences Meeting and Exhibit\/} (1995),
  p.~117.

\bibitem{liu2017}
{\sc Liu, Y., Cheng, M.-M., Hu, X., Wang, K., and Bai, X.}
\newblock Richer convolutional features for edge detection.
\newblock In {\em Proceedings of the IEEE conference on computer vision and
  pattern recognition\/} (2017), pp.~3000--3009.

\bibitem{LIU20191}
{\sc Liu, Y., Lu, Y., Wang, Y., Sun, D., Deng, L., Wang, F., and Lei, Y.}
\newblock A cnn-based shock detection method in flow visualization.
\newblock {\em Computers \& Fluids 184\/} (2019), 1 -- 9.

\bibitem{lovely1999shock}
{\sc Lovely, D., and Haimes, R.}
\newblock Shock detection from computational fluid dynamics results.
\newblock In {\em 14th Computational Fluid Dynamics Conference\/} (1999),
  p.~3285.

\bibitem{NIPS2017_7203}
{\sc Lu, Z., Pu, H., Wang, F., Hu, Z., and Wang, L.}
\newblock {The Expressive Power of Neural Networks: A View from the Width}.
\newblock In {\em Advances in Neural Information Processing Systems 30},
  I.~Guyon, U.~V. Luxburg, S.~Bengio, H.~Wallach, R.~Fergus, S.~Vishwanathan,
  and R.~Garnett, Eds. Curran Associates, Inc., 2017, pp.~6231--6239.

\bibitem{lv2016entropy}
{\sc Lv, Y., See, Y.~C., and Ihme, M.}
\newblock An entropy-residual shock detector for solving conservation laws
  using high-order discontinuous {G}alerkin methods.
\newblock {\em Journal of Computational Physics 322\/} (2016), 448--472.

\bibitem{maas2013}
{\sc Maas, A.~L., Hannun, A.~Y., and Ng, A.~Y.}
\newblock Rectifier nonlinearities improve neural network acoustic models.
\newblock In {\em Proc. icml\/} (2013), vol.~30, p.~3.

\bibitem{10.1007/978-3-319-61358-1_16}
{\sc Monfort, M., Luciani, T., Komperda, J., Ziebart, B., Mashayek, F., and
  Marai, G.~E.}
\newblock {A Deep Learning Approach to Identifying Shock Locations in Turbulent
  Combustion Tensor Fields}.
\newblock In {\em Modeling, Analysis, and Visualization of Anisotropy\/} (Cham,
  2017), T.~Schultz, E.~{\"{O}}zarslan, and I.~Hotz, Eds., Springer
  International Publishing, pp.~375--392.

\bibitem{morgan2020machine}
{\sc Morgan, N.~R., Tokareva, S., Liu, X., and Morgan, A.}
\newblock A machine learning approach for detecting shocks with high-order
  hydrodynamic methods.
\newblock In {\em AIAA Scitech 2020 Forum\/} (2020), p.~2024.

\bibitem{paciorri2009shock}
{\sc Paciorri, R., and Bonfiglioli, A.}
\newblock A shock-fitting technique for 2d unstructured grids.
\newblock {\em Computers \& Fluids 38}, 3 (2009), 715--726.

\bibitem{pagendarm1993algorithm}
{\sc Pagendarm, H.-G., and Seitz, B.}
\newblock An algorithm for detection and visualization of discontinuities in
  scientific data fields applied to flow data with shock waves.
\newblock {\em Scientific Visualization: Advanced Software Techniques\/}
  (1993), 161--177.

\bibitem{persson2006}
{\sc Persson, P.-O., and Peraire, J.}
\newblock Sub-cell shock capturing for discontinuous {G}alerkin methods.
\newblock In {\em 44th AIAA Aerospace Sciences Meeting and Exhibit\/} (2006),
  p.~112.

\bibitem{pirozzoli2011numerical}
{\sc Pirozzoli, S.}
\newblock Numerical methods for high-speed flows.
\newblock {\em Annual review of fluid mechanics 43\/} (2011), 163--194.

\bibitem{premasuthan2014computation}
{\sc Premasuthan, S., Liang, C., and Jameson, A.}
\newblock Computation of flows with shocks using the spectral difference method
  with artificial viscosity, i: basic formulation and application.
\newblock {\em Computers \& Fluids 98\/} (2014), 111--121.

\bibitem{doi:10.1137/04061372X}
{\sc Qiu, J., and Shu, C.-W.}
\newblock {A Comparison of Troubled-Cell Indicators for Runge--Kutta
  Discontinuous Galerkin Methods Using Weighted Essentially Nonoscillatory
  Limiters}.
\newblock {\em SIAM Journal on Scientific Computing 27}, 3 (2005), 995--1013.

\bibitem{QIU2005642}
{\sc Qiu, J., and Shu, C.-W.}
\newblock {Hermite WENO schemes and their application as limiters for
  Runge–Kutta discontinuous Galerkin method II: Two dimensional case}.
\newblock {\em Computers {\&} Fluids 34}, 6 (2005), 642--663.

\bibitem{rabault2019artificial}
{\sc Rabault, J., Kuchta, M., Jensen, A., R{\'e}glade, U., and Cerardi, N.}
\newblock Artificial neural networks trained through deep reinforcement
  learning discover control strategies for active flow control.
\newblock {\em Journal of Fluid Mechanics 865\/} (2019), 281--302.

\bibitem{ray2018}
{\sc Ray, D., and Hesthaven, J.~S.}
\newblock An artificial neural network as a troubled-cell indicator.
\newblock {\em Journal of computational physics 367\/} (2018), 166--191.

\bibitem{ray2019}
{\sc Ray, D., and Hesthaven, J.~S.}
\newblock Detecting troubled-cells on two-dimensional unstructured grids using
  a neural network.
\newblock {\em Journal of Computational Physics 397\/} (2019), 108845.

\bibitem{riemann1860fortpflanzung}
{\sc Riemann, B., et~al.}
\newblock {\em {\"U}ber die Fortpflanzung ebener Luftwellen von endlicher
  Schwingungsweite}.
\newblock Verlag der Dieterichschen Buchhandlung, 1860.

\bibitem{rumelhart1986learning}
{\sc Rumelhart, D.~E., Hinton, G.~E., and Williams, R.~J.}
\newblock Learning representations by back-propagating errors.
\newblock {\em nature 323}, 6088 (1986), 533.

\bibitem{rusanov1973processing}
{\sc Rusanov, V.}
\newblock Processing and analysis of computation results for multidimensional
  problems of aerohydrodynamics.
\newblock In {\em Proceedings of the Third International Conference on
  Numerical Methods in Fluid Mechanics\/} (1973), Springer, pp.~154--162.

\bibitem{schmidhuber2015deep}
{\sc Schmidhuber, J.}
\newblock Deep learning in neural networks: An overview.
\newblock {\em Neural networks 61\/} (2015), 85--117.

\bibitem{schulz1993numerical}
{\sc Schulz-Rinne, C.~W., Collins, J.~P., and Glaz, H.~M.}
\newblock Numerical solution of the {R}iemann problem for two-dimensional gas
  dynamics.
\newblock {\em SIAM Journal on Scientific Computing 14}, 6 (1993), 1394--1414.

\bibitem{sheshadri2014shock}
{\sc Sheshadri, A., and Jameson, A.}
\newblock Shock detection and capturing methods for high order
  discontinuous-{G}alerkin finite element methods.
\newblock In {\em 32nd AIAA Applied Aerodynamics Conference\/} (2014), p.~2688.

\bibitem{shu2003high}
{\sc Shu, C.-W.}
\newblock High-order finite difference and finite volume {WENO }schemes and
  discontinuous {G}alerkin methods for cfd.
\newblock {\em International Journal of Computational Fluid Dynamics 17}, 2
  (2003), 107--118.

\bibitem{shu1988efficient}
{\sc Shu, C.-W., and Osher, S.}
\newblock Efficient implementation of essentially non-oscillatory
  shock-capturing schemes.
\newblock {\em Journal of computational physics 77}, 2 (1988), 439--471.

\bibitem{sonntagdiss}
{\sc Sonntag, M.}
\newblock {\em Shape derivatives and shock capturing for the {N}avier-{S}tokes
  equations in discontinuous {G}alerkin methods}.
\newblock Dissertation, University of Stuttgart, 2017.

\bibitem{sonntag2017efficient}
{\sc Sonntag, M., and Munz, C.-D.}
\newblock Efficient parallelization of a shock capturing for discontinuous
  {G}alerkin methods using finite volume sub-cells.
\newblock {\em Journal of Scientific Computing 70}, 3 (2017), 1262--1289.

\bibitem{nnzoo}
{\sc Van~Veen, F., and Leijnen, S.}
\newblock The neural network zoo.
\newblock http://www.asimovinstitute.org/neural-network-zoo.

\bibitem{vorozhtsov1987shock}
{\sc Vorozhtsov, E.}
\newblock On shock localization by digital image processing techniques.
\newblock {\em Computers \& Fluids 15}, 1 (1987), 13--45.

\bibitem{wang2013high}
{\sc Wang, Z.~J., Fidkowski, K., Abgrall, R., Bassi, F., Caraeni, D., Cary, A.,
  Deconinck, H., Hartmann, R., Hillewaert, K., Huynh, H.~T., et~al.}
\newblock High-order cfd methods: current status and perspective.
\newblock {\em International Journal for Numerical Methods in Fluids 72}, 8
  (2013), 811--845.

\bibitem{werbos}
{\sc Werbos, P.~J.}
\newblock Backpropagation through time: what it does and how to do it.
\newblock {\em Proceedings of the IEEE 78}, 10 (Oct 1990), 1550--1560.

\bibitem{woodward1984numerical}
{\sc Woodward, P., and Colella, P.}
\newblock The numerical simulation of two-dimensional fluid flow with strong
  shocks.
\newblock {\em Journal of computational physics 54}, 1 (1984), 115--173.

\bibitem{WU2013501}
{\sc Wu, Z., Xu, Y., Wang, W., and Hu, R.}
\newblock Review of shock wave detection method in {CFD} post-processing.
\newblock {\em Chinese Journal of Aeronautics 26}, 3 (2013), 501 -- 513.

\bibitem{xie2015}
{\sc Xie, S., and Tu, Z.}
\newblock Holistically-nested edge detection.
\newblock In {\em Proceedings of the IEEE international conference on computer
  vision\/} (2015), pp.~1395--1403.

\bibitem{yang2009parameter}
{\sc Yang, M., and Wang, Z.-J.}
\newblock A parameter-free generalized moment limiter for high-order methods on
  unstructured grids.
\newblock In {\em 47th AIAA Aerospace Sciences Meeting Including The New
  Horizons Forum and Aerospace Exposition\/} (2009), p.~605.

\bibitem{yu2018data}
{\sc Yu, J., Hesthaven, J.~S., and Yan, C.}
\newblock A data-driven shock capturing approach for discontinuous {G}alerkin
  methods.
\newblock Tech. rep., 2018.

\bibitem{ZHONG2013397}
{\sc Zhong, X., and Shu, C.-W.}
\newblock {A simple weighted essentially nonoscillatory limiter for
  Runge–Kutta discontinuous {G}alerkin methods}.
\newblock {\em Journal of Computational Physics 232}, 1 (2013), 397--415.

\end{thebibliography}

\end{document}